\documentclass[reqno,final]{amsart}

\usepackage{fancyhdr}
\usepackage{color}
\usepackage{hyperref}
\usepackage{graphicx}

\usepackage[english]{babel}
\usepackage[utf8]{inputenc}
\usepackage[T1]{fontenc}
\usepackage{dsfont}
\usepackage{latexsym}
\usepackage{marvosym}

\usepackage{enumerate}
\usepackage{bbold}
\usepackage{verbatim}
\usepackage{cprotect}
\usepackage[vlined]{algorithm2e}
\SetAlgoCaptionSeparator{:}
\renewcommand\AlCapFnt{\normalfont}
\usepackage{tikz}

\usepackage{amssymb}
\usepackage{amsmath}
\usepackage{amsfonts}
\usepackage{enumitem}
\usepackage{titletoc}
\usepackage{caption}
\usepackage{subfig}
\usepackage{listings}
\usepackage{tikz}
\usetikzlibrary{decorations.text}

\definecolor{pcolor}{rgb}{0.8,0.2,0.3}
\hypersetup{
	breaklinks,
	colorlinks=true,
	linkcolor=pcolor,
	urlcolor=pcolor,
	citecolor=pcolor}
\renewcommand{\headheight}{11pt}
\theoremstyle{plain}
\newtheorem{theorem}{Theorem}[section]                                          
\newtheorem{proposition}[theorem]{Proposition}                          
\newtheorem{lemma}[theorem]{Lemma}
\newtheorem{corollary}[theorem]{Corollary}

\theoremstyle{definition}

\theoremstyle{remark}
\newtheorem{remark}[theorem]{Remark}

\makeatletter \@addtoreset{equation}{section} \makeatother
\renewcommand\theequation{\thesection.\arabic{equation}}
\renewcommand\thefigure{\thesection.\arabic{figure}}
\renewcommand\thetable{\thesection.\arabic{table}}

\makeatletter
\newcommand*{\centerfloat}{%
	\parindent \z@
	\leftskip \z@ \@plus 1fil \@minus \textwidth
	\rightskip\leftskip
	\parfillskip \z@skip}
\makeatother

\newcommand{\N}{\mathbf{N}} 
\newcommand{\E}{\mathbf{E}}
\newcommand{\R}{\mathbf{R}}
\newcommand{\prob}{\mathbf{P}}
\newcommand{\ee}{\mathbf{e}}
\newcommand{\e}{\epsilon}
\newcommand{\dd}{{\rm d}}
\newcommand{\Var}{\mathbf{V}\text{\normalfont ar}}

\newcommand{\iF}{F^{-1}_{\Lambda_{\infty}}}
\renewcommand{\dd}{{\rm d}}
\newcommand{\norm}[1]{|#1 |}

\DeclareMathOperator*{\argmin}{arg\,min}

\newcommand{\purcontour}{
\begin{tikzpicture}[xscale=0.4,yscale=0.5]

\tikzstyle{fleche}=[-,thin]

\tikzstyle{noeud}=[fill=black,circle,draw,scale=0.4]
\tikzstyle{noeudtext}=[]
\tikzstyle{phantom}=[-,scale=0.05]

\tikzstyle{pointt}=[dashed,thin,draw=black!70!white]
\tikzstyle{fleched}=[dashed,thin,draw=black]
\tikzstyle{fleche}=[-,draw,thick]
\tikzstyle{flech}=[-,draw]

\tikzstyle{arroww}=[->,draw,thick]

\def\DistanceInterNiveaux{1}
\def\DistanceInterFeuilles{1}

\def\NiveauA{(1)*\DistanceInterNiveaux}
\def\NiveauB{(2)*\DistanceInterNiveaux}
\def\NiveauC{(3*\DistanceInterNiveaux}
\def\NiveauD{(4*\DistanceInterNiveaux}
\def\InterFeuilles{(1)*\DistanceInterFeuilles}

\draw[pointt] ({-2},{\NiveauA})--({18},{\NiveauA}) {};
\draw[pointt] ({-2},{\NiveauB})--({18},{\NiveauB}) {};
\draw[pointt] ({-2},{\NiveauC})--({18},{\NiveauC}) {};
\draw[pointt] ({-2},{\NiveauD})--({18},{\NiveauD}) {};

\node[noeud] (A) at ({(0)*\InterFeuilles},{\NiveauA}) {}; 
\node[noeud] (A0) at ({(-1)*\InterFeuilles},{\NiveauB}) {};
\node[noeud] (A1) at ({(1)*\InterFeuilles},{\NiveauB}) {};

\node[noeud] (A00) at ({(-1)*\InterFeuilles},{\NiveauC}) {};
\node[noeud] (A10) at ({(0.5)*\InterFeuilles},{\NiveauC}) {};
\node[noeud] (A11) at ({(1.5)*\InterFeuilles},{\NiveauC}) {};

\node[noeud] (A100) at ({(0.5)*\InterFeuilles},{\NiveauD}) {};

\draw[fleche,thick] (A)--(A0) {};
\draw[fleche,thick] (A)--(A1) {};
\draw[fleche,thick] (A0)--(A00) {};
\draw[fleche,thick] (A1)--(A10) {};
\draw[fleche,thick] (A1)--(A11) {};
\draw[fleche,thick] (A10)--(A100) {};

\node[noeud] (B0) at ({(3)*\InterFeuilles},{\NiveauA-1}) {}; 
\node[noeud] (B1) at ({(4)*\InterFeuilles},{\NiveauA}) {}; 
\node[noeud] (B2) at ({(5)*\InterFeuilles},{\NiveauB}) {}; 
\node[noeud] (B3) at ({(6)*\InterFeuilles},{\NiveauC}) {};
\node[noeud] (B4) at ({(7)*\InterFeuilles},{\NiveauB}) {}; 
\node[noeud] (B5) at ({(8)*\InterFeuilles},{\NiveauA}) {}; 
\node[noeud] (B6) at ({(9)*\InterFeuilles},{\NiveauB}) {}; 
\node[noeud] (B7) at ({(10)*\InterFeuilles},{\NiveauC}) {}; 
\node[noeud] (B8) at ({(11)*\InterFeuilles},{\NiveauD}) {}; 
\node[noeud] (B9) at ({(12)*\InterFeuilles},{\NiveauC}) {}; 
\node[noeud] (B10) at ({(13)*\InterFeuilles},{\NiveauB}) {}; 
\node[noeud] (B11) at ({(14)*\InterFeuilles},{\NiveauC}) {}; 
\node[noeud] (B12) at ({(15)*\InterFeuilles},{\NiveauB}) {}; 
\node[noeud] (B13) at ({(16)*\InterFeuilles},{\NiveauA}) {}; 
\node[noeud] (B14) at ({(17)*\InterFeuilles},{\NiveauA-1}) {}; 

\draw[fleche] (B0)--(B1) {};
\draw[fleche] (B1)--(B2) {};
\draw[fleche] (B2)--(B3) {};
\draw[fleche] (B3)--(B4) {};
\draw[fleche] (B4)--(B5) {};
\draw[fleche] (B5)--(B6) {};
\draw[fleche] (B6)--(B7) {};
\draw[fleche] (B7)--(B8) {};
\draw[fleche] (B8)--(B9) {};
\draw[fleche] (B9)--(B10) {};
\draw[fleche] (B10)--(B11) {};
\draw[fleche] (B11)--(B12) {};
\draw[fleche] (B12)--(B13) {};
\draw[fleche] (B13)--(B14) {};

\draw[arroww] (3,0)--(3,5) {};
\node[noeudtext] (T1) at (3,-0.5) {\scriptsize$0$};
\draw[arroww] (3,0)--(18,0) {};
\node[noeudtext] (T2) at (17,-0.5) {\scriptsize$14$};
\node[noeudtext] (T3) at (4,-0.5) {\scriptsize$1$};

\draw[red,thick,->] (-0.3,0) -- (-0.3,1) -- (-1.25,1.9) -- (-1.3,2) -- (-1.3,3) -- (-1.15,3.2)-- (-1,3.25) --(-0.85,3.2)-- (-0.7,3) -- (-0.7,2) -- (0,1.3) -- (0.7,2) -- (0.25,2.85) -- (0.2,3) -- (0.2,4) -- (0.35,4.2) -- (0.5,4.25) -- (0.65,4.2) -- (0.8,4) -- (0.8,3) -- (1,2.5) -- (1.2,3) -- (1.35,3.2) -- (1.5,3.25) -- (1.65,3.2) -- (1.8,3) -- (1.78,2.9) -- (1.3,2) -- (1.23,1.9) -- (0.3,1) -- (0.3,0.25){};

\draw[flech] (4,{\NiveauA-1.2})--(4,{\NiveauA-0.8}) {};
\draw[flech] (5,{\NiveauA-1.2})--(5,{\NiveauA-0.8}) {};
\draw[flech] (6,{\NiveauA-1.2})--(6,{\NiveauA-0.8}) {};
\draw[flech] (7,{\NiveauA-1.2})--(7,{\NiveauA-0.8}) {};
\draw[flech] (8,{\NiveauA-1.2})--(8,{\NiveauA-0.8}) {};
\draw[flech] (9,{\NiveauA-1.2})--(9,{\NiveauA-0.8}) {};
\draw[flech] (10,{\NiveauA-1.2})--(10,{\NiveauA-0.8}) {};
\draw[flech] (11,{\NiveauA-1.2})--(11,{\NiveauA-0.8}) {};
\draw[flech] (12,{\NiveauA-1.2})--(12,{\NiveauA-0.8}) {};
\draw[flech] (13,{\NiveauA-1.2})--(13,{\NiveauA-0.8}) {};
\draw[flech] (14,{\NiveauA-1.2})--(14,{\NiveauA-0.8}) {};
\draw[flech] (15,{\NiveauA-1.2})--(15,{\NiveauA-0.8}) {};
\draw[flech] (16,{\NiveauA-1.2})--(16,{\NiveauA-0.8}) {};

\draw[red,thick,->] (3,0.5) -- (6,3.5) -- (8,1.5) -- (11,4.5) -- (13,2.5) -- (14,3.5) -- (16.75,0.75){};

\end{tikzpicture}
}
\newcommand{\purarbredanscontour}{
\begin{tikzpicture}[xscale=0.35,yscale=0.5]
\tikzstyle{fleche}=[-,thin]

\tikzstyle{noeud}=[fill=black,circle,draw,scale=0.4]
\tikzstyle{noeudrrr}=[circle,draw,scale=0.4, very thick]
\tikzstyle{noeudred}=[fill=red,circle,draw,scale=0.4]
\tikzstyle{noeudtext}=[]
\tikzstyle{phantom}=[-,scale=0.05]

\tikzstyle{pointt}=[dashed,thin,draw=black!70!white]
\tikzstyle{fleched}=[dashed,thin,draw=black]
\tikzstyle{fleche}=[-,draw,thick]

\tikzstyle{arroww}=[->,draw,thick]

\def\DistanceInterNiveaux{1}
\def\DistanceInterFeuilles{1}

\def\NiveauA{(1)*\DistanceInterNiveaux}
\def\NiveauB{(2)*\DistanceInterNiveaux}
\def\NiveauC{(3*\DistanceInterNiveaux}
\def\NiveauD{(4*\DistanceInterNiveaux}
\def\InterFeuilles{(1)*\DistanceInterFeuilles}

\node[noeud] (B0) at ({(3)*\InterFeuilles},{\NiveauA-1}) {}; 
\node[noeud] (B1) at ({(4)*\InterFeuilles},{\NiveauA}) {}; 
\node[noeud] (B2) at ({(5)*\InterFeuilles},{\NiveauB}) {}; 
\node[noeud] (B3) at ({(6)*\InterFeuilles},{\NiveauC}) {};
\node[noeud] (B4) at ({(7)*\InterFeuilles},{\NiveauB}) {}; 
\node[noeud] (B5) at ({(8)*\InterFeuilles},{\NiveauA}) {}; 
\node[noeud] (B6) at ({(9)*\InterFeuilles},{\NiveauB}) {}; 
\node[noeud] (B7) at ({(10)*\InterFeuilles},{\NiveauC}) {}; 
\node[noeud] (B8) at ({(11)*\InterFeuilles},{\NiveauD}) {}; 
\node[noeud] (B9) at ({(12)*\InterFeuilles},{\NiveauC}) {}; 
\node[noeud] (B10) at ({(13)*\InterFeuilles},{\NiveauB}) {}; 
\node[noeud] (B11) at ({(14)*\InterFeuilles},{\NiveauC}) {}; 
\node[noeud] (B12) at ({(15)*\InterFeuilles},{\NiveauB}) {}; 
\node[noeud] (B13) at ({(16)*\InterFeuilles},{\NiveauA}) {}; 
\node[noeud] (B14) at ({(17)*\InterFeuilles},{\NiveauA-1}) {}; 

\draw[fleche] (B0)--(B1) {};
\draw[fleche] (B1)--(B2) {};
\draw[fleche] (B2)--(B3) {};
\draw[fleche] (B3)--(B4) {};
\draw[fleche] (B4)--(B5) {};
\draw[fleche] (B5)--(B6) {};
\draw[fleche] (B6)--(B7) {};
\draw[fleche] (B7)--(B8) {};
\draw[fleche] (B8)--(B9) {};
\draw[fleche] (B9)--(B10) {};
\draw[fleche] (B10)--(B11) {};
\draw[fleche] (B11)--(B12) {};
\draw[fleche] (B12)--(B13) {};
\draw[fleche] (B13)--(B14) {};

\draw[red,fleche,very thick] (6,3) -- (6,1) {};
\draw[red,fleche,very thick] (11,4) -- (11,2) {};
\draw[red,fleche,very thick] (14,3) -- (14,2) {};
\draw[red,fleche,very thick] (13,2) -- (13,1) {};
\draw[red,fleche,very thick] (7,1) -- (7,0) {};
\draw[red,fleche,very thick] (5.95,1) -- (13.05,1) {};
\draw[red,fleche,very thick] (10.95,2) -- (14.05,2) {};
\draw[red,fleche,very thick] (5.85,2) -- (6.15,2) {};
\draw[red,fleche,very thick] (10.85,3) -- (11.15,3) {};

\def\Shift{14}

\draw[red,fleche,very thick] (\Shift+6,3) -- (\Shift+6,1) {};
\draw[red,fleche,very thick] (\Shift+11,4) -- (\Shift+11,2) {};
\draw[red,fleche,very thick] (\Shift+14,3) -- (\Shift+14,2) {};
\draw[red,fleche,very thick] (\Shift+13,2) -- (\Shift+13,1) {};
\draw[red,fleche,very thick] (\Shift+7,1) -- (\Shift+7,0) {};
\draw[red,fleche,very thick] (\Shift+5.95,1) -- (\Shift+13.05,1) {};
\draw[red,fleche,very thick] (\Shift+10.95,2) -- (\Shift+14.05,2) {};
\draw[red,fleche,very thick] (\Shift+5.85,2) -- (\Shift+6.15,2) {};
\draw[red,fleche,very thick] (\Shift+10.85,3) -- (\Shift+11.15,3) {};

\node[noeudrrr] (AAA0) at ({(\Shift+7)*\InterFeuilles},{\NiveauA-0.5}) {};
\node[noeudrrr] (AAA1) at ({(\Shift+6)*\InterFeuilles},{\NiveauB-0.5}) {};
\node[noeudrrr] (AAA2) at ({(\Shift+6)*\InterFeuilles},{\NiveauC-0.5}) {};
\node[noeudrrr] (AAA3) at ({(\Shift+13)*\InterFeuilles},{\NiveauB-0.5}) {};
\node[noeudrrr] (AAA4) at ({(\Shift+11)*\InterFeuilles},{\NiveauC-0.5}) {};
\node[noeudrrr] (AAA5) at ({(\Shift+11)*\InterFeuilles},{\NiveauD-0.5}) {};
\node[noeudrrr] (AAA6) at ({(\Shift+14)*\InterFeuilles},{\NiveauC-0.5}) {};
\draw[fleche,very thick] (AAA0)--(AAA1) {};
\draw[fleche, very thick] (AAA1)--(AAA2) {};
\draw[fleche, very thick] (AAA0)--(AAA3) {};
\draw[fleche, very thick] (AAA3)--(AAA4) {};
\draw[fleche, very thick] (AAA4)--(AAA5) {};
\draw[fleche, very thick] (AAA3)--(AAA6) {};

\end{tikzpicture}
}

\begin{document}
	
\title[Inference for conditioned Galton-Watson trees]{Inference for conditioned Galton-Watson trees\\
from their Harris path}
	
\author{Romain Aza\"{\i}s}
\author{Alexandre Genadot}
\author{Benoit Henry}
	
\address{Laboratoire Reproduction et D\'eveloppement des Plantes, Univ Lyon, ENS de Lyon, UCB Lyon 1, CNRS, INRA, Inria, F-69342, Lyon, France.}
\address{Institut de Math\'ematiques de Bordeaux, Univ Bordeaux, CNRS, UMR 5251 and INRIA Bordeaux-Sud Ouest, Team CQFD, F-33400 Talence, France.}
\address{IMT Lille Douai, Universit\'e de Lille, Villeneuve d'Ascq, France }
	


\begin{abstract}
Tree-structured data naturally appear in various fields, particularly in biology where plants and blood vessels may be described by trees, but also in computer science because \verb+XML+ documents form a tree structure.
This paper is devoted to the estimation of the relative scale parameter of conditioned Galton-Watson trees.
New estimators are introduced and their consistency is stated. A comparison is made with an existing approach of the literature.
A simulation study shows the good behavior of our procedure on finite-sample sizes and from missing or noisy data. An application to the analysis of revisions of Wikipedia articles is also considered through real data.
\end{abstract}

\maketitle

\section{Introduction}

Many data are naturally modeled by an ordered tree structure: from blood vessels in biology to \verb+XML+ files
in computer science through the secondary structure of RNA in biochemistry. The statistical analysis
of a dataset of hierarchical records is thus of great interest.
In this paper, our aim is to propose new methods to estimate the scale parameter arising in Galton-Watson trees conditioned on their number of nodes from various statistical settings.

A Galton-Watson tree is the genealogical tree of a population starting from one initial ancestor (the root) in which each individual gives birth to
a random number of children according to the same probability distribution, independently of each other. In this article, we focus on Galton-Watson trees conditional on their number of nodes. Several main classes of random trees can be seen as conditioned Galton-Watson trees \cite{Dev12,Janson12}.
For instance, an ordered tree picked uniformly at random in the set of all ordered trees of a given size is a conditioned Galton-Watson tree with offspring distribution the geometric law with parameter $1/2$. In addition, an ordered tree picked uniformly at random in the set of $d$-ary trees, i.e., trees in which each node has no more than $d$ children, is a conditioned Galton-Watson tree with offspring distribution the binomial law with parameter $d$ and $1/d$. In particular, binary trees but also full binary trees (taking the uniform law on the set $\{0,2\}$ as offspring distribution) are thus encoded by a conditioned Galton-Watson model. Binary trees are widely used in computer science, through binary search trees \cite{knu81vol3} and Huffman coding \cite{knu85} commonly used for data compression. They also appear in biology in the approximation of phylogenetic trees \cite{aldous1996} for example. One also refers the reader to \cite[10. Examples of simply generated random trees]{Janson12} for other examples of conditioned Galton-Watson trees arising from random trees (that can even be unordered and labelled). To sum up, conditioned Galton-Watson trees model a large variety of random hierarchical structures. Developing specific statistical methods for this stochastic model is thus of first importance.

Any ordered tree may be encoded by its Harris path which returns height of nodes in depth-first order (see Subsection \ref{ss:harris}, Algorithm \ref{algo:dfs} and Figure \ref{fig:tree2contour}).
In \cite[Theorem 23]{aldous1993}, Aldous stated the following asymptotic property of the Harris path $\mathcal{H}[\tau_n]$ of a Galton-Watson tree $\tau_n$ conditioned on having
$n$ nodes,
\begin{equation}
\label{eq:intro:cv}
\left(\frac{\mathcal{H}[\tau_n](2nt)}{\sqrt{n}},~t\in[0,1]\right)\stackrel{(d)}{\longrightarrow} \left(\frac{2}{\sigma}\ee_t,~t\in[0,1]\right),
\end{equation}
in the uniform topology of $\mathcal{C}([0,1],\R)$, when $n$ goes to infinity  whenever the offspring distribution is $1$ on average with standard deviation $\sigma$ and $\mathbf{e}$ denotes
the normalized Brownian excursion. This means that conditioned Galton-Watson trees asymptotically share a common form (the so-called \textit{continuum random tree})
given by the Brownian excursion, and can be differentiated only by the scale parameter of interest $\sigma^{-1}$. This unknown quantity is to be estimated from only one tree or from a forest of independent trees generated from the same birth distribution.

Estimating (functions of) $\sigma$ from a forest of independent conditioned Galton-Watson trees has only been considered in a recent paper. The authors of \cite{KPDRV14} exploit a corollary of the weak convergence \eqref{eq:intro:cv} providing the asymptotic distribution of the height of a uniformly sampled node in the tree \cite[Proposition 4]{KPDRV14} to construct estimators of the variance $\sigma^2$ and develop asymptotic tests. It should be already noticed that estimation strategies based on the convergence in distribution \eqref{eq:intro:cv} can only lead to weak convergence results for estimators computed from a unique tree. The aim of the present paper is twofold. First, we establish in Theorem \ref{cv:varempir} that the empirical variance of the numbers of children is a consistent estimator of $\sigma^2$ (in particular even from the observation of only one tree), whereas, even if a Galton-Watson tree is generated from a sequence of i.i.d. random variables, this is not the case for the conditioned structure. This new result based on a corollary of Bartlett's formula shows that the empirical variance provides a better estimate of $\sigma^2$ than any other statistical method based on Aldoustheorem \eqref{eq:intro:cv}. Secondly, we propose two new estimation strategies for $\sigma^{-1}$ from a forest of independent conditioned Galton-Watson trees based on the weak convergence established by Aldous and we compare them with the procedure developed in \cite{KPDRV14}.

These estimation strategies rely on the idea motivated by the weak convergence \eqref{eq:intro:cv} that, on average, the normalized Harris paths of the forest should look like the expected process
$(2\sigma^{-1}E(t),t\in[0,1])$ at least asymptotically, where $E(t)=\mathbf{E}[\mathbf{e}(t)]$. The parameter $\sigma^{-1}$ can thus be expressed as the solution of a least square problem. Our first method consists in computing the least square estimator of $\sigma^{-1}$ from the concatenation of the Harris paths of the forest. We establish two results of convergence in Subsection \ref{ss:result:ls}. For only one Galton-Watson tree $\tau_n$ conditioned on having $n$ nodes, this estimator of $\sigma^{-1}$ is given (see Subsection \ref{ssec:adequacy}) by
$$
\widehat{\lambda}[\tau_n] = \frac{ \langle \mathcal{H}[\tau_n](2n\cdot) , E\rangle}{2\sqrt{n}\|E\|_2^2},
$$
where $\langle \cdot , \cdot\rangle$ is the scalar product of $\mathbf{L}^2([0,1],\mathbf{R})$. By virtue of the weak convergence \eqref{eq:intro:cv}, one may remark (see Corollary \ref{lim:lambda:n}) that
$$\widehat{\lambda}[\tau_n] \stackrel{(d)}{\longrightarrow} \sigma^{-1}\Lambda_\infty,$$
where $\Lambda_\infty= \frac{\langle \mathbf{e},E\rangle}{\|E\|^2_2}$.
Actually, the aforementioned least square estimator only exploits the average behavior of $\Lambda_\infty$
(in other words, the average asymptotic behavior of Harris paths) and not its complete distribution.
Our second strategy takes into account the shape of the distribution of $\Lambda_\infty$: we estimate $\sigma^{-1}$ by the parameter $x$ that aligns the theoretical distribution of $x\Lambda_\infty$ and the empirical measure of the $\widehat{\lambda}[\tau_{n_i}^i]$'s in terms of Wasserstein distance, the considered forest being composed of $N$ trees $\tau_{n_i}^i$. Convergence results are stated in Subsection \ref{ss:result:w}. We point out that the theoretical properties of $\Lambda_\infty$ are far from obvious. In particular, we establish by Malliavin calculus that $\Lambda_\infty$ is absolutely continuous w.r.t.\ the Lebesgue measure in Proposition \ref{prop:theoryLambdaInfty}, which is required in some proofs.

The authors of \cite{KPDRV14} do not focus on the problem of estimating $\sigma^{-1}$ but, for the sake of comparison, we rely on their approach to provide another estimator of this quantity. We compare these alternative strategies from both theoretical and numerical points of view. In particular, we show in Subsection \ref{ssec:adequacy} that the variances of our estimators are approximately $4$ times lower than the one of the estimator based on this competitive approach of the literature. Our results are better in terms of dispersion because the estimators take into account all the behavior of the tree and not only the behavior of a randomly chosen node. We also point out that the theoretical setting of \cite{KPDRV14} is slightly different because investigations are directly based on infinite trees (i.e., \textit{continuum random trees}, unobservable in practice) and not on large but finite trees.

At this step, one may wonder whether an approach that only yields weak convergence results is relevant considering the empirical variance is a consistent estimator of $\sigma^2$. Our idea is to explore statistical inference for trees from coding processes, i.e., via functional data analysis. This connection has been first established in the recent paper \cite{doi:10.1080/10618600.2013.786943}. In the present article we aim at investigating this strategy when the data have been generated from the stochastic model of conditioned Galton-Watson trees. In Subsection \ref{ss:missingnoisydata}, we prove from simulations that our estimators based on the weak convergence of Harris paths provide good results even from missing or noisy data, in particular when the empirical variance presents a large bias or can not be computed, showing the great interest of this approach.

The application of our estimators on simulated and real data in Section \ref{sec:simu} appears to be a non trivial task, in particular because it requires important preliminary computations. For this reason and to provide a turnkey solution, we have developed a \verb+Matlab+ toolbox that enables users to quickly and easily apply our methods to data. This toolbox as well as a detailed user documentation are available from the authors upon request.
The numerical experiments presented in Subsection \ref{ss:simulateddata} show that both our estimators and the approach developed in \cite{KPDRV14} are intrinsically biased for binary trees because of the approximation of the Harris paths of finite trees by the average Brownian excursion. Indeed, we empirically observe on simulation examples that Harris paths of binary trees weakly converge to the Brownian excursion from below (see Figure \ref{fig:biasvssigma}). As a consequence, we introduce a numerical correction of this negative bias, also implemented in the toolbox. The simulation study illustrates the good behavior of the corrected estimates on finite-sample sizes.

Visualizing the evolution of historical hierarchical data is a difficult issue in particular because such objects have no representation in a Euclidean space. This problem occurs in the study of the sequence of revisions of a given Wikipedia article. Indeed, the famous free Internet encyclopedia allows its users (the \textit{Wikipedians}) to edit almost any articles. Starting from the creation of a given article, the history of revisions is accessible and can be investigated to understand how the contributors agree on its structure, or to automatically detect vandalism\footnote{It frequently happens that malicious people willingly disrupt the content of an article, for instance, for political or ideological reasons.} \cite{Adler2011,Mola2010}. \textit{IBM's History Flow} is a visualization tool for documents in various stages of their development which has been applied to Wikipedia articles \cite{Viegas:2004:SCC:985692.985765,Viegas:2007:TBY:1255480.1255643}. We think that our method may be a complementary tool to this famous technique. Indeed the structure of \verb+HTML+ documents, such as Wikipedia articles, may be encoded by an ordered tree structure (see Figure \ref{fig:htmltree}). Furthermore, all the Wikipedia webpages share the same template, i.e., standardized \verb+HTML/CSS+ files, and thus can be differentiated by their relative scale. In Section \ref{ss:realdata}, we apply our estimators to the analysis of two Wikipedia articles. We highlight that Wikipedia articles undergo  ``running in'period before reaching some kind of steady state in which the contributors had agreed on the structure of the article. In addition, we show that our techniques may be used to detect improper editions of an article.

The organization of the paper is as follows. Section \ref{s:gw} is devoted to the formulation of the problem at hand: definition of conditioned Galton-Watson trees in Subsection \ref{ss:gw:def}, definition of Harris paths in Subsection \ref{ss:harris}, asymptotic behavior of Harris paths of conditioned Galton-Watson trees in Subsection \ref{ss:gw:as}. In addition, we state in Subsection \ref{ss:emp} the consistency of the empirical variance of the numbers of children. The two estimation procedures from Harris paths are presented in Section \ref{s:estim}, while Section \ref{sec:mainResults} focuses on the results of convergence.
Simulation techniques for conditioned Galton-Watson trees, numerical experiments and application to real data are presented in Sections \ref{sec:simu} and \ref{ss:realdata}. In particular, Subsection \ref{ss:missingnoisydata} is dedicated to the difficult context of missing or noisy data in which the empirical variance performs less well than our Harris paths-based estimators or even can not be computed.


\section{Conditioned Galton-Watson trees}
\label{s:gw}

\subsection{Definition}
\label{ss:gw:def}

Trees are connected graphs with no cycles. A rooted tree $\tau$ is a tree in which one node has been distinguished as the root, denoted by $r(\tau)$ (always drawn at the bottom of the tree in this paper). In this case, the edges are assigned a natural orientation, away from the root towards the leaves. One obtains a directed rooted tree in which there exists a parent-child relationship: the parent of a node $v$ is the first vertex met on the path to the root starting from $v$. The length of this path (in number of nodes) is called the height $h(v)$ of $v$. The set $c(v)$ of children of a vertex $v$ is the set of nodes that have $v$ as parent. An ordered or plane tree is a rooted tree in which an ordering has been specified for the set of children of each node, conventionally drawn from left to right. In this paper we consider ordered rooted trees simply referred to as trees. In addition, for any node $v$, $\tau[v]$ denotes the subtree of $\tau$ composed of $v$ and all of its descendants in $\tau$.

Intuitively, a Galton-Watson tree can be seen as a tree encoding the dynamic of a population generated from some offspring distribution $\mu$ on $\N$. A Galton-Watson tree $\tau$ with offspring distribution $\mu$ is a random ordered rooted tree constructed recursively as follows.
\begin{itemize}[label={$\diamond$}]
\item The number of children $\#c(r(\tau))$ emanating from the root is a random variable with law $\mu$. The first generation consists thus in $\#c(r(\tau))$ vertices.
\item Assume that the $n^\text{th}$ generation of children has been constructed and consists in a list of vertices $\mathcal{V}_n$. Then, the generation $n+1$ is constructed such that $\{\#c(v)\,:\,v\in\mathcal{V}_n\}$ is a collection of independent random variables with law $\mu$.
\end{itemize}
The asymptotic behavior of Galton-Watson trees may exhibit different regimes depending on the average number of children per capita,
$$\overline{\mu} = \sum_{k\geq0} k\mu(k),$$
where $\mu(k)$ is the measure of the singleton $\{k\}$ by $\mu$.
\begin{itemize}[label={$\diamond$}]
\item The subcritical case: $\overline{\mu}<1$. In this case, the average number of nodes is finite. This means that the population goes extinct almost surely.
\item The critical case: $\overline{\mu}=1$. The fact that the offspring distribution $\mu$ is critical also ensures the almost sure finiteness of the tree, except when $\mu(1)=1$ where the number of nodes is almost surely infinite. When $\mu(1)<1$, in contrary to the sub-critical case, the expected number of nodes is infinite.
\item The supercritical case: $\overline{\mu}>1$. In this case, the number of vertices is infinite with positive probability.
\end{itemize}
We use the notation $\text{GW}_n(\mu)$ for the distribution of Galton-Watson trees with offspring distribution $\mu$ conditioned on having $n$ nodes.

\begin{remark}
In this paper, we will always state our results in terms of critical Galton-Watson trees. However, this is not really a restriction since, as noted in \cite[6.3 Brownian asymptotics for conditioned Galton-Watson trees]{Pit06}, for any offspring distribution $\mu$, there exists a critical law $\mu'$ such that
$$\text{\normalfont GW}_n(\mu) \stackrel{(d)}{=} \text{\normalfont GW}_n(\mu').$$
In particular, this means that the average number of children $\overline{\mu}$ is not identifiable from conditioned Galton-Watson trees without some additional assumptions on $\mu$.
\end{remark}

\subsection{From ordered trees to Harris paths}
\label{ss:harris}

The Harris walk $\mathcal{H}[\tau]$ of an ordered rooted tree $\tau$ is defined from the depth-first search algorithm and the notion of height of nodes already presented in Subsection \ref{ss:gw:def}. Depth-first search is an algorithm for traversing a tree which one explores as far as possible along each branch before backtracking. The version of the algorithm used to define the Harris walk of a tree is presented in Algorithm \ref{algo:dfs}.

\SetKwFunction{FRecurs}{DFS}%
\SetKwProg{Fn}{Function}{\string:}{}
\begin{algorithm}[H]
\DontPrintSemicolon
\Fn(){\FRecurs{$\tau$\,,\,$l=\emptyset$}}{
\KwData{an ordered tree $\tau$}
\KwResult{vertices of $\tau$ in depth-first order}
add $r(\tau)$ to $l$\\
  \For{$v$ {\bf{in}} $c(r(\tau))$}{
    \If{r(t[v]) \bf{is not in} $l$}{
	call \FRecurs($t[v]$,$l$)\\
	add again $r(\tau)$ to $l$}
    }
\KwRet{l}
}
\vspace{2mm}
\caption{Recursive depth-first search.}
\label{algo:dfs}
\end{algorithm}

\begin{figure}
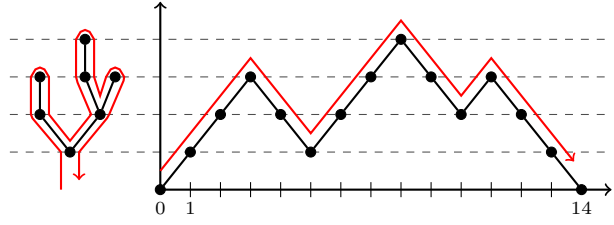

\purcontour
\caption{Construction of the Harris path (right) from $0$ to $2n=14$ as the contour of an ordered tree (left) with $n=7$ nodes.
}
\label{fig:tree2contour}
\end{figure}

\begin{remark}
In Algorithm \ref{algo:dfs}, each node $v$ appears $\#c(v)+1$ times. Starting from the root of a tree $\tau$, the result is thus a sequence of length
$$\sum_{v\in\tau} (\#c(v)+1) = \#\tau+\sum_{v\in\tau}\#c(v) = 2\#\tau-1,$$
because the root is the only vertex not to be counted.
\end{remark}

\begin{figure}[!h]
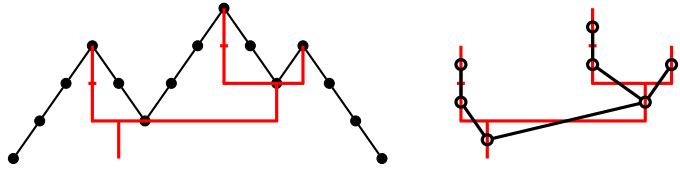

\centering
\purarbredanscontour
\caption{The ordered tree of Figure \ref{fig:tree2contour} in its Harris path (left): each vertical axis represents a node of the original structure (right). A common picture helping to see how to recover the tree from the contour is to imagine putting glue under the contour and then squeezing the contour together horizontally such that the inner parts of the contour which face each other are glued.}
\label{fig:contour2tree}
\end{figure}

The Harris walk $\mathcal{H}[\tau]$ of $\tau$ is defined as a sequence of integers indexed by the set $\{0,\dots,2\#\tau\}$ as follows:
\begin{itemize}[label={$\diamond$}]
\item $\mathcal{H}[\tau](0)=\mathcal{H}[\tau](2\#\tau)=0$,
\item for $1\leq k<2\#\tau$, $\mathcal{H}[\tau](k)=h(v)+1$ where $v$ is the $k^{\text{th}}$ node in depth-first traversal of $\tau$.
\end{itemize}
The Harris process is then defined as the linear interpolation of the Harris walk (see example in Figure \ref{fig:tree2contour}). Note that, as displayed in Figure \ref{fig:contour2tree}, the tree can be recovered from its Harris process such that the correspondence is one to one.

\subsection{Asymptotic behavior of Harris paths}
\label{ss:gw:as}

Let $\tau_n\sim\text{GW}_n(\mu)$ with $\overline{\mu}=1$. The variance of the offspring distribution $\mu$ is denoted by $\sigma^2$,
$$\sigma^2=\sum_{k\geq1} (k-1)^2\mu(k).$$
We focus on the asymptotic behavior of the Harris process $\mathcal{H}[\tau_n](2n\cdot)$ when $n$ tends to infinity. The convergence in distribution has been stated in \cite[Theorem 23]{aldous1993}.

\begin{figure}[!h]
	\centering
	\includegraphics[scale=0.25]{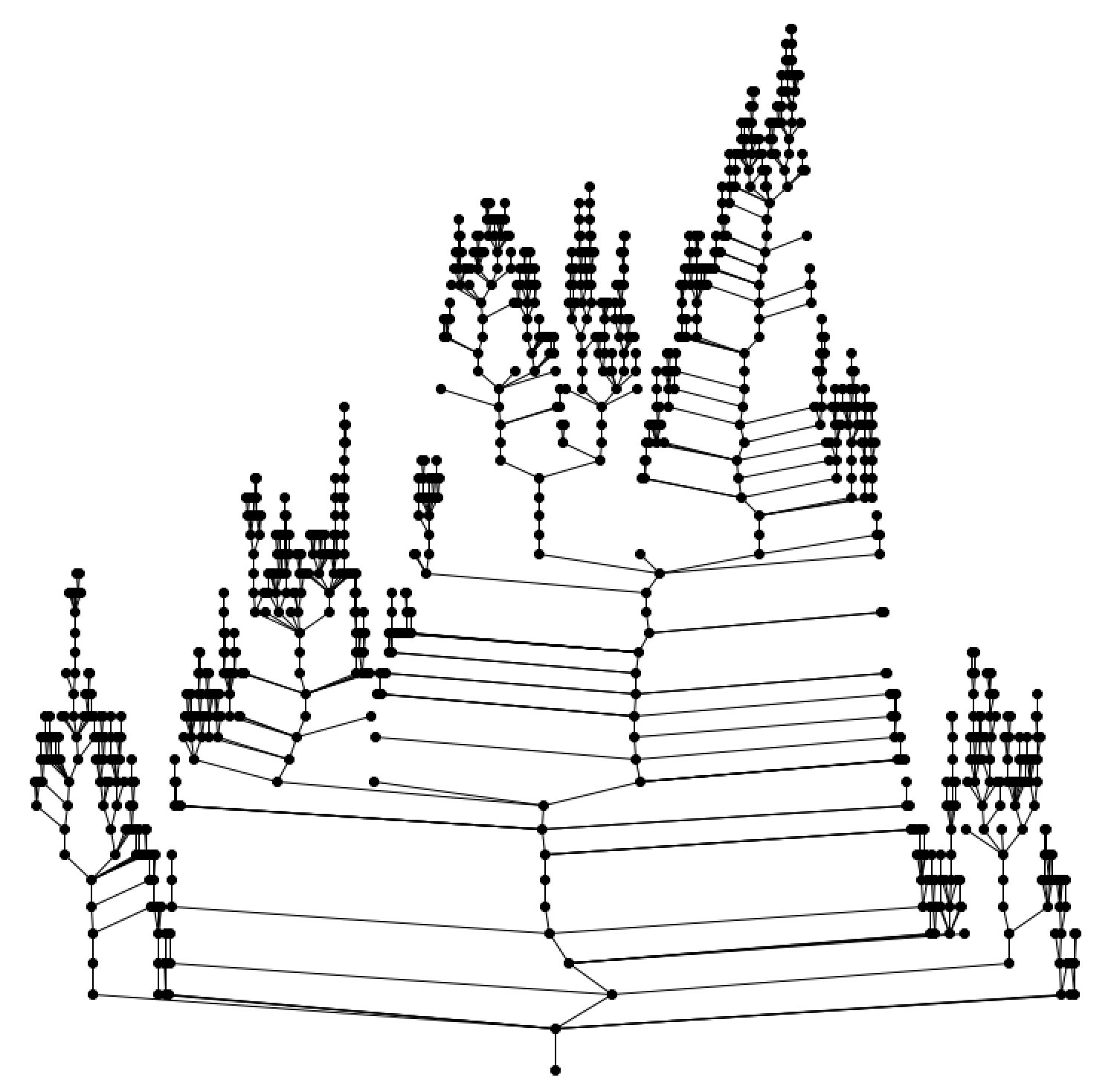}\\
	\includegraphics[scale=0.25,trim= 20mm 0mm 0mm 0mm]{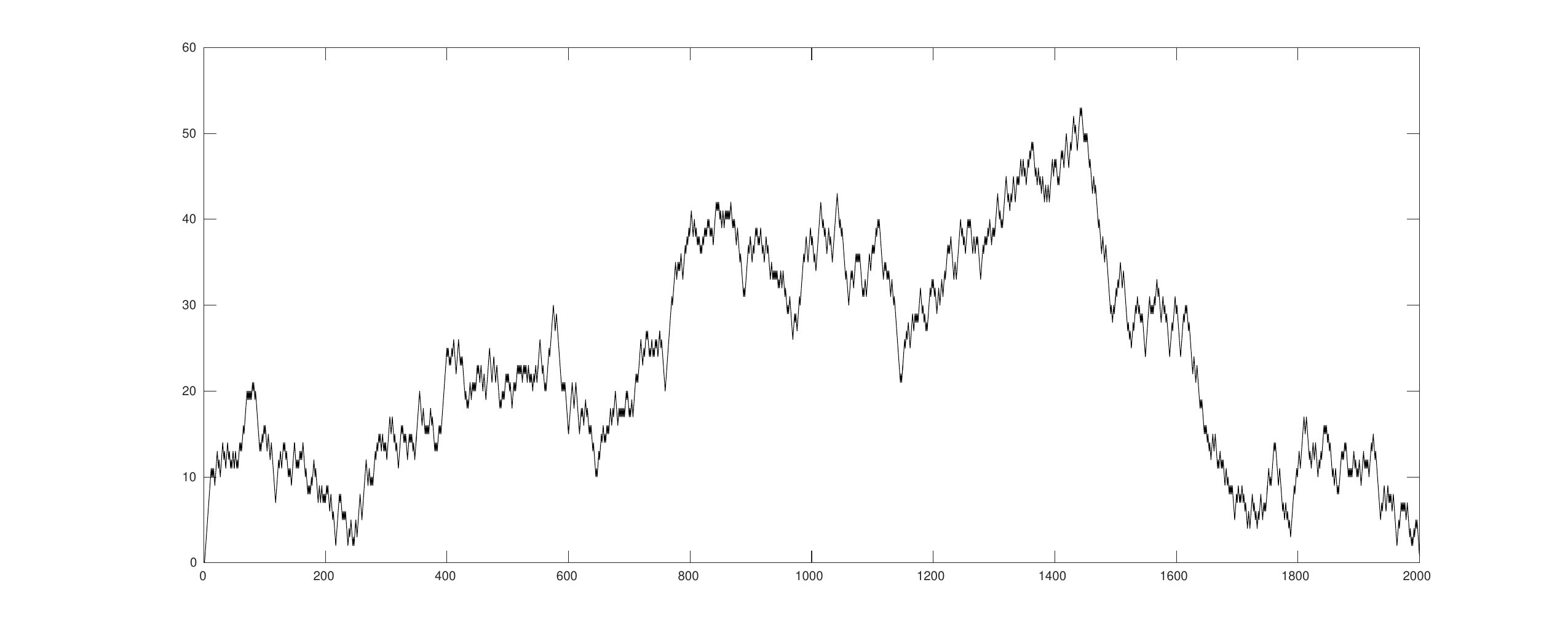}
	\caption{A Galton-Watson tree conditional on having $1000$ nodes generated from the geometric birth distribution with variance $\sigma^2=2$ (top) and its Harris path (bottom).}
	\label{fig:exempleGW}
\end{figure}

\begin{theorem}\label{thm:aldous}
When $n$ goes to infinity, we have
$$
\left(\frac{\mathcal{H}[\tau_n](2nt)}{\sqrt{n}},~t\in[0,1]\right)\stackrel{(d)}{\longrightarrow} \left(\frac{2}{\sigma}\ee_t,~t\in[0,1]\right),
$$
where $\ee$ is a standard Brownian excursion, the convergence holding in law in the space $\mathcal{C}([0,1],\R)$. An illustration of this convergence in distribution may be found in Figure \ref{fig:exempleGW}.
\end{theorem}

Let us simply recall that a standard Brownian excursion is a Brownian motion conditioned on being positive and on taking the value $0$ at time $1$. The density of $\ee_t$, for $0\leq t\leq 1$, is given in \cite[XI. 3. Bessel Bridges]{RY} and writes
\begin{equation*}
\forall\,x\in\R,~ f_{\ee_t}(x)=\sqrt{\frac2\pi}\frac{x^2}{\sqrt{t(1-t)}^3}\exp\left(-\frac{x^2}{2t(1-t)}\right)\mathbf{1}_{\R_+}(x).
\end{equation*}
From this, we can compute some simple functionals of the excursion. For instance, we have,
\begin{equation}\label{eq:Et}
\forall\,0\leq t\leq 1,\quad E_t=\E[\ee_t]= 4\sqrt{\frac{t(1-t)}{2\pi}}\quad\text{and}\quad \E\left[\ee_t^2\right]=3t(1-t).
\end{equation}
The easiest way to simulate a Brownian excursion is certainly from its identity in law with a three-dimensional Bessel bridge \cite[Theorem XII.4.2]{RY}, which is simply the Euclidean norm of a three-dimensional Brownian bridge,
\begin{equation}\label{eq:rep:3d}
(\ee_t,~t\in[0,1])\stackrel{\rm (d)}{=}\left(\sqrt{\sum_{i=1}^3\left(B^{i}_t-tB^{i}_1\right)^2
},~t\in[0,1]\right),
\end{equation}
where the $B^{i}$'s are three independent Brownian motions. The convergence presented in Theorem \ref{thm:aldous} also holds in expectation \cite[Theorem 1]{Marck05}.
\begin{theorem}\label{thm:cv:expectation}
When $n$ goes to infinity, we have,
$$
\forall\,0\leq t\leq 1,\quad\E\left[\frac{\mathcal{H}[\tau_n](2n t)}{\sqrt{n}}\right]\longrightarrow\frac{2}{\sigma}E_t ,
$$
where the function $(E_t,~0\leq t\leq 1)$ has been defined in \eqref{eq:Et}.
\end{theorem}

\begin{remark}
Theorem \ref{thm:aldous} establishes that, in the asymptotic regime, the shape of a conditioned Galton-Watson tree is given by the normalized Brownian excursion, regardless of the offspring distribution $\mu$. However, there is one scale parameter given by the inverse of the standard deviation of $\mu$. As a consequence, when $\mu$ is unknown, the only quantity of interest that one may access by asymptotic inference from Theorem \ref{thm:aldous} is $\sigma^{-1}$. From Section \ref{s:estim}, we shall focus on the estimation of $\sigma^{-1}$.
\end{remark}

\subsection{Empirical estimators}
\label{ss:emp}

The purpose of this section is to study the behavior of empirical estimators of the mean $\overline{\mu}$ and of the variance $\sigma^2$. These estimators are privileged candidates in the case of standard, i.e., unconditioned, Galton-Watson trees. Unfortunately, the lack of independency and homogeneity in the numbers of children of the nodes in conditioned Galton-Watson trees suggests such methods should not work in this more complex framework.

\subsubsection{Empirical mean}

We begin our study with the empirical mean. Let $\tau_{n}$ be a Galton-Watson tree with birth distribution $\mu$ conditioned on having $n$ nodes. Denote, for any $i$ in $\{1,\dots,n \}$, $X_{i}$ the number of children of the $i^\text{th}$ individual indexed in depth-first order in $\tau_{n}$. Set
\[
M[\tau_{n}]=\frac{1}{n}\sum_{i=1}^{n}X_{i}
\]
the empirical mean of the number of children of the individuals in $\tau_n$. However, it is easily seen that, whatever the underlying stochastic model,
\[
M[\tau_{n}]=\frac{\sharp\tau_n-1}{n}=1-\frac{1}{n} .
\]
As a consequence, this estimator is deterministic and always estimates $1$ asymptotically whatever the real mean of the birth distribution.

\subsubsection{Empirical variance}

This section is devoted to the study of the empirical variance,
$$V[\tau_n]=\frac{1}{n}\sum_{i=1}^{n}\left(X_{i}-M[\tau_n]\right)^2.$$
Our main result of convergence is given below.
\begin{theorem}
\label{cv:varempir}
When $n$ goes to infinity,
$$\mathbf{E}\left[\left|V[\tau_{n}] - \sigma^2\right|\right] \to 0 .$$
In particular, $V[\tau_{n}]$ converges to $\sigma^2$ in probability.
\end{theorem}
\textit{Proof.}
Let $(\xi_i)_{1\leq i\leq n}$ be a sequence of i.i.d. random variables distributed according to $\mu$. Using the work of \cite[Turning to random walks]{Dev12}, it is known that there exists a random permutation $\Sigma$ such that the random vector $\left(\xi_{\Sigma(1)},\dots,\xi_{\Sigma(n)} \right)$
conditioned to $\sum_{i=1}^{n}\xi_{i}=n-1$ is equal in distribution to $(X_{1},\dots,X_{n})$. However, since the empirical estimators are invariant up to permutation, it follows that we can work directly with the vector $\left(\xi_{1},\dots,\xi_{n} \right)$ conditioned to $\sum_{i=1}^{n}\xi_{i}=n-1$.
Consequently, our main goal is simply to prove, keeping in mind that $M[\tau_n]=1-1/n$, that
$$R(n) = \mathbf{E}\left[\left|\frac{1}{n}\sum_{i=1}^{n}\left(\xi_{i}-1+\frac{1}{n}  \right)^{2}-\sigma^{2} \right| \Bigg| \frac{1}{n-1}\sum_{i=1}^{n}\xi_{i} =1 \right]$$
converges to $0$ as $n$ goes to infinity.
The proof lies on the asymptotic behavior of conditional probabilities which were obtained in \cite{condExpec} and allows to get the expected convergence. Let $F$ be a measurable and bounded real valued function. Since $\xi_1/(n-1)$ goes to $0$ in probability as $n$ goes to infinity,
it follows, according to \cite[Theorem 1]{condExpec} (using the alternative hypothesis of Remark 2.9 of this paper), that the conditional expectation
\[
\mathbf{E}\left[F(\xi_{1})  \Bigg| \frac{1}{n-1}\sum_{i=1}^{n}\xi_{i} =1  \right]
\]
converges to $\mathbf{E}\left[F(\xi_{1}) \right]$, as $n$ goes to infinity.
Since, $\mu$ does not have necessarily third order moments, we need to use a truncation method. Hence, let us consider, for any positive integer $k$,
$$ R(n) = \mathbf{E}\left[\left|\frac{1}{n}\sum_{i=1}^{n}\left(\xi_{i}\mathbf{1}_{\xi_{i}>k}+\xi_{i}\mathbf{1}_{\xi_{i}\leq k}-1+\frac{1}{n}  \right)^{2}-\sigma^{2} \right| \Bigg| \frac{1}{n-1}\sum_{i=1}^{n}\xi_{i} =1 \right].$$
Now, since
\begin{align*}
& \left(\xi_{i}\mathbf{1}_{\xi_{i}>k}+\xi_{i}\mathbf{1}_{\xi_{i}\leq k}-1+\frac{1}{n}  \right)^{2}\\
& =\left(\xi_{i}\mathbf{1}_{\xi_{i}\leq k}-1+\frac{1}{n}  \right)^{2}+\left(\xi_{i}\mathbf{1}_{\xi_{i}>k} \right)^{2}+2\xi_{i}\mathbf{1}_{\xi_{i}>k}\left(\frac{1}{n}-1 \right),
\end{align*}
we get
$$R(n) \leq T_1(n,k) + T_2(n,k) + 2\left(1-\frac{1}{n}\right)T_3(n,k) ,$$
where
\begin{eqnarray*}
T_1(n,k) &=& \mathbf{E}\left[\left|\frac{1}{n}\sum_{i=1}^{n}\left(\xi_{i}\mathbf{1}_{\xi_{i}\leq k}-1+\frac{1}{n}  \right)^{2}-\sigma^{2} \right| \Bigg| \frac{1}{n-1}\sum_{i=1}^{n}\xi_{i} =1 \right],\\
T_2(n,k) &=& \mathbf{E}\left[\left|\frac{1}{n}\sum_{i=1}^{n}\left(\xi_{i}\mathbf{1}_{\xi_{i}> k} \right)^{2} \right| \Bigg| \frac{1}{n-1}\sum_{i=1}^{n}\xi_{i} =1 \right] ,\\
T_3(n,k) &=&\mathbf{E}\left[\left|\frac{1}{n}\sum_{i=1}^{n}\xi_{i}\mathbf{1}_{\xi_{i}>k} \right| \Bigg| \frac{1}{n-1}\sum_{i=1}^{n}\xi_{i} =1 \right] .
\end{eqnarray*}
In order to treat $T_{1}(n,k)$, we consider
$$Q(n,k) = \mathbf{E}\left[\left(\frac{1}{n}\sum_{i=1}^{n}\left(\xi_{i}\mathbf{1}_{\xi_{i}\leq k}-1+\frac{1}{n}  \right)^{2}-\sigma^{2} \right)^{2} \Bigg| \frac{1}{n-1}\sum_{i=1}^{n}\xi_{i} =1 \right] .$$
We have
\begin{align*}
Q(n,k)&=\frac{1}{n^{2}}\sum_{1\leq i,j\leq n}\mathbf{E}\Bigg[\left(\left(\xi_{i}\mathbf{1}_{\xi_{i}\leq k}-1+\frac{1}{n}  \right)^{2}-\sigma^{2} \right)\\
& \times\left(\left(\xi_{j}\mathbf{1}_{\xi_{j}\leq k}-1+\frac{1}{n}  \right)^{2}-\sigma^{2} \right)\Bigg| \frac{1}{n-1}\sum_{i=1}^{n}\xi_{i} =1 \Bigg].
\end{align*}
Using the exchangeability of the vector $(\xi_{1},\dots,\xi_{n})$ under $\mathbf{P}\left(\cdot \mid \sum_{i=1}^{n}\xi_{i}=n-1 \right)$, we get
\begin{align*}
Q(n,k) &=\frac{1}{n}\mathbf{E}\left[\left(\left(\xi_{1}\mathbf{1}_{\xi_{1}\leq k}-1+\frac{1}{n}  \right)^{2}-\sigma^{2} \right)^{2}\Bigg| \frac{1}{n-1}\sum_{i=1}^{n}\xi_{i} =1 \right]\\
&+\frac{n(n-1)}{n^2}\mathbf{E}\Bigg[\left(\left(\xi_{1}\mathbf{1}_{\xi_{1}\leq k}-1+\frac{1}{n}  \right)^{2}-\sigma^{2} \right)\\
&\times\left(\left(\xi_{2}\mathbf{1}_{\xi_{2}\leq k}-1+\frac{1}{n}  \right)^{2}-\sigma^{2} \right)\Bigg| \frac{1}{n-1}\sum_{i=1}^{n}\xi_{i} =1 \Bigg].
\end{align*}
Then, \cite[Theorem 1]{condExpec} allows to understand the asymptotic behavior of both terms in the above sum,
\begin{align*}
& \lim_{n\to\infty}\mathbf{E}\left[\left(\left(\xi_{1}\mathbf{1}_{\xi_{1}\leq k}-1+\frac{1}{n}  \right)^{2}-\sigma^{2} \right)^{2}\Bigg| \frac{1}{n-1}\sum_{i=1}^{n}\xi_{i} =1 \right] \\
& = \mathbf{E}\left[\left(\left(\xi_{1}\mathbf{1}_{\xi_{1}\leq k}-1  \right)^{2}-\sigma^{2} \right)^{2} \right]>0,
\end{align*}
and
\[
\mathbf{E}\left[\left(\left(\xi_{1}\mathbf{1}_{\xi_{1}\leq k}-1+\frac{1}{n}  \right)^{2}-\sigma^{2} \right)\left(\left(\xi_{2}\mathbf{1}_{\xi_{2}\leq k}-1+\frac{1}{n}  \right)^{2}-\sigma^{2} \right)\Bigg| \frac{1}{n-1}\sum_{i=1}^{n}\xi_{i} =1 \right]
\]
goes, as $n$ tends to infinity, to
\begin{align*}
& \mathbf{E}\left[\left(\left(\xi_{1}\mathbf{1}_{\xi_{1}\leq k}-1 \right)^{2}-\sigma^{2} \right)\left(\left(\xi_{2}\mathbf{1}_{\xi_{2}\leq k}-1  \right)^{2}-\sigma^{2} \right)\right]\\
& =\left(\mathbf{E}\left[\left(\xi_{1}\mathbf{1}_{\xi_{1}\leq k}-1 \right)^{2} \right]-\sigma^{2} \right)^2\to0,
\end{align*}
when $k$ goes to infinity. Hence, we get, for any positive integer $k$,
\begin{equation}
\label{eq:vareq1}
\lim_{n\to\infty} Q(n,k) = \left(\mathbf{E}\left[\left(\xi_{1}\mathbf{1}_{\xi_{1}\leq k}-1 \right)^{2} \right]-\sigma^{2} \right)^2,
\end{equation}
without any assumption on the moments of $\mu$ because of the truncation $\xi_{1}\mathbf{1}_{\xi_{1}\leq k}$. 
First, the Cauchy-Schwarz inequality entails $T_1(n,k) \leq \sqrt{Q(n,k)}$,
which gives, according to \eqref{eq:vareq1},
\[
\limsup_{n\to\infty}T_{1}(n,k)\leq\left|\mathbf{E}\left[\left(\xi_{1}\mathbf{1}_{\xi_{1}\leq k}-1 \right)^{2} \right]-\sigma^{2}\right|.
\]
In addition, we have according to \cite[Theorem 1]{condExpec},
\[
\limsup_{n\to\infty}T_{2}(n,k)\leq\limsup_{n\to\infty} \mathbf{E}\left[\left(\xi_{i}\mathbf{1}_{\xi_{i}> k} \right)^{2} \Bigg| \frac{1}{n-1}\sum_{i=1}^{n}\xi_{i} =1 \right]=\mathbf{E}\left[\left(\xi_{i}\mathbf{1}_{\xi_{i}> k} \right)^{2}  \right].
\]The case of $T_{3}(n,k)$ can be treated similarly and leads to
$$
\limsup_{n\to\infty} R(n) \leq\left|\mathbf{E}\left[\left(\xi_{1}\mathbf{1}_{\xi_{1}\leq k}-1 \right)^{2} \right]-\sigma^{2}\right|+\mathbf{E}\left[\left(\xi_{1}\mathbf{1}_{\xi_{1}> k} \right)^{2}  \right]+2\mathbf{E}\left[\left|\xi_{1}\mathbf{1}_{\xi_{1}>k} \right|\right].
$$
Now, letting $k$ going to infinity leads to the result.\hfill$\Box$


\section{Estimation procedure}
\label{s:estim}

In the previous part, we have shown that the empirical variance is a consistent estimator of $\sigma^2$ in critical conditioned Galton-Watson trees despite the lack of independency in the numbers of children. In this section, we aim at developing estimation procedures for $\sigma^2$ that exploit the weak convergence towards the Brownian excursion stated by Aldous and presented in Theorem \ref{thm:aldous}.

\subsection{Adequacy of the Harris path with the expected contour}
\label{ssec:adequacy}
Let $\tau_n\sim\text{GW}_n(\mu)$ with $\overline{\mu}=1$. We assume that the offspring distribution $\mu$ is unknown. By virtue of Theorem \ref{thm:cv:expectation}, the asymptotic average behavior of the normalized Harris process $(n^{-1/2}\mathcal{H}[\tau_n](2n t),~0\leq t\leq1)$ is given by $(2\sigma^{-1}E_t,~0\leq t\leq1)$, where $\sigma^{-1}$ is obviously also unknown. We propose to estimate $\sigma^{-1}$ by minimizing the $\mathbf{L}^2$-error defined by
$$\lambda\mapsto\left\|\frac{\mathcal{H}[\tau_n](2n\cdot)}{\sqrt{n}} - 2\lambda E\right\|_2^2,$$
where, and in all the sequel, $\mathbf{L}^2=\mathbf{L}^2([0,1],\mathbf{R})$ and its usual norm is denoted $\|\cdot\|_{2}$ for the sake of readability.
The solution of this least square problem is well-known and is given by
\begin{equation}\label{eq:def:lambdahat}
\widehat{\lambda}[\tau_n] = \frac{ \langle \mathcal{H}[\tau_n](2n\cdot) , E\rangle}{2\sqrt{n}\|E\|_2^2},
\end{equation}
where $\langle\cdot,\cdot\rangle$ is the scalar product of $\mathbf{L}^{2}$.
\begin{corollary}\label{lim:lambda:n}
When $n$ goes to infinity, we have
$$\widehat{\lambda}[\tau_n] \stackrel{(d)}{\longrightarrow} \sigma^{-1}\Lambda_\infty ,$$
where the random variable $\Lambda_\infty$ is defined by
\begin{equation}\label{eq:lambdainfty}
\Lambda_\infty = \frac{\langle \mathbf{e},E\rangle}{\|E\|^2_2}.\end{equation}
\end{corollary}
\textit{Proof.} The result directly follows from Theorem \ref{thm:aldous} because the functional $x\mapsto\langle x,E\rangle$ is continuous on $\mathcal{C}([0,1],\mathbf{R})$.
\hfill$\Box$

\begin{remark}
The convergence in distribution stated in Corollary \ref{lim:lambda:n} seems quite unsatisfactory because this means that $\widehat{\lambda}[\tau_n]$ is not a consistent estimator of $\sigma^{-1}$ and the least square strategy thus seems like inadequate in regards to the consistency of $V[\tau_n]$. In the sequel, we shall focus on the estimation of the parameter of interest $\sigma^{-1}$ from a forest of conditioned Galton-Watson trees as in \cite{KPDRV14}, only chance to get consistent estimates from Aldoustheorem. As mentioned in the introduction, our goal in this paper is to explore statistical inference for trees via functional data analysis of their Harris paths.
\end{remark}

Computing $\widehat{\lambda}[\tau_n]$ is a first step in the estimation of the inverse standard deviation from a large number of conditioned Galton-Watson trees. As a consequence, the distribution of the limit variable $\Lambda_\infty$ is of first importance.


\begin{lemma}\label{lem:cov:e}
For any $0\leq t<u\leq 1$, we have
$$
\E[\ee_t\ee_u]=\frac{2}{\pi}\left[3\sqrt{t(u-t)(1-u)}+(2t(1-u)+u(1-t))\arcsin\left(\sqrt{\frac{t(1-u)}{u(1-t)}}\right)\right].
$$
\end{lemma}

\textit{Proof.} This identity is derived from the joint density of $(\ee_t,\ee_u)$ given in \cite[XI. 3. Bessel Bridges]{RY}. The density of $(\ee_t,\ee_u)$ for $0\leq t<u\leq 1$ is given, for any positive numbers $x$ and $y$, by
\begin{align*}
& f_{t,u}(x,y)\\
&=\frac{2xy}{\pi\sqrt{(u-t)t^3(1-u)^3}}\sinh\left(\frac{xy}{u-t}\right)\exp\left(-\frac{x^2u}{2t(u-t)}\right)\exp\left(-\frac{y^2(1-t)}{2(1-u)(u-t)}\right).
\end{align*}
Thus,
\begin{align*}
& \E[\ee_t\ee_u] \\
& =\int_0^\infty \!\!\! \int_0^\infty \!\! \frac{2}{\pi}\frac{x^2y^2\sinh\left(\frac{xy}{u-t}\right)}{\sqrt{(u-t)t^3(1-u)^3}}\exp\left(-\frac{x^2u}{2t(u-t)}\right)\exp\left(-\frac{y^2(1-t)}{2(1-u)(u-t)}\right)\dd x\dd y.
\end{align*}
At this point, one can use the power series of $\sinh$ to separate the variables $x$ and $y$ and obtain
$$
\E[\ee_t\ee_u]= \frac{2}{\pi}\frac{\sqrt{(u-t)^5}}{\sqrt{u^3(1-t)^3}}\sum_{k=0}^\infty \frac{4^{k+1}((k+1)!)^2}{(2k+1)!}\left(\sqrt{\frac{t(1-u)}{u(1-t)}}\right)^{2k+1}.
$$
Note that when $0\leq t<u\leq 1$ we have indeed $\sqrt{\frac{t(1-u)}{u(1-t)}}<1$. Then, using
$$
\sum_{k=0}^\infty \frac{4^{k+1}((k+1)!)^2}{(2k+1)!}x^{2k+1}=\frac{3x\sqrt{1-x^2}+(2x^2+1)\arcsin(x)}{\sqrt{(1-x^2)^5}},
$$
we obtain the desired expression for $\E[\ee_t\ee_u]$.\hfill$\Box$


\begin{proposition}
\label{prop:theoryLambdaInfty}
The random variable $\Lambda_\infty$ admits a density $f_{\Lambda_\infty}$ w.r.t.\ the Lebesgue measure. Furthermore,
\begin{equation}
\label{eq:variance:lambdainfty}
\E[\Lambda_\infty]=1\quad\text{and}\quad\Var(\Lambda_\infty)=\frac{1}{\|E\|^4_2}\int_0^1\int_0^1 g(s,u)\,E_s\,E_u\,\dd s\,\dd u\,-\,1,
\end{equation}
where the mapping $g:[0,1]^2\to\R_+$ is defined from
$$g(t,u) = \frac{2}{\pi}\left[3\sqrt{t(u-t)(1-u)}+(2t(1-u)+u(1-t))\arcsin\left(\sqrt{\frac{t(1-u)}{u(1-t)}}\right)\right]$$
if $0\leq t\leq u\leq 1$ and $g(t,u)=g(u,t)$ otherwise.
\end{proposition}


\textit{Proof.}
We consider the probability space $(\mathcal{C}([0,1],\R^{3}),\mathcal{F},\mathbf{W})$, where $\mathcal{C}([0,1],\R^{3})$ is endowed with the uniform topology, $\mathcal{F}$ is the corresponding Borel $\sigma$-field and $\mathbf{W}$ is the Wiener measure. Let $T$ be the continuous linear operator defined by
\[
\begin{array}{cccc}
T:&\mathcal{C}([0,1],\R^{3}) & \to & \mathcal{C}([0,1],\R^{3}), \\
& \varphi & \mapsto & \left(T\varphi(s)=\varphi(s)-s\varphi(1) \right). \\
\end{array}
\]
Let also $\Gamma$ be the following function,
\[
\Gamma:\varphi\mapsto\int_{0}^{1}\norm{\varphi(s)}_{2}\frac{E_{s}}{\|E\|^2_2}\dd s.
\]
where $\norm{x}_{2}$ denotes the Euclidian norm on $\R^{3}$. With these notations and \eqref{eq:rep:3d}, we have that the pushforward measure of $\mathbf{W}$ through the application
\[F:\varphi\mapsto\Gamma(T\varphi),
\]
is the law of $\Lambda_{\infty}$. In other words, the random variable $F$ is equal in distribution to $\Lambda_{\infty}$. Now for every $\varphi$ in $\mathcal{C}([0,1],\R^{3})$ such that $Leb\left(\overline{\{t\in\R_{+}\,:\,\varphi(t)=0 \} }\right)=0$, we have that $\Gamma$ is Fr\'echet differentiable at the point $\varphi$. Moreover, the derivative at such point $\varphi$ is given by  \[
\begin{array}{cccc}
D_{\varphi}\Gamma:&\mathcal{C}([0,1],\R^{3}) & \to & \R, \\
& h & \mapsto & \int_{0}^{1}\frac{\left(\varphi(s),h(s)\right)}{\norm{\varphi(s)}_{2}}{\frac{E_{s}}{\|E\|^2_2}}\ \dd s, \\
\end{array}
\]
where $(\cdot,\cdot)$ denotes the Euclidean scalar product on $\mathbf{R}^{3}$.
Indeed, some straightforward manipulations give
\begin{align*}
& \int_{0}^{1}\!\!\left[\norm{\varphi(s)\!+\!h(s)}_{2}-\norm{\varphi(s)}_{2}\!-\!\frac{(\varphi(s),h(s))}{\norm{\varphi(s)}_{2}}\right]\!\!\frac{E_{s}}{\|E\|^2_2}\dd s\\
&=\!\!\int_{0}^{1}\!\!\left[\frac{\norm{h(s)}^{2}_{2}+(\varphi(s),h(s))\left(1\!-\!\frac{\norm{\varphi(s)+h(s)}_{2}}{\norm{\varphi(s)}_{2}} \right)}{\norm{\varphi(s)+h(s)}_{2}+\norm{\varphi(s)}_{2}}\right]\!\! {\frac{E_{s}}{\|E\|^2_2}}\dd s.
\end{align*}
Now, since $\frac{E_{s}}{\|E\|^2_2}\leq\frac{3\sqrt\pi}{2\sqrt2} $ and using the Cauchy-Schwarz inequality, we obtain
\begin{align*}
& \Bigg|\int_{0}^{1}\Bigg[\norm{\varphi(s)+h(s)}_{2}-\norm{\varphi(s)}_{2}-\frac{(\varphi(s),h(s))}{\norm{\varphi(s)}_{2}}\Bigg]\frac{E_{s}}{\|E\|^2_2} \dd s\Bigg|\\
&\leq\frac{3\sqrt\pi}{2\sqrt2} \int_{0}^{1}\left[\frac{\norm{h(s)}_{2}^{2}+\norm{h(s)}_{2}\Bigg|\norm{\varphi(s)}_{2}-\norm{\varphi(s)+h(s)}_{2} \Bigg|}{\norm{\varphi(s)+h(s)}_{2}+\norm{\varphi(s)}_{2}}\right] \dd s\\
&\leq\frac{3\sqrt\pi}{2\sqrt2}\|h\|_{\infty} \int_{0}^{1}\left[\frac{\norm{h(s)}_{2}+\bigg|\norm{\varphi(s)}_{2}-\norm{\varphi(s)+h(s)}_{2} \bigg|}{\norm{\varphi(s)+h(s)}_{2}+\norm{\varphi(s)}_{2}}\right] \dd s,
\end{align*}
with $\|h\|_{\infty}=\sup_{s\in[0,1]}| h(s)|_{2}$. Since 
\[
\int_{0}^{1}\left[\frac{\norm{h(s)}_{2}+\bigg|\norm{\varphi(s)}_{2}-\norm{\varphi(s)+h(s)}_{2} \bigg|}{\norm{\varphi(s)+h(s)}_{2}+\norm{\varphi(s)}_{2}}\right] \dd s
\]
is well-defined (because the integrand is bounded by $2$) and goes to zero as $\|h\|_{\infty}$ goes to zero, this proves that $D_{\varphi}\Gamma$ is the Fr\'echet derivative of $\Gamma$ at point $\varphi$. The functional $T$ being linear, $F$ is also Fr\'echet differentiable with Fr\'echet derivative given by
\[
\begin{array}{cccc}
D_{\varphi}F:&(\mathcal{C}([0,1],\R^{3}) & \to & \R, \\
& h & \mapsto & \int_{0}^{1}\frac{\left(T\varphi(s),Th(s)\right)}{\norm{T\varphi(s)}_{2}}{\frac{E_{s}}{\|E\|^2_2}}\ \dd s. \\
\end{array}
\]
Moreover, let $h$ be an element of $\mathbf{L}^{2}([0,1],\R^3)$, we have, since $\|E\|^2_2=\frac{4}{3\pi}$,
\begin{align*}
\left|F\left(\omega +\int_{0}^{\cdot}h(s)\dd s\right)-F(\omega)\right|&\leq \frac{3\pi}{4}\int_{0}^{1}\left\{\left|{\int_{0}^{t}h(s)\dd s}\right|_2+t\left|{\int_{0}^{1}h(s)\dd s}\right|_2\right\} E_{t}\, \dd t.
\end{align*}
But in the right hand side of the last inequality, we have, using Jensen's inequality,
\begin{align*}
& \int_{0}^{1}\left\{\sqrt{\sum_{i=1}^{3}\left(\int_{0}^{t}h^{i}(s)\dd s\right)^{2}}+t\sqrt{\sum_{i=1}^{3}\left(\int_{0}^{1}h^{i}(s)\dd s\right)^{2}}\right\}E_{t}\ \dd t\\&\leq\int_{0}^{1}\sqrt{\sum_{i=1}^{3}\left(\int_{0}^{1}h^{i}(s)^{2}\dd s\right)}(1+t)E_{t}\ \dd t\\&=\int_{0}^{1}\|h\|_{\mathbf{L}^{2}([0,1],\R^{3})} (1+t)E_{t}\ \dd t.
\end{align*}
From this, using the results of \cite[p.\,35]{nualard}, we have that $F$ belongs to the space $\mathbf{D}^{1,2}$, which is the domain of the Malliavin operator $D$ in $\mathbf{L}^{2}([0,1],\mathbf{R}^{3})$ (see \cite[pp.\,25--27]{nualard} for more details). Before going further let us recall some facts on Malliavin derivative. When working with the probability space $(\mathcal{C}([0,1],\R^{3}),\mathcal{F},\mathbf{W})$, it is known \cite[1.2.1 The derivative operator in the white noise case]{nualard} that there exist strong connections between Malliavin derivative and Fr\'echet derivative for a random variable $G$ of $\mathbf{D}^{1,2}$ defined from $(\mathcal{C}([0,1],\R^{3}),\mathcal{F},\mathbf{W})$ to $\mathbf{R}$. Since, the Fr\'echet derivative $D_{\omega}G$ at point $\omega$ of $G$ is a continuous linear form from $\mathcal{C}([0,1],\R^{3})$ into $\mathbf{R}$, it can be identified to a triple $(\mu^{\omega}_{1},\mu^{\omega}_{2},\mu^{\omega}_{3})$ of $\sigma$-finite measures on $\mathbf{R}$ such that,
\[
\forall\,h\in \mathcal{C}([0,1],\R^{3}),~D_{\varphi}Gh=\sum_{i=1}^{3}\int_{[0,1]}h^{i}(s)\ \mu^{\omega}_{i}(ds).
\]
In such a case, the Malliavin derivative of $G$ is the random process belonging to $\mathbf{L}^{2}([0,1],\R^{3})$ given by
$$\big\{\left(\mu^{\omega}_{1}(u,1],\mu^{\omega}_{2}(u,1],\mu^{\omega}_{3}(u,1]\right)~:~u\in[0,1] \big\}.$$
In our case, it follows that the Malliavin derivative of $F$ is given by
\[
DF(\omega)=\left(\int_{0}^{1}\frac{(\omega_{s}-s\omega_{1})E_{s}}{\norm{\omega_{s}-s\omega_{1}}_{2}\|E\|^2_2}(\mathbf{1}_{s>u}-s)\dd s,\ u\in[0,1]\right)\in\mathbf{L}^{2}([0,1],\R^{3}).
\]
Now, since $DF$ is not zero in $\mathbf{L}^{2}([0,1],\R^{3})$ for $\mathbf{W}$-almost every $\omega$, we get, together with \cite[Theorem 2.1.2]{nualard}, the existence of a density for $F$ w.r.t.\ the Lebesgue measure. The calculation of the variance is derived from the expectation of $\ee_t\ee_s$, $(s,t)\in[0,1]^2$, stated in Lemma \ref{lem:cov:e}.
\hfill$\Box$

\begin{remark}
The existence of a density was already known for the random variable $\int_{0}^{1}\ee_{s}\dd s$ \cite{louchard,jansonEx} but to the best of our knowledge no paper investigates the existence of a density for $\Lambda_{\infty}$. In these papers the study is performed thanks to the analysis of the double Laplace transform
\[
\lambda\mapsto\int_{0}^{\infty}\exp(-\lambda t)\E\left[\exp\left({-t \int_{0}^{1}\ee_{s}\dd s}\right)\right]\dd t.
\]
Thanks to the Feynman-Kac formula, the authors express this quantity in terms of Airy functions. Then, they take the inverse of the Laplace transform via analytical methods. 
Unfortunately, their method does not extend to our case. Indeed, in their case, an expression of the double Laplace transform given above is derived from the Feynman-Kac formula for standard Brownian motion which tells us that the function
\[
u(t,x)=\mathbf{E}_{x}\left[f(B_{t})\exp\left(\int_{0}^{t}B_{s}\dd s\right) \right], \quad \forall\,(t,x)\in\mathbf{R}_{+}\times\mathbf{R},
\]
is solution of the partial differential equation
\[
\left\{
\begin{array}{ll}
\partial_{t}u(t,x)=\frac{1}{2}\Delta u(t,x)+xu(t,x)& \quad \forall\,x\in \mathbf{R},\ t\in\mathbf{R}_{+},\\
u(0,x)=f(x)&\quad \forall\,x\in\mathbf{R}.
\end{array}\right.
\]
In this case, taking the Laplace transform in time of $u$ leads to an ordinary differential equation whose solution can be expressed in terms of Airy functions \cite{Ja07}. In our problem, this partial differential equation becomes inhomogeneous in time which prevents us to use this Laplace transform. As a consequence, we think that one can not obtain information by this method. This is why we have established that $\Lambda_\infty$ admits a density using Malliavin calculus and the representation of the Brownian excursion as a three-dimensional Bessel bridge \eqref{eq:rep:3d}.
\end{remark}

\begin{figure}[!h]
\centering
\includegraphics[scale=0.55]{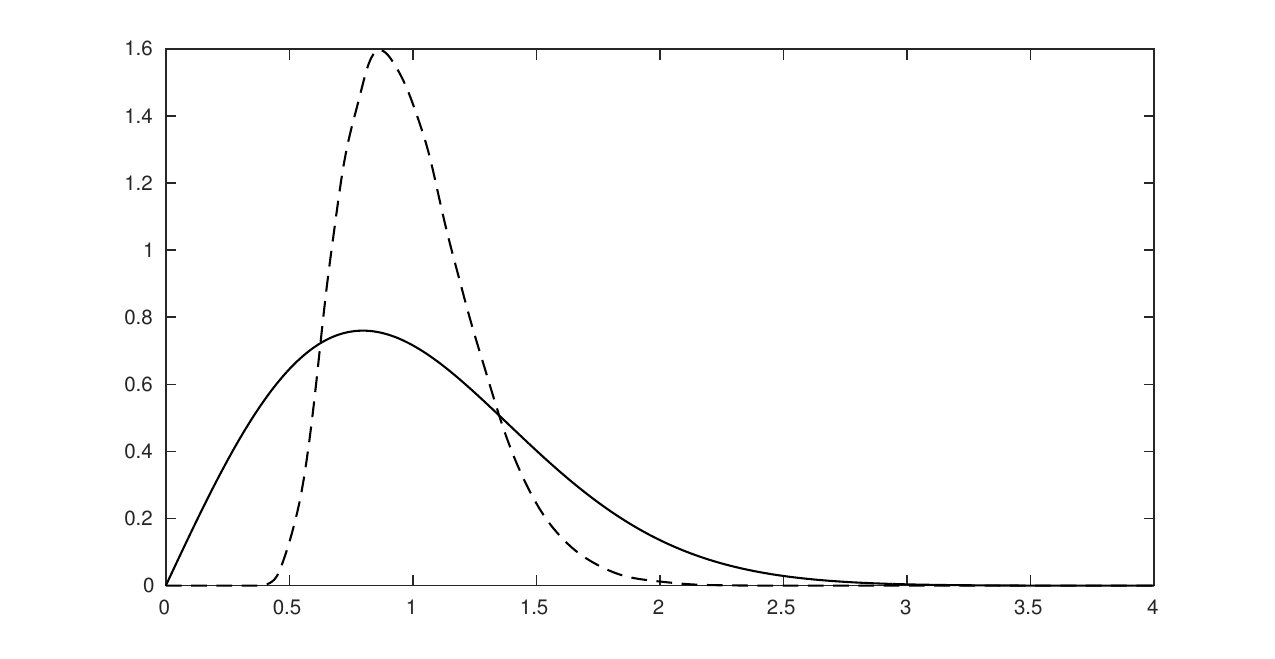}
\caption{Densities of $\sqrt{\frac{2}{\pi}}\Delta_{\infty}$ (full line) where $\Delta_\infty$ follows the Rayleigh distribution given by $f(x)=\frac{\pi}{2} x \exp\left(-\frac{\pi x^2}{4}\right)$ for $x\in\R_+$ and of $\Lambda_{\infty}$ (dashed line) estimated from $1\,000\,000$ simulated Brownian excursions.}
\label{fig:densities}
\end{figure}

Of course, $\widehat{\lambda}[\tau_n]$ is not a consistent estimator of $\sigma^{-1}$ but it should be noted that its weak limit is unbiased by \eqref{eq:variance:lambdainfty} and Corollary \ref{lim:lambda:n}.
The expression \eqref{eq:variance:lambdainfty} of the variance of $\Lambda_\infty$ is an explicit but quite intractable formula. Nevertheless, it may be at least evaluated numerically to compute the variance of $\Lambda_\infty$. Otherwise, we can also use Monte Carlo simulations to produce a sample with the same law as $\Lambda_\infty$ to achieve this task. Both methods lead to
$$
\Var(\Lambda_\infty) \simeq 0.0690785.
$$
At this point, it is quite interesting to compare our approach to the one developed in \cite{KPDRV14}. They construct estimators for the variance of the offspring distribution of a forest of conditioned critical Galton-Watson trees. Their strategy relies on the distance to the root of a uniformly sampled node $v$ of the considered tree $\tau_n\sim\text{GW}_n(\mu)$,
\begin{equation}\label{eq:delta1}
\delta[\tau_n]=\frac{h(v)}{\sqrt{n}} .
\end{equation}
Using Theorem \ref{thm:aldous}, it has been shown that $\delta[\tau_n]$ converges in law, when the number of nodes $n$ goes to infinity, towards $\sigma^{-1}\Delta_\infty$ where the random variable $\Delta_\infty$ follows the Rayleigh distribution with scale $1$ \cite[Proposition 4]{KPDRV14} with density,
$$\forall\,x\in\R_+,~f_{\Delta_\infty}(x)=x\exp\left(-\frac12 x^2\right).$$
We emphasize that $\delta[\tau_n]$ is somehow biased because $\E[\Delta_\infty]=\sqrt{\frac{\pi}{2}}\neq 1$. Nevertheless, one may avoid this issue by considering the quantity
\begin{equation}\label{eq:delta2}\widehat{\delta}[\tau_n]=\sqrt{\frac{2}{\pi}}\delta[\tau_n]\end{equation}
that converges to $\sigma^{-1}\sqrt{\frac{2}{\pi}}\Delta_\infty$ which is $\sigma^{-1}$ on average. As a consequence, $\widehat{\lambda}[\tau_n]$ and $\widehat{\delta}[\tau_n]$ are two quantities directly computable from the tree $\tau_n$ and that may be used to construct an estimator of the inverse standard deviation of interest. We propose to compare them from their respective asymptotic dispersion which should be as small as possible in order to get an accurate estimator. A first comparison may be done by computing the variances of $\Lambda_\infty$ and $\sqrt{\frac{2}{\pi}}\Delta_\infty$. One has
$$\Var\left(\sqrt{\frac{2}{\pi}}\Delta_\infty\right) \simeq 0.2732395 \qquad\text{and}\qquad \Var(\Lambda_\infty) \simeq 0.0690785 .$$
This difference in the dispersions is quite apparent in Figure \ref{fig:densities} where the densities of $\sqrt{\frac{2}{\pi}}\Delta_\infty$ and $\Lambda_\infty$ have been displayed. Consequently, one may expect better results in terms of dispersion from our approach.

\subsection{Interpretation in functional principal component analysis}
\label{ss:FPCA}
\begin{figure}[!h]
\centering
\includegraphics[scale=0.45]{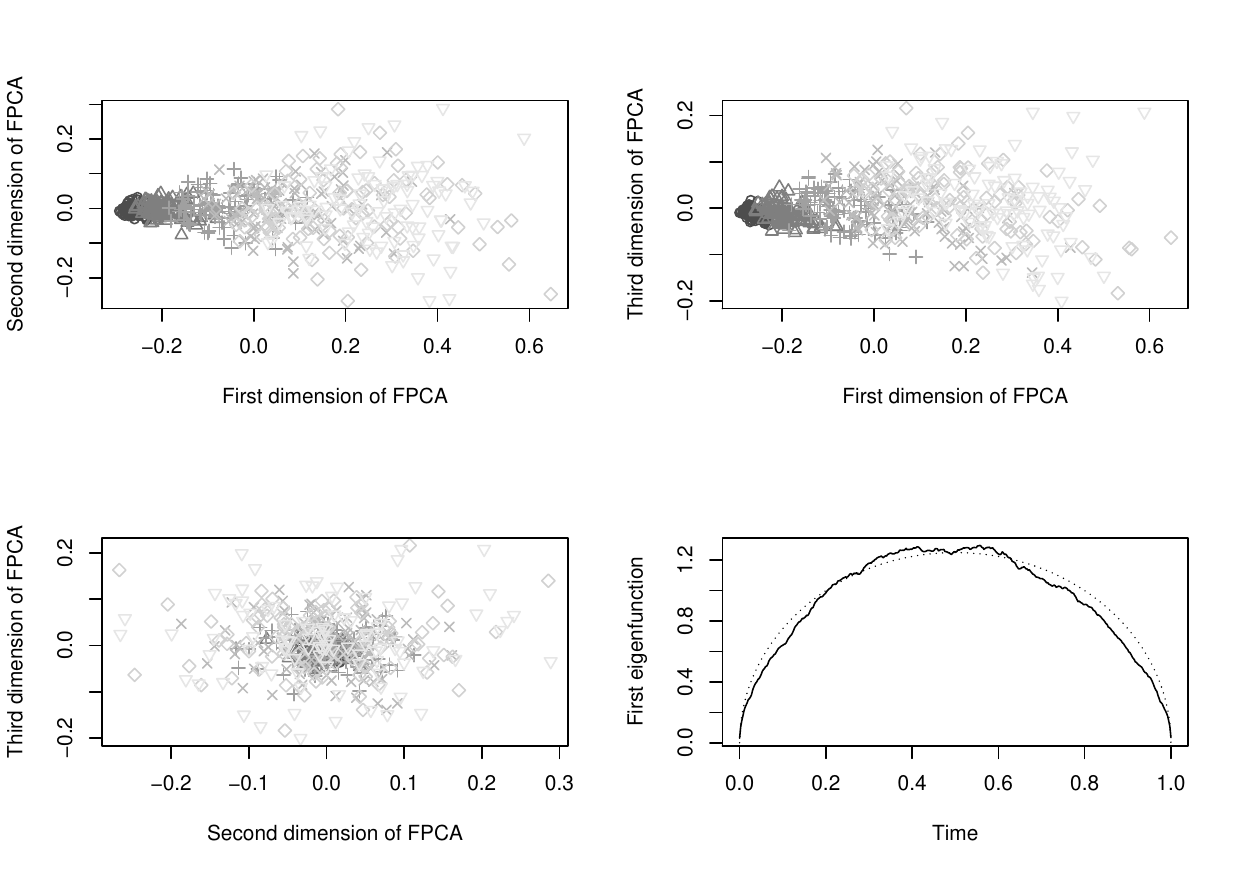}
\caption{FPCA performed from $600$ large conditioned Galton-Watson trees with different standard deviations: from dark to bright gray, $\sigma\in\{0.1,0.2,0.4,0.6,0.8,0.95\}$. Coefficients of projection on the three first eigenfunctions (top left, top right and bottom left) and first eigenfunction (full line, bottom right) compared to the average Brownian excursion (dashed line).}
\label{fig:fpca}
\end{figure}

On the suggestion of a reviewer, we performed a functional principal component analysis (FPCA) from normalized Harris paths of large conditioned Galton-Watson trees with different values of $\sigma$. We refer the reader to \cite{FPCA} and the references therein for explanations on this statistical tool and to Subsection \ref{ss:simugw} for simulation methods. The FPCA has been carried out with function \verb+FPCA+ from the \verb+R+ package \verb+fdapace+. The results are presented in Figure~\ref{fig:fpca}.

These numerical experiments show that the first eigenfunction is very close to the average Brownian excursion (see Figure~\ref{fig:fpca} (bottom right)). As a consequence, $\widehat{\lambda}[\tau]$ can be interpreted as the first eigenvalue in FPCA of the normalized Harris path of $\tau$. The first eigenspace expresses approximately $70\%$ of the total dataset inertia. For the sake of comparison, the second eigenspace gets only $7\%$ of the inertia. In addition, one can see on Figure~\ref{fig:fpca} (top left, top right and bottom left) that only the first dimension of FPCA captures information on the value of $\sigma$. These results highlight that $\widehat{\lambda}[\tau]$ is a relevant quantity in our estimation problem as well as one can not expect more significant information from the projection on the other eigenspaces.

Alternative functionals of the Brownian excursion can be investigated in order to develop statistical methods for conditioned Galton-Watson trees. For instance, the vector of peaks and valleys at random \cite{Derrida2004} or fixed \cite{pitman1999} times could be considered. If the distribution of such objects is simple enough, this could enable the use of maximum likelihood methods. However, as shown above, the functional $\widehat{\lambda}[\tau]$ considered in our work seems to be one of the best in our setting in terms of quantity of information.

\subsection{Estimation strategies}

\noindent
In this section, we give details on two ideas in order to estimate $\sigma^{-1}$ from a forest of conditioned Galton-Watson trees. A forest is defined as a tuple of trees. Let $N$ be a positive integer. In this section, we consider a forest $\mathcal{F}$ made of $N$ independent trees $\tau^{1},\dots,\tau^{N}$ with respective sizes $n_{1},\dots,n_{N}$ and respective laws $\text{GW}_{n_{1}}(\mu),\dots,\text{GW}_{n_{N}}(\mu)$.

\subsubsection{Least square estimation}

\noindent
This first strategy lies on the goodness of fit between the Harris path of the forest with the expected limiting contour. This adequacy is measured thanks to an $\mathbf{L}^{2}([0,N],\mathbf{R})$-norm.
More precisely, we denote $(\mathcal{H}[\mathcal{F}](t),\ t\in[0,N])$ the Harris path of the forest $\mathcal{F}$. This process is defined by
$$\forall\,0\leq t\leq N,~\mathcal{H}[\mathcal{F}](t) = \sum_{i=1}^N \frac{1}{\sqrt{n_{i}}}\mathcal{H}[\tau^i](2n_i (t-i+1)) \mathbf{1}_{[i-1,i)}(t).$$
The Harris path of a forest is simply the concatenation of the Harris paths of each tree, in the natural order. 
We propose to estimate $\sigma^{-1}$ by $\widehat{\lambda}_{ls}[\mathcal{F}]$ that minimizes the $\mathbf{L}^2([0,N],\mathbf{R})$-error
$$\lambda\mapsto\|\mathcal{H}[\mathcal{F}](\cdot) - \lambda H(\cdot-\lfloor\cdot\rfloor)\|_{\mathbf{L}^2([0,N],\mathbf{R})}^2 ,$$
the function $H(\cdot-\lfloor\cdot\rfloor)$ mapping $x\in[0,N]$ to $H(x-\lfloor x\rfloor)$. As aforementioned in \eqref{eq:def:lambdahat}, $\widehat{\lambda}_{ls}[\mathcal{F}]$ can be explicitly computed. Indeed, one can check that
\begin{equation*}\label{eq:def:lambda:ls}
\widehat{\lambda}_{ls}[\mathcal{F}] = \frac{\langle \mathcal{H}[\mathcal{F}](\cdot) , H(\cdot-\lfloor\cdot\rfloor)\rangle}{\|H(\cdot-\lfloor\cdot\rfloor)\|_{{2}}^2}.
\end{equation*}
We remark that $\widehat{\lambda}_{ls}[\mathcal{F}]$ is only the average of the quantities $\widehat{\lambda}[\tau^i]$ (defined in \eqref{eq:def:lambdahat}),
$$\widehat{\lambda}_{ls}[\mathcal{F}]  = \frac{1}{N} \sum_{i=1}^N \widehat{\lambda}[\tau^i].$$
Thus, according to Theorems \ref{thm:cv:expectation} and \ref{lim:lambda:n}, one can expect that $\widehat{\lambda}_{ls}[\mathcal{F}]$ tends to $\sigma^{-1}$ in some sense, when both $N$ and $n_i$ go to infinity, by virtue of the law of large numbers.

\subsubsection{Estimation by minimal Wasserstein distance}

In the preceding method, we did not use our knowledge of the limiting distribution of the random variable of type $\lambda[\tau^{n}]$. In order to take this into account, one may want to test the goodness of fit between the empirical measure $ \widehat{\mathcal{P}}$ defined by
\begin{equation}
\label{eq:empiricalMeasure}
\widehat{\mathcal{P}}=\frac{1}{N}\sum_{i=1}^N \delta_{\widehat{\lambda}[\tau^i]}
\end{equation}
and the the law of $\Lambda_{\infty}$.
 Using Wasserstein metrics to align distributions is rather natural since it corresponds to the transportation cost between two probability laws. In particular, this feature appears to be useful in a statistical framework \cite{CzadoMunk,gallon:hal-01135106}. In our case, $\widehat{\mathcal{P}}$ is expected to be close in terms of Wasserstein distance to $\sigma^{-1}\Lambda_{\infty}$ in the asymptotic regime of an infinite forest of infinite trees. That is why, we propose to estimate $\sigma^{-1}$ with the real number $\lambda$ which minimizes the distance between $\widehat{\mathcal{P}}$ and $\lambda\Lambda_{\infty}$. More precisely, our estimator $\widehat{\lambda}_W[\mathcal{F}]$ is defined by
\begin{equation*}
\label{eq:def:lambdaW}
\widehat{\lambda}_W[\mathcal{F}] = \argmin_{\lambda>0} d_W\left(\widehat{\mathcal{P}} \, , \, \mathbf{P}_{\lambda\Lambda_{\infty}}\right) ,
\end{equation*}
where $d_W$ denotes the Wasserstein distance of order $2$ and $\mathbf{P}_{\lambda\Lambda_{\infty}}$ denotes the law of $\lambda\Lambda_{\infty}$.

The Wasserstein distance of order $2$, denoted $d_W(\nu_1,\nu_2)$, between two probability measures $\nu_1$ and $\nu_2$ on $\mathbf{R}$ can be defined from their cumulative distribution functions $F_1$ and $F_2$ as follows,
 \begin{equation}
 \label{eq:wassersteinCumulative}
 d_W(\nu_1,\nu_2)=\|F^{-1}_1-F^{-1}_2\|_{2}.
 \end{equation}
 Let $\widehat{F}$ be the cumulative function of the empirical measure $\widehat{\mathcal{P}}$, while $F_{\lambda\Lambda_\infty}$ stands for the cumulative function of the random variable $\lambda\Lambda_\infty$. As a consequence of \eqref{eq:wassersteinCumulative}, one has
 \begin{eqnarray*}
 d_W\left(\frac{1}{N}\sum_{i=1}^N \delta_{\widehat{\lambda}[\tau^i]} \, , \, \mathbf{P}_{\lambda\Lambda_{\infty}}\right)^2&=&\int_0^1\left(\widehat{F}^{-1}(s)-F_{\lambda\Lambda_\infty}^{-1}(s) \right)^{2}\dd s\\
 &=&\int_0^1\left( \widehat{F}^{-1}(s)-\lambda F^{-1}_{\Lambda_\infty}(s) \right)^{2}\dd s,
 \end{eqnarray*}
 thanks to the fact that $F^{-1}_{\lambda\Lambda_\infty}=\lambda F^{-1}_{\Lambda_\infty}$. It follows that minimizing the Wasserstein distance boils down to solving a least square minimization problem. Hence, it comes that
 \begin{eqnarray}
 \widehat{\lambda}_W[\mathcal{F}] &=&\frac{\langle\widehat{F}^{-1}, F^{-1}_{\Lambda_\infty}\rangle}{\|F^{-1}_{\Lambda_\infty}\|^{2}_{{2}}} \nonumber\\
 &=&\frac{1}{\|F^{-1}_{\Lambda_\infty}\|^{2}_{{2}}}\sum_{i=1}^{N}\widehat{\lambda}[\tau^{(i)}]\int_{\frac{i-1}{N}}^{\frac{i}{N}}F^{-1}_{\Lambda_\infty}(s)\dd s, \label{eq:wasserstein:formula}
 \end{eqnarray}
 where $(\widehat{\lambda}[\tau^{(i)}])_{1\leq i\leq N}$ denotes the order statistic associated with the family $(\widehat{\lambda}[\tau^{i}])_{1\leq i\leq N}$. 

\begin{remark}
\label{rem:lastRem}
We point out the fact that there is no problem of definition in the above quantities because both $\widehat{F}^{-1}$ and $F^{-1}_{\Lambda_{\infty}}$ belong to $\mathbf{L}^{2}$. In the first case, this follows from the fact that $\widehat{F}^{-1}$ is bounded (because $\widehat{\mathcal{P}}$ has compact support). For $F^{-1}_{\Lambda_{\infty}}$, this comes from the uniform sampling principle which entails that
\[
\int_{0}^{1}F^{-1}_{\Lambda_{\infty}}(u)^{2}\ \dd u=\mathbf{E}[\Lambda_{\infty}^{2}].
\]
\end{remark}

\begin{remark}
The proposed methodology consists in identifying the best parameter $\lambda$ that allows to align the distributions $\widehat{\mathcal{P}}$ and $\mathbf{P}_{\lambda\Lambda_{\infty}}$. The Wasserstein distance is well-adapted to this problem because it is computed from the inverse cumulative distribution functions together with the fact that $F^{-1}_{\lambda\Lambda_\infty}=\lambda F^{-1}_{\Lambda_\infty}$. As a consequence, one may get the optimal parameter $\widehat{\lambda}_W[\mathcal{F}]$ from only a numerical estimate of $F^{-1}_{\Lambda_\infty}$. The same trick does not hold for the maximum likelihood method: one can not express the likelihood of $\lambda\Lambda_\infty$ as a function of the two variables $\lambda$ and $f_{\Lambda_\infty}$. Thus this alternative method is not adequate without  an explicit formula for $f_{\Lambda_\infty}$, which seems to be out of our reach.
\end{remark}


\section{Main results}
\label{sec:mainResults}

\subsection{Increasing sequences of random forests}
\label{ss:statfram}

Before going further, the statistical framework needs to be precisely formulated. In the sequel, the set of integer sequences is denoted by $\mathbf{S}$. For any positive real number $A$, we denote by $\mathbf{S}_{A}$ the subset of $\mathbf{S}$ defined by
\[
\mathbf{S}_{A}=\left\{u\in\mathbf{S}\, :\, \min_{i\geq 1}u_{i}\geq A \right\}.
\]
In addition, for any sequence $u$ in $\mathbf{S}$ and any positive integer $N$, $\vec{u}_{N}$ is the multi-integer made of the $N$ first components of $u$, that is
\[
\vec{u}_{N}=\left(u_{1},\dots, u_{N} \right).
\]

Now, let us introduce our probabilistic framework.
Let $(\tau_{n}^{k})_{n,k\geq 1}$ be a family of independent conditioned Galton-Watson trees such that, for a given $n$, the family $(\tau_{n}^{k})_{k\geq 1}$ is i.i.d.\ $\text{GW}_{n}(\mu)$. From this family, we define, for any multi-integer $\vec{u}_N=(u_1,\dots,u_N)$, the random forest $\mathcal{F}_{\vec{u}_N}$ made of the trees $(\tau_{u_{1}}^{1},\dots,\tau_{u_{N}}^{N})$. 

The idea of this construction is to consider increasing (in the sense of inclusion) sequences of random forests. Indeed, assume we are given a sequence $(u_{n})_{n\geq 1}$ of integers (corresponding with the sizes of our trees), then the $N$ first trees of the forest $\mathcal{F}_{\vec{u}_{N+1}}$ are the same as the trees of the forest $\mathcal{F}_{\vec{u}_{N}}$.

We point out that the hypothesis of independence may be thought to be too strong in some applications. Exchangeability is a weaker assumption that could be considered. In such a statistical setting, the reference \cite{Haulk} is particularly relevant.

\subsection{Least square estimation}
\label{ss:result:ls}

This first result focuses on the large trees regime and gives the asymptotic unbiasedness of the least square estimator in this regime.

\begin{proposition}
\label{cor:ls:law}
The least square estimator is asymptotically unbiased in the large trees regime, that is
\[
\forall\,\epsilon>0,~\exists\,A\in\N,~\forall\,u\in\mathbf{S}_{A},~\forall\,N\in\mathbf{N}~\left|\mathbf{E}\left[\widehat{\lambda}_{ls}[\mathcal{F}_{\vec{u}_{N}}]\right]-\sigma^{-1}\right|<\epsilon.
\]
This means that the expectation of the least square estimator converges to $\sigma^{-1}$ as the sizes of the trees increase.
\end{proposition}
\textit{Proof.} Since the family $(\tau_{u_i}^i)_{1\leq i\leq N}$ is made of independent random variables, its follows from Theorem \ref{thm:cv:expectation} and the definition \eqref{eq:def:lambdahat} of $\widehat{\lambda}[\tau_{{\bf}n_{i}}]$ that the proof of this last statement boils down to proving that, when $n$ goes to infinity,
\[
\int_{0}^{1}\E\left[\frac{\mathcal{H}[\tau_n](2n s)}{\sqrt{n}}\right]E_{s}\ \dd s\longrightarrow\frac{2}{\sigma}\int_{0}^{1}E^{2}_s \ \dd s, 
\]
where $\tau_{n}$ is some tree with law $\text{GW}_{n}(\mu)$. It is known from \cite[Lemma 4]{Marck05} that, for any positive integer $n$ and real number $0< t<1$,
\begin{equation}
\label{eq: niceUnifEstimate}
 \forall\,x\in\mathbf{R}_{+},~\mathbf{P}\left( \frac{\mathcal{H}[\tau_{n}](2n t)}{\sqrt{n}} > x \right)\leq \frac{C}{t}\exp\left(\frac{-Dx}{\sqrt{t}}\right).
\end{equation}
From this last estimate, one can easily show that $\E\left[\frac{\mathcal{H}[\tau_n](2n \cdot)}{\sqrt{n}}\right]$ is uniformly bounded (w.r.t.\ $n$) by an integrable function. Finally, the result follows from Theorem \ref{thm:cv:expectation} and the dominated convergence theorem.
\hfill$\Box$

The spirit of the following result is that, given an increasing sequence of random forests, the least square estimator can not be too far from $\sigma^{-1}$ as soon as the sizes of the trees are large enough. In particular, due to the weakness of the convergence of conditioned Galton-Watson trees given in Theorem \ref{thm:cv:expectation}, one can not expect a stronger result of convergence.

\begin{proposition}We have,
\label{prop:CVLS1}
\[
\forall\,\epsilon>0,~\exists\,A\in\N,~\forall\, u\in\mathbf{S}_{A},~\mathbf{P}\left(\limsup_{N\to\infty}\left|\widehat{\lambda}_{ls}[\mathcal{F}_{\vec{u}_{N}}]-\sigma^{-1}\right|<\epsilon\right)=1.
\]
\end{proposition}

\textit{Proof.}
We begin the proof by showing that the family $(\widehat{\lambda}[\tau_{n}^{k}])_{n,k\geq 1}$ has uniformly bounded fourth moments.
By Jensen's inequality, there exists a positive constant $c$ such that
\begin{eqnarray}
\label{eq:newEqRef}
 \mathbf{E}\left[\left(\widehat{\lambda}[\tau^{i}_{n}]\right)^{4}\right]&\leq & c \int_0^1\mathbf{E}\left[\left(\frac{\mathcal{H}[\tau^{i}_{n}](2n s)}{\sqrt{n}} \right)^{4} \right]\dd s\notag\\
 &=&4c\int_0^1\int_{\mathbf{R}_{+}}x^{3}\ \mathbf{P}\left( \frac{\mathcal{H}[\tau^{i}_{n}](2n s)}{\sqrt{n}} > x \right)\dd x\ \dd s.
 \end{eqnarray}
Finally, using again equation \eqref{eq: niceUnifEstimate} gives the desired bound,
\begin{equation}\label{unif:bound}
 \mathbf{E}\left[\widehat{\lambda}[\tau_{n}^i]^{4}\right]\leq \frac{12\,c\,C}{D^{4}}.
\end{equation}
From this point we consider a sequence $u$ of integers. This sequence corresponds to the sizes of the trees in our increasing sequence of random forests $(\mathcal{F}_{\vec{u}_{N}})_{N\geq 1}$. We recall according to the definitions given in the beginning of this section that the random forest $\mathcal{F}_{\vec{u}_{N}}$ is composed of the trees $(\tau_{u_{1}}^{1},\dots,\tau_{u_{N}}^{N})$.

Let $m_{u_i}^i$ be the expectation of $\widehat{\lambda}[\tau_{u_i}^i]$. It is worth noting that this expectation depends only on the integer $u_{i}$. Now, using the uniform bound on the fourth moment \eqref{unif:bound}, Markov's inequality applied to the fourth power of
	\[
	\frac{1}{N}\sum_{i=1}^{N}\left(\widehat{\lambda}[\tau_{u_{i}}^{i}]-m_{u_{i}}^i\right)
	\]
gives the convergence in probability of the above sum to zero at rate $N^{-2}$ which implies, in light of the Borel-Cantelli lemma, that
\begin{equation}\label{bh:lgn}
\frac{1}{N}\sum_{i=1}^{N}\left(\widehat{\lambda}[\tau_{u_{i}}^{i}]-m_{u_{i}}^i\right)\xrightarrow{a.s.} 0 ,
\end{equation}
when $N$ goes to infinity. Moreover, using Theorem \ref{thm:cv:expectation}, we have that $m_{u_i}^i$ converges to $\sigma^{-1}$ as $u_i$ goes to infinity, from which it follows that for any $\epsilon>0$, there exists an integer $A$ such that 
\begin{equation}\label{bh:small}
\left|m_{u_i}^i-\sigma^{-1}\right|<\epsilon,
\end{equation}
whenever $u_i>A$. Finally, letting all the $u_i$'s be greater than $A$, we have that there exists a measurable set $\Omega_{u}$ 
, with mass $1$, such that, using (\ref{bh:lgn}) and (\ref{bh:small}), for all $\omega$ in this set,
\begin{align*}
& \limsup_{N\to\infty}\left|\frac{1}{N}\sum_{i=1}^{N} \widehat{\lambda}[\tau_{u_i}^i](\omega)-\sigma^{-1}\right| \\
& \leq\,\limsup_{N\to\infty}\frac{1}{N}\left|\sum_{i=1}^{N}  \widehat{\lambda}[\tau_{u_i}^i](\omega)-m_{u_i}^i \right|+\limsup_{N\to\infty}\frac{1}{N}\sum_{i=1}^{N} \left|m_{u_i}^i-\sigma^{-1}\right| \,\leq\,\epsilon,
\end{align*}
which establishes the expected convergence.\hfill$\Box$

\begin{remark}
According to the proof of the preceding theorem, it would be very interesting to control the rate of convergence in Theorem \ref{thm:cv:expectation}. Indeed, this would enable us to get a control of the error in the convergence stated in Proposition \ref{prop:CVLS1} given in terms of the smallest tree in the increasing sequence of random forests. 
\end{remark}

\begin{remark}
\label{rem:lambda:moment}
Let us point out that equation \eqref{eq: niceUnifEstimate} gives the exponential decay of the tail distribution of $n^{-1/2}\mathcal{H}[\tau_{n}](2n t)$ uniformly w.r.t.\ $n$. In particular, one can apply the method used in equation \eqref{eq:newEqRef} to obtain uniform (w.r.t.\ $n$) bounds like
\[
\mathbf{E}\left[\widehat{\lambda}[\tau_{n}]^{k} \right]\leq \frac{(k-1)!\,c\,C}{D^{k}},
\] 
for any positive integer $k$ and some positive constant $c$.
\end{remark}

\subsection{Estimation by minimal Wasserstein distance}
\label{ss:result:w}

As in the preceding section, we begin by looking at the asymptotic bias of the considered estimator.

\begin{proposition}\label{prop:w:law}
The Wasserstein estimator is asymptotically unbiased in the large trees and large forests regime. That is,
\[
\forall\,\epsilon>0,~\exists\,(\mathfrak{N},A)\in\mathbf{N}^{2},~\forall\,u\in\mathbf{S}_{A},~\forall\,N\geq \mathfrak{N},~ \left|\mathbf{E}\left[\widehat{\lambda}_W[\mathcal{F}_{\vec{u}_{N}}]\right]-\sigma^{-1} \right|<\epsilon.
\]
\end{proposition}
\textit{Proof.} Let $u$ in $\mathbf{S}$ and $N$ in $\mathbf{N}$. Let $(\Lambda_{\infty,i})_{i\geq 1}$ be an i.i.d.\ sequence of random variables with the same distribution as $\Lambda_{\infty}$.  The first step of the proof is to show that 
\[
\eta_{\vec{u}_N} = \left|\mathbf{E}\left[\widehat{\lambda}_W[\mathcal{F}_{\vec{u}_{N}}]-\frac{1}{\sigma \|F^{-1}_{\Lambda_\infty}\|_2^{2}}\sum_{i=1}^{N}\mathbf{E}\left[\Lambda_{\infty,(i)}\right] \int_{\frac{i-1}{N}}^{\frac{i}{N}}F^{-1}_{\Lambda_\infty}(s)\dd s\right]\right|
\]
converges to $0$ as $\min(\vec{u}_{N})$ goes to infinity uniformly w.r.t.\ $N$, where the order statistic in the above formula has to be understood w.r.t.\ the random vector $\left(\Lambda_{\infty,1},\dots,\Lambda_{\infty,N} \right)$.
The convergence easily follows from the results of \cite{Marck05}. However, getting it uniformly w.r.t.\ $N$ requires to give some new insights. To this end let $(\mathcal{U}_i)_{i\geq 1}$ be a sequence of i.i.d.\ random variables with uniform distribution on $[0,1]$. Denote, for any positive integer $n$, $F_{n}^{-1}$ the right inverse of the cumulative distribution function associated to the random variable $\widehat{\lambda}[\tau_{n}^{1}]$. According to the right inverse principle, the random vector $(F^{-1}_{u_{i}}(\mathcal{U}_{i}))_{1\leq i\leq N}$ is equal in distribution to $(\widehat{\lambda}[\tau^{i}_{u_{i}}])_{1\leq i\leq N}$ and the vector $(F^{-1}_{\Lambda_{\infty}}(\mathcal{U}_{i}))_{1\leq i\leq N}$ is a vector of i.i.d.\ random variables with the same distribution as $\Lambda_{\infty}$. At this point, let us highlight that the sequences $(\widehat{\lambda}[\tau^{i}_{n}])_{n\geq 1}$ and $(F^{-1}_{n}(\mathcal{U}_{i}))_{n \geq 1}$ do not have the same distribution since the first one is made of independent random variables whereas this is clearly not the case for the second one. However, this feature is not important since, in the above formula, we only look at averaged behaviors. 
First note that, by definition of the Wasserstein estimator and the right inverse sampling principle, we have
\begin{equation}
\label{eq:lasteq}
\eta_{\vec{u}_N}=\frac{1}{\sigma \|F^{-1}_{\Lambda_\infty}\|_2^{2}}\left|\int_{0}^{1}F^{-1}_{\Lambda_{\infty}}(s)\mathbf{E}\left[\widehat{G}^{-1}_{\vec{u}_{N}}(s)-\widehat{H}^{-1}_{N}(s) \right]\dd s \right|,
\end{equation}
where $\widehat{G}^{-1}_{\vec{u}_{N}}$ and $\widehat{H}^{-1}_{N}$ denote the inverse distribution functions of the empirical measures respectively associated to the vectors $(F^{-1}_{u_{i}}(\mathcal{U}_{i}))_{1\leq i\leq N}$ and $(F^{-1}_{\Lambda_{\infty}}(\mathcal{U}_{i}))_{1\leq i\leq N}$.
Now, the Cauchy-Schwartz inequality entails that
\[
\left|\int_{0}^{1}F^{-1}_{\Lambda_{\infty}}(s)\mathbf{E}\left[\widehat{G}^{-1}_{\vec{u}_{N}}(s)-\widehat{H}^{-1}_{N}(s) \right]\dd s \right|\leq \|F^{-1}_{\Lambda_{\infty}}\|_{2}\sqrt{\int_{0}^{1} \mathbf{E}\left[\widehat{G}^{-1}_{\vec{u}_{N}}(s)-\widehat{H}^{-1}_{N}(s) \right]^{2}\dd s}.
\]
By the definition of the inverse distribution function, we get
\begin{align*}
& \int_{0}^{1} \mathbf{E}\left[\widehat{G}^{-1}_{\vec{u}_{N}}(s)-\widehat{H}^{-1}_{N}(s) \right]^{2}\dd s\\
& =\int_{0}^{1}\left(\sum_{i=1}^{N} \mathbf{1}_{\left[\frac{i-1}{N},\frac{i}{N}\right)}(s)\mathbf{E}\left[F^{-1}_{u_{(i)}}(\mathcal{U}_{(i)})-F^{-1}_{\Lambda_{\infty}}(\mathcal{U}_{(i)})\right]\right)^{2}\dd s,
\end{align*}
where $(F^{-1}_{n_{(i)}}(\mathcal{U}_{(i)}))_{1\leq i \leq N}$ and $(F^{-1}_{\Lambda_{\infty}}(\mathcal{U}_{(i)}))_{1\leq i \leq N}$ denote the order statistics respectively associated with the vectors $(F^{-1}_{n_{i}}(\mathcal{U}_{i}))_{1\leq i\leq N}$ and $(F^{-1}_{\Lambda_{\infty}}(\mathcal{U}_{i}))_{1\leq i\leq N}$.
Using two times Jensen's inequality leads to
\[
\int_{0}^{1} \mathbf{E}\left[\widehat{G}^{-1}_{\vec{u}_{N}}(s)-\widehat{H}^{-1}_{N}(s) \right]^{2}\dd s\leq \frac{1}{N}\sum_{i=1}^{N}\mathbf{E}\left[\left|F^{-1}_{n_{(i)}}(\mathcal{U}_{(i)})-F^{-1}_{\Lambda_{\infty}}(\mathcal{U}_{(i)})\right|^{2} \right].
\]
Finally, using that the order function $(x_{1},\dots,x_{n})\mapsto(x_{(1)},\dots,x_{(n)})$ is 1-Lipschitz w.r.t.\ the Euclidean norm (as a consequence of the rearrangement inequality), we have
\begin{align}
&\left|\int_{0}^{1}F^{-1}_{\Lambda_{\infty}}(s)\mathbf{E}\left[\widehat{G}^{-1}_{\vec{u}_{N}}(s)-\widehat{H}^{-1}_{N}(s) \right]\dd s \right| \nonumber \\
&\leq \|F^{-1}_{\Lambda_{\infty}}\|_{2}\sqrt{\frac{1}{N}\sum_{i=1}^{N}\mathbf{E}\left[\left|F^{-1}_{n_{i}}(\mathcal{U}_{i})-F^{-1}_{\Lambda_{\infty}}(\mathcal{U}_{i})\right|^{2} \right]}.\label{eq:lasteq-1}
\end{align}
Now by construction, for any $i$, $F^{-1}_{n}(\mathcal{U}_{i})$ converges almost surely to $F^{-1}_{\Lambda_{\infty}}(\mathcal{U}_{i})$. Moreover, the uniform square-integrability of the laws of the $\widehat{\lambda}[\tau^i_{n_i}]$'s provided by Remark \ref{rem:lambda:moment} gives that
\[
\lim_{n\to\infty}\mathbf{E}\left[\left|F^{-1}_{n}(\mathcal{U}_{i})-F^{-1}_{\Lambda_{\infty}}(\mathcal{U}_{i})\right|^{2} \right]=0.
\]
Now, let $0<\epsilon<1$ be some positive real number and  let $A$ such that for any $n\geq A$, we have
\[
\mathbf{E}\left[\left|F^{-1}_{n}(\mathcal{U}_{i})-F^{-1}_{\Lambda_{\infty}}(\mathcal{U}_{i})\right|^{2} \right]\leq \sigma^{2}\|F^{-1}_{\Lambda_{\infty}}\|^{2}_{2}\,\epsilon^{2}.
\]
Hence, as soon as $\min(u)\geq A$, we have, together with \eqref{eq:lasteq-1},
\[
\frac{1}{\sigma \|F^{-1}_{\Lambda_{\infty}}\|_{2}^{2}}\left|\int_{0}^{1}F^{-1}_{\Lambda_{\infty}}(s)\mathbf{E}\left[\widehat{G}^{-1}_{\vec{u}_{N}}(s)-\widehat{H}^{-1}_{N}(s) \right]\dd s \right|\leq\epsilon.
\]
Finally, \eqref{eq:lasteq} gives the desired uniform convergence. It remains to prove that
\[
\frac{1}{\sigma \|F^{-1}_{\Lambda_\infty}\|_2^{2}}\sum_{i=1}^{N}\mathbf{E}\left[\Lambda_{\infty,(i)}\right] \int_{\frac{i-1}{N}}^{\frac{i}{N}}F^{-1}_{\Lambda_\infty}(s)\dd s
\]
converges to $\sigma^{-1}$ as $N$ goes to infinity. It is well known, since $\Lambda_{\infty}$ has a density, that, for any $1\leq i\leq N$, one has (see for instance \cite{DN03})
$$\E\left[\Lambda_{\infty,(i)}\right] = N\dbinom{N-1}{i-1}\int_0^\infty x F_{\Lambda_\infty}(x)^{i-1} (1-F_{\Lambda_\infty}(x))^{N-i} f_{\Lambda_\infty}(x) \dd x.$$
Hence, 
\begin{align*}
&\E\left[\sum_{i=1}^{N}\Lambda_{\infty,(i)} \int_{\frac{i-1}{N}}^{\frac{i}{N}}F^{-1}_{\Lambda_\infty}(s)\dd s\right]\\
&= N \int_{0}^{\infty}\!\!\!\!xf_{\Lambda_{\infty}}(x)\sum_{i=1}^{N}\binom{N-1}{i-1}F_{\Lambda_\infty}(x)^{i-1} (1-F_{\Lambda_\infty}(x))^{N-i}\\
& \times \int_{0}^{\frac{1}{N}} F^{-1}_{\Lambda_\infty}\left(s+\frac{i-1}{N}\right)\dd s\ \dd x.
\end{align*}
This rewrites thanks to the right inverse sampling principle as
\[
\E\left[\sum_{i=1}^{N}\Lambda_{\infty,(i)} \int_{\frac{i-1}{N}}^{\frac{i}{N}}F^{-1}_{\Lambda_\infty}(s)\dd s\right]=\int_{0}^{1}\!\!\!\iF(y) K_{N}\left(\iF\right)(y)\ \dd y,
\]
where $K_{N}$ is defined for all function $\varphi$ in $\mathbf{L}^{2}$ by
\[
K_{N}\left(\varphi\right)(y)=N\sum_{i=1}^{N}\binom{N-1}{i-1}y^{i-1} (1-y)^{N-i} \int_{0}^{\frac{1}{N}}\varphi\left(s+\frac{i-1}{N}\right)\dd s,\quad \forall\, y\in[0,1].
\]
The operators $K_{N}$ are known as Berstein-Kantorovich operators which were introduce in the 30's by Kantorovich in order to extend the properties of Berstein polynomials to non-continuous functions \cite{kantorovich}. In particular, it is known that, for all $\varphi$ in $\mathbf{L}^{2}$, $K_{N}(\varphi)$ converges to $\varphi$ in $\mathbf{L}^{2}$ \cite[Theorem 2.1.2 and p.\,33]{lorentz}. Now, according to the Cauchy-Schwarz inequality we have that
\[
\left|\int_{0}^{1}\!\!\!\iF(y) K_{N}\left(\iF\right)(y)\ \dd y-\int_{0}^{1}\iF(y)^{2} \ \dd y\right|\leq \left\| \iF\right\|_{2}\, \left\|K_{N}(\iF)-{\iF}\right\|_{2}.
\]
But since $K_{N}(\iF)$ converges to $\iF$ in $\mathbf{L}^{2}$, we finally obtain
\[
\E\left[\sum_{i=1}^{N}\Lambda_{\infty,(i)} \int_{\frac{i-1}{N}}^{\frac{i}{N}}F^{-1}_{\Lambda_\infty}(s)\dd s\right]\longrightarrow\|F^{-1}_{\Lambda_\infty}\|_{2}^{2},
\]
when $N$ goes to infinity. This gives the result. \hfill$\Box$

We also have a stronger convergence result for this estimator. It relies on the fact that the empirical measure $\widehat{\mathcal{P}}$ defined in \eqref{eq:empiricalMeasure} must be close (in Wasserstein distance) to the law of $\sigma^{-1}\Lambda_{\infty}$ as soon as the trees are large enough. More precisely, we have the following result of consistency.

\begin{proposition}
\label{lem:cvEmpricalMeasure}
Let $\mathcal{P}$ be the law of $\sigma^{-1}\Lambda_{\infty}$. Let also $\widehat{\mathcal{P}}_{\vec{u}_N}$ be the empirical distribution defined for any multi-integer $\vec{u}_N$ by
\[
\widehat{\mathcal{P}}_{\vec{u}_N}={\frac1N}\sum_{i=1}^{N}\delta_{\widehat{\lambda}\left[\tau^{i}_{u_{i}} \right]}.
\]
Then, the following statement holds,
\begin{equation*}
\forall\,\epsilon>0,~\exists\,A\in\N,~\forall\, u\in\mathbf{S}_{A},\quad\mathbf{P}\left(\limsup_{N\to\infty}d_{W}\left(\widehat{\mathcal{P}}_{\vec{u}_{N}},\mathcal{P} \right)<\epsilon\right)=1.
\end{equation*}
\end{proposition}
\textit{Proof.}
Let $\Pi_{\delta}$ be the canonical projection of $\R$ on $[-\delta,\delta]$, for a positive real number $\delta$. We have
\begin{align}
& d_{W}\left(\widehat{\mathcal{P}}_{\vec{u}_{N}}\left(\omega\right),\mathcal{P} \right) \nonumber \\
& \leq d_{W}\left(\widehat{\mathcal{P}}_{\vec{u}_{N}}\left(\omega\right),\Pi_{\delta}\widehat{\mathcal{P}}_{\vec{u}_{N}}\left(\omega\right) \right)+d_{W}\left(\Pi_{\delta}\widehat{\mathcal{P}}_{\vec{u}_{N}}\left(\omega\right),\Pi_{\delta}\mathcal{P}\right)+d_{W}\left(\mathcal{P},\Pi_{\delta}\mathcal{P} \right), \label{eq:mainbound}
\end{align}
where $\Pi_{\delta}\mu$ denotes the image measure of $\mu$ by $\Pi_{\delta}$.
To obtain the desired result, we need to control each of the three terms in the right hand side of \eqref{eq:mainbound}.

\textbf{Third term.} First, it is clear, for any probability measure $\mu$, that $\Pi_{\delta}$ is a transport of $\mu$ on $\Pi_{\delta}\mu$ which need not be optimal \cite[2. Generalities on Kantorovich transport distances]{BL14}. Hence,
\[
d_{W}\left(\mu,\Pi_{\delta}\mu \right)\leq \sqrt{ \int_{\R}\left|x-\Pi_{\delta}(x) \right|^{2} \mu(dx)}.
\]
It follows, since $x\mapsto x^{2}$ is integrable w.r.t.\ $\mathcal{P}$, that $\delta$ can be chosen in order to have
\begin{equation}
\label{eq:boundWass}
d_{W}\left(\mathcal{P},\Pi_{\delta}\mathcal{P} \right)\leq\sqrt{\E\left[\left(\sigma^{-1}\Lambda_{\infty}\right)^{2}\mathbf{1}_{|\sigma^{-1}\Lambda_{\infty}|>\delta} \right]}<\frac{\e}{3}.
\end{equation}

\textbf{First term.} On the other hand, following the same lines as in the proof of Proposition \ref{prop:CVLS1} (and using Remark \ref{rem:lambda:moment}), one can show that, for any $\epsilon>0$, there exists $A\in\N$ such that, for any $u\in\mathbf{S}_{A}$,
\begin{equation}
\label{eq:bound2}
\mathbf{P}\left(\limsup_{N\to\infty}\left|\frac{1}{N}\sum_{i=1}^{N}\widehat{\lambda}[\tau_{u_{i}}^{i}]^{2}\mathbf{1}_{|\widehat{\lambda}[\tau_{u_{i}}^{i}]|>\delta}-\E\left[\sigma^{-2}\Lambda_{\infty}^{2}\mathbf{1}_{|\sigma^{-1}\Lambda_{\infty}|>\delta} \right] \right|<\epsilon\right)=1.
\end{equation}
This bound allows us to control the first term in the right hand side of \eqref{eq:mainbound} since
\begin{eqnarray*}
d_{W}\left(\widehat{\mathcal{P}}_{\vec{u}_{N}}\left(\omega\right),\Pi_{\delta}\widehat{\mathcal{P}}_{\vec{u}_{N}}\left(\omega\right) \right)&\leq&  \sqrt{ \int_{\R}\left|x-\Pi_{\delta}(x) \right|^{2} \widehat{\mathcal{P}}_{\vec{u}_{N}}(\omega)(dx)}\\
&\leq&\sqrt{\frac{1}{N}\sum_{i=1}^{N}\widehat{\lambda}[\tau_{n_i}^i](\omega)^{2}\mathbf{1}_{|\widehat{\lambda}[\tau_{n_i}^i](\omega)|>\delta}}.
\end{eqnarray*}
Hence, it remains to control the second term.

\textbf{Second term.} Since $\Pi_{\delta}\widehat{\mathcal{P}}_{\vec{u}_{N}}\left(\omega\right)$ and $\Pi_{\delta}\mathcal{P}$ are compactly supported measures, we have
\[
d_{W}\left(\Pi_{\delta}\widehat{\mathcal{P}}_{\vec{u}_{N}}\left(\omega\right),\Pi_{\delta}\mathcal{P} \right)\leq C \sqrt{d_W^{(1)}\left(\Pi_{\delta}\widehat{\mathcal{P}}_{\vec{u}_{N}}\left(\omega\right),\Pi_{\delta}\mathcal{P} \right)},
\]
where $d_W^{(1)}$ is the first order Wasserstein metric. As a consequence, if one gets the result for $d_W^{(1)}$, it gives the result for $d_{W}$. First of all, we have
\begin{equation}
\label{eq:ineqWass}
d_W^{(1)}\left(\Pi_{\delta}\widehat{\mathcal{P}}_{\vec{u}_{N}}\left(\omega\right),\Pi_{\delta}\mathcal{P}\right)\leq d_W^{(1)}\left(\Pi_{\delta}\widehat{\mathcal{P}}_{\vec{u}_{N}}\left(\omega\right),\overline{\Pi_{\delta}\mathcal{P}_{\vec{u}_{N}}}\right)+d_W^{(1)}\left(\overline{\Pi_{\delta}\mathcal{P}_{\vec{u}_{N}}},\Pi_{\delta}\mathcal{P}\right),
\end{equation}
with
\[
\overline{\Pi_{\delta}\mathcal{P}_{\vec{u}_{N}}}=\frac{1}{N}\sum_{i=1}^{N}\mathbf{P}_{\Pi_{\delta}(\widehat{\lambda}[\tau^{i}_{n_{i}}])},
\]
where $\mathbf{P}_{\Pi_{\delta}(\widehat{\lambda}[\tau^{i}_{n_{i}}])}$ denotes the law of $\Pi_{\delta}(\widehat{\lambda}[\tau^{i}_{n_{i}}])$. Since the space $\mathcal{C}([-\delta,\delta],\mathbf{R})$ of continuous functions on $[-\delta,\delta]$, endowed with the uniform topology is separable, there exists a countable dense subset $(f_{k})_{k\geq 1}$ of $\mathcal{C}([-\delta,\delta],\mathbf{R})$. Once again, using the method developed in Proposition \ref{prop:CVLS1}, it is easy to get that, for any positive integer $k$ and all $\omega$ in a set $\Omega_{k}$ of mass $1$,
\[
\lim_{N\to\infty}\left|\Pi_{\delta}\widehat{\mathcal{P}}_{\vec{u}_{N}}\left(\omega\right)f_{k}-\overline{\Pi_{\delta}\mathcal{P}_{\vec{u}_{N}}}f_{k}\right| = 0,
\]
where $\mu f$ denotes $\int f(x) \mu(\dd x)$, for any a measure $\mu$ and any $\mu$-integrable function $f$. Now, take $\omega$ in $\bigcap_{k\geq 1}\Omega_{k}$ and $f$ in $\mathcal{C}([-\delta,\delta],\mathbf{R})$. Since $(f_{k})_{k\geq 1}$ is dense in $\mathcal{C}([-\delta,\delta],\mathbf{R})$, there exists for any $\epsilon>0$ an integer $k$ such that $\|f-f_{k}\|<\epsilon/2$. This implies that
\begin{eqnarray*}
\limsup_{N\to\infty}\left|\Pi_{\delta}\widehat{\mathcal{P}}_{\vec{u}_{N}}\left(\omega\right)f-\overline{\Pi_{\delta}\mathcal{P}_{\vec{u}_{N}}}f\right|&\leq& \limsup_{N\to\infty}\left|\Pi_{\delta}\widehat{\mathcal{P}}_{\vec{u}_{N}}\left(\omega\right)f-\Pi_{\delta}\widehat{\mathcal{P}}_{\vec{u}_{N}}\left(\omega\right)f_{k}\right|\\
&+&\limsup_{N\to\infty}\left|\Pi_{\delta}\widehat{\mathcal{P}}_{\vec{u}_{N}}\left(\omega\right)f_{k}-\overline{\Pi_{\delta}\mathcal{P}_{\vec{u}_{N}}}f_{k}\right|\\
&+&\limsup_{N\to\infty}\left|\overline{\Pi_{\delta}\mathcal{P}_{\vec{u}_{N}}}\left(\omega\right)f_{k}-\overline{\Pi_{\delta}\mathcal{P}_{\vec{u}_{N}}}f\right|\qquad<\epsilon.
\end{eqnarray*}
Since this last inequality holds for any $\epsilon>0$, $\Pi_{\delta}\widehat{\mathcal{P}}_{\vec{u}_{N}}-\overline{\Pi_{\delta}\mathcal{P}_{\vec{u}_{N}}}$ converges weakly to $0$ in the space $\mathcal{M}_{s}([-\delta,\delta])$ of signed measures on $[-\delta,\delta]$ with probability $1$. Consequently, since $d_W^{(1)}$ metricizes the subspace of probability measures, we get that, almost surely,
\begin{equation}
\label{eq:ascv}
\lim_{N\to\infty}d_W^{(1)}\left(\Pi_{\delta}\widehat{\mathcal{P}}_{\vec{u}_{N}}\left(\omega\right),\overline{\Pi_{\delta}\mathcal{P}_{\vec{u}_{N}}}\right)=0 .
\end{equation}
In order to control the second term and, hence, end the proof, it remains to show that
\begin{equation*}
\forall\,\epsilon>0,~\exists\,A\in\N,~\forall\, u\in\mathbf{S}_{A},~d_W^{(1)}\left(\overline{\Pi_{\delta}\mathcal{P}_{\vec{u}_{N}}},\Pi_{\delta}\mathcal{P}\right)<\epsilon.
\end{equation*}
To get this, we use the duality formula for the first order Wasserstein distance,
\begin{equation}
\label{eq:duality}
d_W^{(1)}\left(\overline{\Pi_{\delta}\mathcal{P}_{\vec{u}_{N}}},\Pi_{\delta}\mathcal{P}\right)=\sup_{ \phi\in Lip_{1}\left([-\delta,\delta],\mathbf{R}\right)}\left|\overline{\Pi_{\delta}\mathcal{P}_{\vec{u}_{N}}}\phi-\Pi_{\delta}\mathcal{P}\phi\right|,
\end{equation}
where $Lip_{1}\left([-\delta,\delta],\mathbf{R}\right)$ denotes the set of $1$-Lipschitz continuous functions on $[-\delta,\delta]$. Now, let $\phi$ be an element of $Lip_{1}\left([-\delta,\delta],\mathbf{R}\right)$, we have that
\[
\left|\overline{\Pi_{\delta}\mathcal{P}_{\vec{u}_{N}}}\phi-\Pi_{\delta}\mathcal{P}\phi\right|\leq \frac{1}{N}\sum_{i=1}^{N} \left|\mathbf{E}\left[\phi\left(\Pi_{\delta}(\widehat{\lambda}[\tau_{n_{i}}^{i}])\right)\right]-\mathbf{E}\left[\phi\left(\Pi_{\delta}(\Lambda_{\infty})\right)\right]\right|.
\]
To prove that the supremum taken in the above inequality will be small as soon as the trees are large enough, we use the following lemma that gives the uniform convergence of the expectation of Lipschitz functionals of the $\widehat{\lambda}[\tau^i_{n_i}]$'s.
\begin{lemma}For any $\epsilon>0$, there exists $A\in\mathbf{N}$ such that,
	\label{lem:unifLip}
	\[
	\forall\,n>A,~\forall\,f\in Lip_{1}\left([-\delta,\delta],\mathbf{R}\right),~\left|\mathbf{E}\left[f\left(\Pi_{\delta}(\widehat{\lambda}[\tau_{n_i}^i])\right)\right]-\mathbf{E}\left[f\left(\Pi_{\delta}(\Lambda_{\infty}) \right) \right]\right|<\epsilon.
	\]
\end{lemma}
The proof of the lemma has been postponed to the end of the section. In light of this lemma, we get
\[
\forall\,\epsilon>0,~\exists\,A\in\mathbf{N},~\forall\,n>A,~d_W^{(1)}\left(\overline{\Pi_{\delta}\mathcal{P}_{\vec{u}_{N}}},\Pi_{\delta}\mathcal{P}\right)<\epsilon.
\]
Consequently, using the above result in conjunction with \eqref{eq:ineqWass}, \eqref{eq:ascv} and \eqref{eq:duality}, we have
\[
\forall\,\epsilon>0,~\exists\,A\in\mathbf{N},~\forall\,n>A,~\limsup_{N\to\infty}d_W^{(1)}\left(\Pi_{\delta}\widehat{\mathcal{P}}_{\vec{u}_{N}}(\omega),\Pi_{\delta}\mathcal{P}\right)<\epsilon,
\]
for all $\omega$ in $\bigcup_{k\geq 1}\Omega_{k}$, where the sets $\Omega_{k}$ have been defined above.
Hence, for any $\epsilon>0$, there exists $A\in\N$ such that,
\begin{equation}
\label{eq:bound3}
\forall\, u\in\mathbf{S}_{A},\quad\mathbf{P}\left(\limsup_{N\to\infty}d_{W}\left(\Pi_{\delta}\widehat{\mathcal{P}}_{\vec{u}_{N}}\left(\omega\right),\Pi_{\delta}\mathcal{P} \right)<\epsilon\right)=1.
\end{equation}
To end, using \eqref{eq:boundWass}, \eqref{eq:bound2} and \eqref{eq:bound3} in \eqref{eq:mainbound} leads to the result.\hfill$\Box$

Finally, we get the following consistency result.
\begin{proposition}
\label{prop:limit:w:as}
We have,
\begin{equation*}
\forall\,\epsilon>0,~\exists\,A\in\N,~\forall\, u\in\mathbf{S}_{A},\quad\mathbf{P}\left(\limsup_{N\to\infty}\left|\widehat{\lambda}_W[\mathcal{F}_{\vec{u}_{N}}] -\frac{1}{\sigma}\right|<\epsilon\right)=1.
\end{equation*}
\end{proposition}
\textit{Proof.}
By the Cauchy-Schwarz inequality, the convergence of this estimator is conditioned to the convergence of the Wasserstein distance in the following manner,
\begin{eqnarray*}
\left|\widehat{\lambda}_W[\mathcal{F}_{\vec{u}_N}] -\frac{1}{\sigma}\right| &=&\frac{\left|\langle\widehat{F}^{-1}_{\vec{u}_N}-\sigma^{-1}F^{-1}_{\Lambda_\infty}, F^{-1}_{\Lambda_\infty}\rangle\right|}{\|F^{-1}_{\Lambda_\infty}\|^{2}_2}\\
&\leq&\frac{\left\|\widehat{F}^{-1}_{\vec{u}_N}-\sigma^{-1}F^{-1}_{\Lambda_\infty}\right\| _2\left\|F^{-1}_{\Lambda_\infty} \right\|_2}{\|F^{-1}_{\Lambda_\infty}\|^{2}_2}\\
&=&\frac{d_W\left(\widehat{\mathcal{P}}_{\vec{u}_N} \, , \, \mathcal{P}\right)}{\|F^{-1}_{\Lambda_\infty}\|_2}.
\end{eqnarray*}
The result finally arises from Proposition \ref{lem:cvEmpricalMeasure} concerning the convergence of the empirical measure in the sense of the Wasserstein distance.\hfill$\Box$

\textit{Proof of Lemma \ref{lem:unifLip}.} Let $(n_{i})_{i\geq 1}$ be a strictly increasing sequence of integers and, for any $i\geq1$, let $\tau^i_{n_{i}}$ be a conditioned Galton-Watson tree with size $n_{i}$. By virtue of Skorokhod's representation theorem, there exists a probability space on which is defined a sequence $(X_{n_{i}})_{i\geq1}$ of random variables and a random variable $X_{\infty}$ such that:
	\begin{itemize}[label={$\diamond$}]
		\item $(X_{n_{i}})_{i\geq1}$ has the same distribution as the sequence $(\Pi_{\delta}(\widehat{\lambda}[\tau^i_{n_{i}}]))_{i\geq 1}$;
		\item $X_{\infty}$ has the same distribution as $\Pi_{\delta}(\Lambda_{\infty})$;
		\item $X_{n_{i}}$ converges to $X_{\infty}$ in probability as $i$ goes to infinity.
	\end{itemize}
As a consequence of the uniform integrability properties showed in the proof of Proposition \ref{prop:w:law}, we get
\begin{equation}\label{eq:lim:l1:e}\lim_{i\to\infty}\mathbf{E}\left[\left|X_\infty - X_{n_i}\right|\right] = 0.\end{equation}
Now, choose $f$ in $Lip_{1}\left([-\delta,\delta],\mathbf{R}\right)$, we have
\[
\left|\mathbf{E}\left[f\left(X_{n_{i}}\right)\right]-\mathbf{E}\left[f\left(X_{\infty}\right)\right]\right|\leq\mathbf{E}\left[\left|X_{n_{i}}-X_{\infty}\right| \right].
\]
Together with \eqref{eq:lim:l1:e}, we obtain the expected convergence.\hfill$\Box$


\section{Simulation study}\label{sec:simu}

\subsection{Simulation of conditioned Galton-Watson trees}
\label{ss:simugw}

In order to illustrate our estimation techniques on Galton-Watson forests, we need to make some numerical experiments. However, simulation of conditioned Galton-Watson trees is a difficult problem of independent importance. In this section, we briefly present an algorithm due to \cite{Dev12} allowing to achieve this aim.
Given an integer $n$ and a distribution $\mu$ on the set $\{0,\dots,K \}$, this algorithm provides, in two steps, the simulation of the \L{}ukasciewicz walk $\mathcal{L}[\tau_{n}]$ of a tree $\tau_{n}$ with distribution $\text{GW}_{n}(\mu)$. Three more steps are required to obtain the corresponding Harris path $\mathcal{H}[\tau_{n}]$ through other coding processes (see for example \cite{Duq03}).

\begin{itemize}[label={$\diamond$}]

\item \textbf{Simulation of numbers of children.}  The multinomial distribution of parameters $(\mu(k))_{0\leq k\leq K}$ and $n$ may be defined by its probability mass function,
$$\prob(N_0 = n_0,\dots,N_K=n_K) =
\left\{
\begin{array}{cl}
\displaystyle\frac{n!}{n_0!\,\dots\,n_K!} \mu(0)^{n_0}\dots\mu(K)^{n_K} &\text{if}\quad\displaystyle\sum_{k=0}^K{n_k}=n,\\
\displaystyle0&\text{else.}
\end{array}
\right.
$$
Simulation of the multinomial distribution presents no difficulty. By rejection sampling, we simulate multinomial random variables until obtaining a sequence $(N_k)_{0\leq k\leq K}$ satisfying
$$\sum_{k=0}^K k N_k = n-1 .$$
We define the sequence $(\zeta_i)_{1\leq i\leq n}$ from
$$(\zeta_i)_{1\leq i\leq n} = (\underbrace{0,\dots,0}_{N_0}\,,\,\underbrace{1,\dots,1}_{N_1}\,,~\dots~,\,\underbrace{K,\dots,K}_{N_K} ).$$
Let $(\xi_i)_{1\leq i\leq n}$ be a sequence obtained as a random permutation of $(\zeta_i)_{1\leq i\leq n}$. A suitable technique for random shuffling is presented in \cite[Algorithm P (p.\,139)]{knu81}. The sequence $(\xi_i)_{1\leq i\leq n}$ represents the verticesnumbers of children in the depth-first search order.

\item \textbf{Computation of the \L{}ukasciewicz walk.} Let $L$ be the process defined by $L(0)=0$ and,
$$\forall\,0\leq k\leq n-2,\quad L(k+1) = L(k)+\xi_{k+1} - 1 .$$
Set $l=1+\argmin\,\{ L(k)~:~0\leq k\leq n-1\}$. Then there exists a tree $\tau_n$ with $n$ nodes whose \L{}ukasciecwicz walk is defined by
$$\mathcal{L}[\tau_n](k) =
\left\{
\begin{array}{lll}
L(l+k) + \min L - 1 & \text{if} & 0\leq k\leq n-1-l ,\\
L(k-n+l) + \min L - 1 &\text{if} & n-l\leq k\leq n-1.
\end{array}
\right.
$$

\item \textbf{From the \L{}ukasiewicz walk to the height process.} Now, we compute the corresponding height process \cite[eq.(2)]{Duq03},
$$\forall\,0\leq k\leq n-1,\quad\frak{H}[\tau_n](k) = \#\left\{0\leq j\leq k-1~:~\mathcal{L}[\tau_n](j) = \min_{j\leq l\leq n} \mathcal{L}[\tau_n](l)\right\}.$$

\item \textbf{From the height process to the contour process.} Let $(b_k)_{0\leq k\leq n-1}$ be the sequence defined from $b_k = 2k-\frak{H}[\tau_n](k)$ if $0\leq k\leq n-1$ and $b_n=2(n-1)$. Then the $b_i$'s are sorted in increasing order. The contour process $\mathcal{C}[\tau_n](k)$ is defined for any $0\leq k\leq 2n-2$ in \cite[eq.(1)]{Duq03}
$$
\mathcal{C}[\tau_n](k) =
\left\{
\begin{array}{llll}
\frak{H}[\tau_n](i) - (k-b_i) &\text{if}&\exists\,0\leq i\leq n-2,&b_i\leq k<b_{i+1}-1,\\
k-b_{i+1}+\frak{H}[\tau_n](i+1)&\text{if}&\exists\,0\leq i\leq n-2,&b_{i+1}-1\leq k<b_{i+1},\\
\frak{H}[\tau_n](b_{n-1}) - (k-b_{n-1}) &\text{if}&&b_{n-1}\leq k\leq b_n.
\end{array}
\right.
$$
\item \textbf{From the contour process to the Harris path.} The Harris path is only a small modification of the contour process, defined by $\mathcal{H}[\tau_n](0)=\mathcal{H}[\tau_n](2n)=0$ and
$$\forall\,1\leq k\leq 2n-1,~\mathcal{H}[\tau_n](k) = \mathcal{C}[\tau_n](k-1)+1.$$

\end{itemize}


\subsection{Inference for a forest of binary conditioned Galton-Watson trees}
\label{ss:simulateddata}

The aim of this section is to analyze the finite-sample behavior of both estimators introduced in this paper by means of numerical experiments. The theoretical study achieved in Section \ref{sec:mainResults} shows that we can expect to obtain good numerical results, at least for large trees. To this goal, we consider a forest of independent conditioned Galton-Watson trees with common critical birth distribution $\mu$ such that
 $\mu(k)=0$ for $k\geq3$. Such a distribution satisfies the following linear system of equations,
$$
\left\{\begin{array}{ccc}
\mu(0)+\mu(1)+\mu(2)&=&1\\
\mu(1)+2\mu(2)&=&1\\
\mu(1)+4\mu(2)-1&=&\sigma^2
\end{array}\right.
$$
which is equivalent to
$$
\mu(0)=\mu(2) = \frac{\sigma^2}{2} \qquad\text{and}\qquad \mu(1)=1-\sigma^2.
$$
In other words, $\mu$ is entirely characterized by its variance $\sigma^2$. Simulations of Galton-Watson trees $\text{GW}_n(\mu)$ are performed with the method provided in Subsection \ref{ss:simugw}.

Let $\mathcal{F}=(\tau^i)_{1\leq i\leq N}$ be a forest of $N$ independent trees such that, for any $1\leq i\leq N$, $\tau^i\sim\text{GW}_{n_i}(\mu)$ for some integer $n_i$. From the Harris process of each tree $\tau^i$, one first computes the quantity
\begin{eqnarray*}
\widehat{\lambda}\left[\tau^{i} \right] &=&\frac{\langle\mathcal{H}[\tau^{i}](2n_{i}\cdot), E\rangle}{2\sqrt{n_{i}}\|E\|_2^2},
\end{eqnarray*}
where $E$ is known and defined in \eqref{eq:Et}.
Then, we propose to estimate $\sigma^{-1}$ in the two following ways, where $(\widehat{\lambda}[\tau^{(i)}])_{1\leq i\leq N}$ denotes the order statistic associated to the family $(\widehat{\lambda}[\tau^{i}])_{1\leq i\leq N}$. 
\begin{center}
\begin{tabular}{l|l}
Least Squares & $\displaystyle\widehat{\lambda}_{ls}[\mathcal{F}]=\frac{1}{N}\sum_{i=1}^N\widehat{\lambda}\left[\tau^{i}\right]$ \\ \hline
Wasserstein & $\displaystyle\widehat{\lambda}_{W}[\mathcal{F}]=\frac{1}{\|F^{-1}_{\Lambda_\infty}\|_2^{2}}\sum_{i=1}^{N}\widehat{\lambda}\left[\tau^{(i)}\right]\int_{\frac{i-1}{N}}^{\frac{i}{N}}F^{-1}_{\Lambda_\infty}(s)\dd s$
\end{tabular}
\end{center}

\begin{remark}\label{rem:Flambdamoins1}
In order to compute $\widehat{\lambda}_{W}[\mathcal{F}]$, we need to be able to perform computations using the function $F_{\Lambda_{\infty}}^{-1}$. Unfortunately, in view of the theoretical study of $\Lambda_{\infty}$ made in Proposition \ref{ssec:adequacy}, one can not expect to have an explicit expression for this function. In the following of this section, we use a numerical estimation of $F^{-1}_{\Lambda_{\infty}}$ by Monte Carlo simulations. To achieve this goal, we perform simulations of $\Lambda_{\infty}$ thanks to formula \eqref{eq:lambdainfty} by simulating Brownian excursions thanks to \eqref{eq:rep:3d}. In order to ensure that the error made on $F^{-1}_{\Lambda_{\infty}}$ does not propagate too much in our results, $F^{-1}_{\Lambda_\infty}$ is estimated from one million simulations of $\Lambda_{\infty}$.
\end{remark}

\begin{figure}[!h]
\centerfloat
	\includegraphics[scale=0.25]{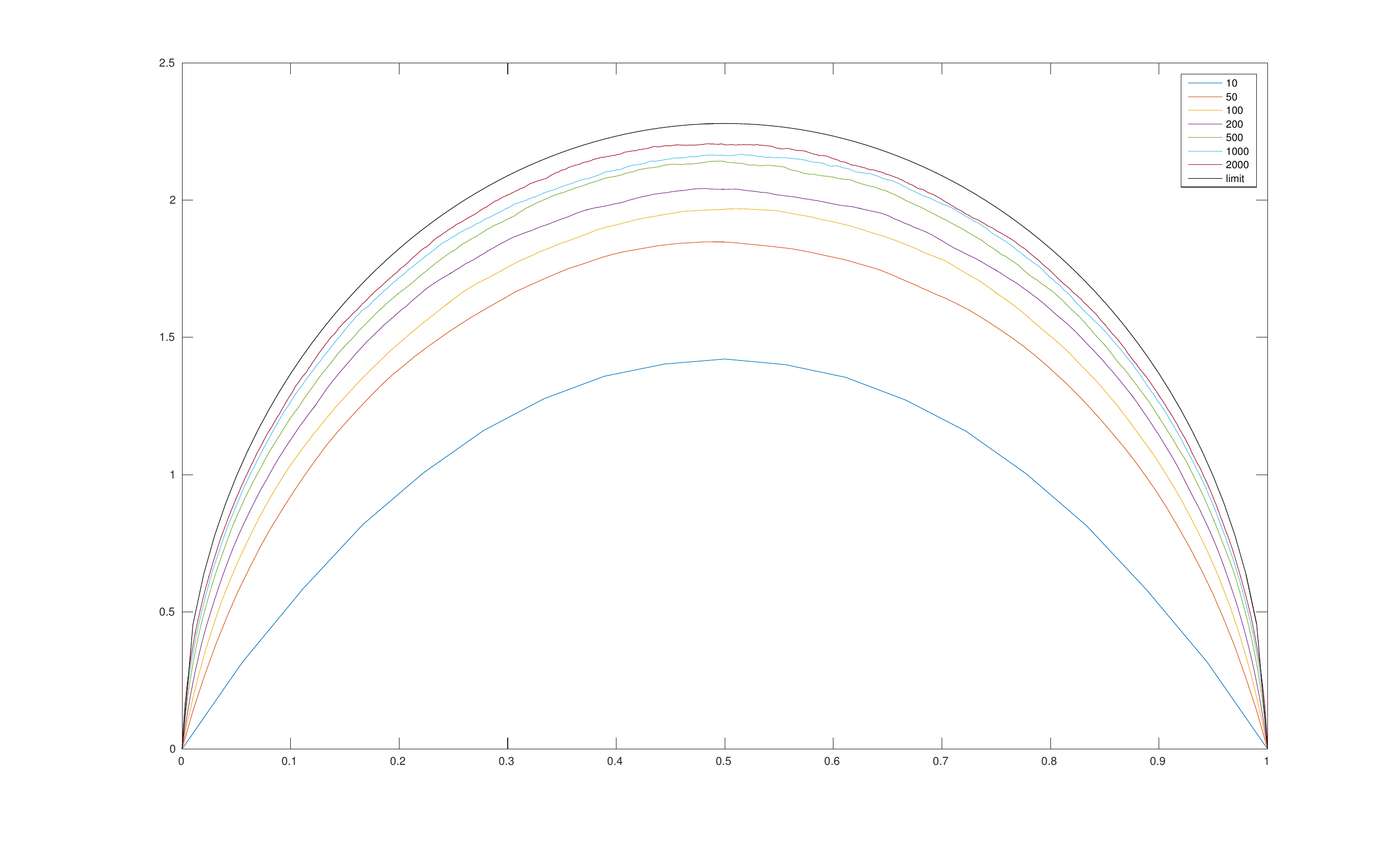}
	\caption{Estimated mean Harris processes of binary conditioned Galton-Watson trees with size $n$ and $\sigma=0.7$ calculated from $2000$ trees for each value of $n$.}
	\label{fig:biasMeanContour}
\end{figure}

The theoretical investigations of Section \ref{sec:mainResults} establish that our estimators are asymptotically unbiased. Nevertheless, the problem is not as simple when working with finite trees. A clear illustration of this comes from the numerical evaluations of the average Harris processes of finite trees. Indeed, the numerical study of Figure \ref{fig:biasMeanContour} shows that the average Harris processes of small trees seem to be lower than the limiting Harris process. Hence, the quantities $\widehat{\lambda}[\tau^i]$ are expected to underestimate the target $\sigma^{-1}$. But any estimator based on the asymptotic behavior of conditioned Galton-Watson trees is expected to present such a bias. In particular, we state in our numerical experiments that the estimator proposed in \cite{KPDRV14} presents the same bias.

\begin{figure}[!h]
	\centerfloat
	\includegraphics[scale=0.25]{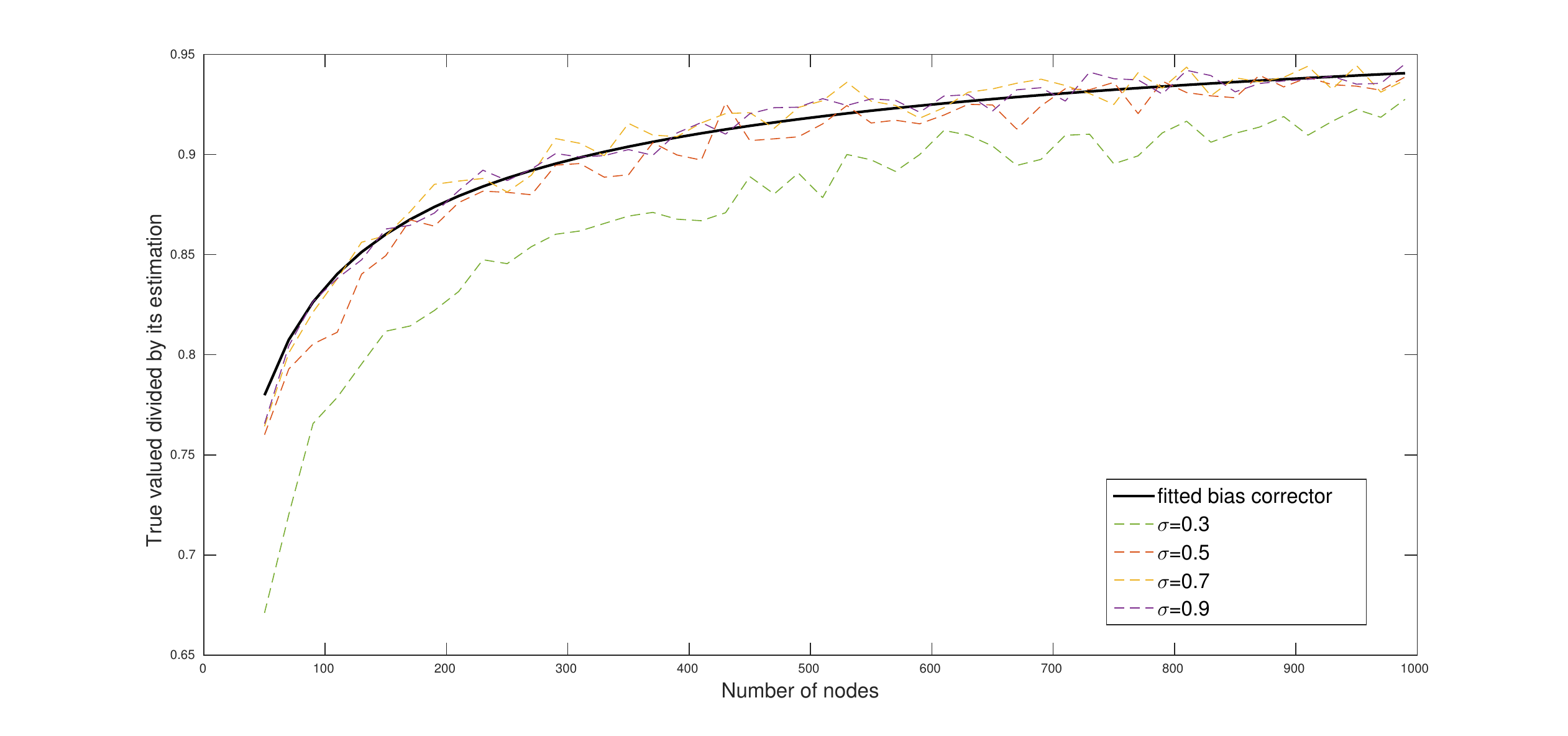}
	\caption{Estimation of the quantity $\eta(n)=\sigma^{-1}\E[\widehat{\lambda}[\tau_n]]^{-1}$, where $\tau_n\sim\text{GW}_n(\mu)$, for different values of $\sigma$ and different numbers of nodes $n$, together with the fitted bias corrector function $\widehat{\eta}$. Estimations have been made by Monte Carlo method with samples of $2000$ trees for each couple $(n,\sigma)$.}
	\label{fig:biasvssigma}
\end{figure}

The natural question arising from the preceding comments is: how is the bias of a conditioned Galton-Watson tree related to its size and/or the unknown parameter $\sigma$ ?
The numerical study presented in Figure \ref{fig:biasvssigma} shows that the quantity $\eta(n) = \sigma^{-1} \E[\widehat{\lambda}[\tau_n]]^{-1}$, where $\tau_n\sim\text{GW}_n(\mu)$, seems close to uncorrelated to $\sigma$ at least when $\sigma$ is large enough. This allows us to construct a bias corrector which is independent of the unknown standard deviation $\sigma$. In addition, the dependency on $n$ may be modeled by the relation $\eta(n)=1-(a\sqrt{n}+b)^{-1}$ . The coefficients appearing in $\eta$ may be estimated from simulated data,
$$\widehat{\eta}(n) = 1-(0.504273\sqrt{n}+0.9754839)^{-1}$$
(see Figure \ref{fig:biasvssigma} again). The correction is obviously expected to be better for large values of $\sigma$. Finally, we construct the following corrected versions of the estimators $\widehat{\lambda}_{ls}[\mathcal{F}]$ and $\widehat{\lambda}_{W}[\mathcal{F}]$.
\begin{center}
\begin{tabular}{l|l}
Corrected LS & $\displaystyle\widehat{\lambda}_{ls}^c[\mathcal{F}]=\frac{1}{N}\sum_{i=1}^N\widehat{\eta}(\#\tau^i)\widehat{\lambda}\left[\tau^{i} \right]$\\ \hline
Corrected Wasserstein & $\displaystyle\widehat{\lambda}_{W}^c[\mathcal{F}]=\frac{1}{\|F^{-1}_{\Lambda_\infty}\|_2^{2}}\sum_{i=1}^{N}\widehat{\eta}\left(\#\tau^{(i)}\right)\widehat{\lambda}\left[\tau^{(i)}\right]\int_{\frac{i-1}{N}}^{\frac{i}{N}}F^{-1}_{\Lambda_\infty}(s)\dd s$
\end{tabular}
\end{center}

\noindent
In light of the previous comments, computing the estimators proposed in this paper is not an easy task. According to Remark \ref{rem:Flambdamoins1}, one needs to perform a significant number of simulations of $\Lambda_\infty$ in order to get an accurate approximation of $F^{-1}_{\Lambda_\infty}$. Moreover, to be able to correct the aforementioned bias, one needs to perform many simulations of finite trees. Together with this work, we propose a \verb+Matlab+ toolbox which already includes these preliminary computations and allows to directly and quickly compute our estimators for forests. This toolbox as well as its documentation and the scripts used in this paper are available from the authors upon request.

\begin{figure}[!h]
	\centerfloat
	\includegraphics[width=5.5cm,height=3cm]{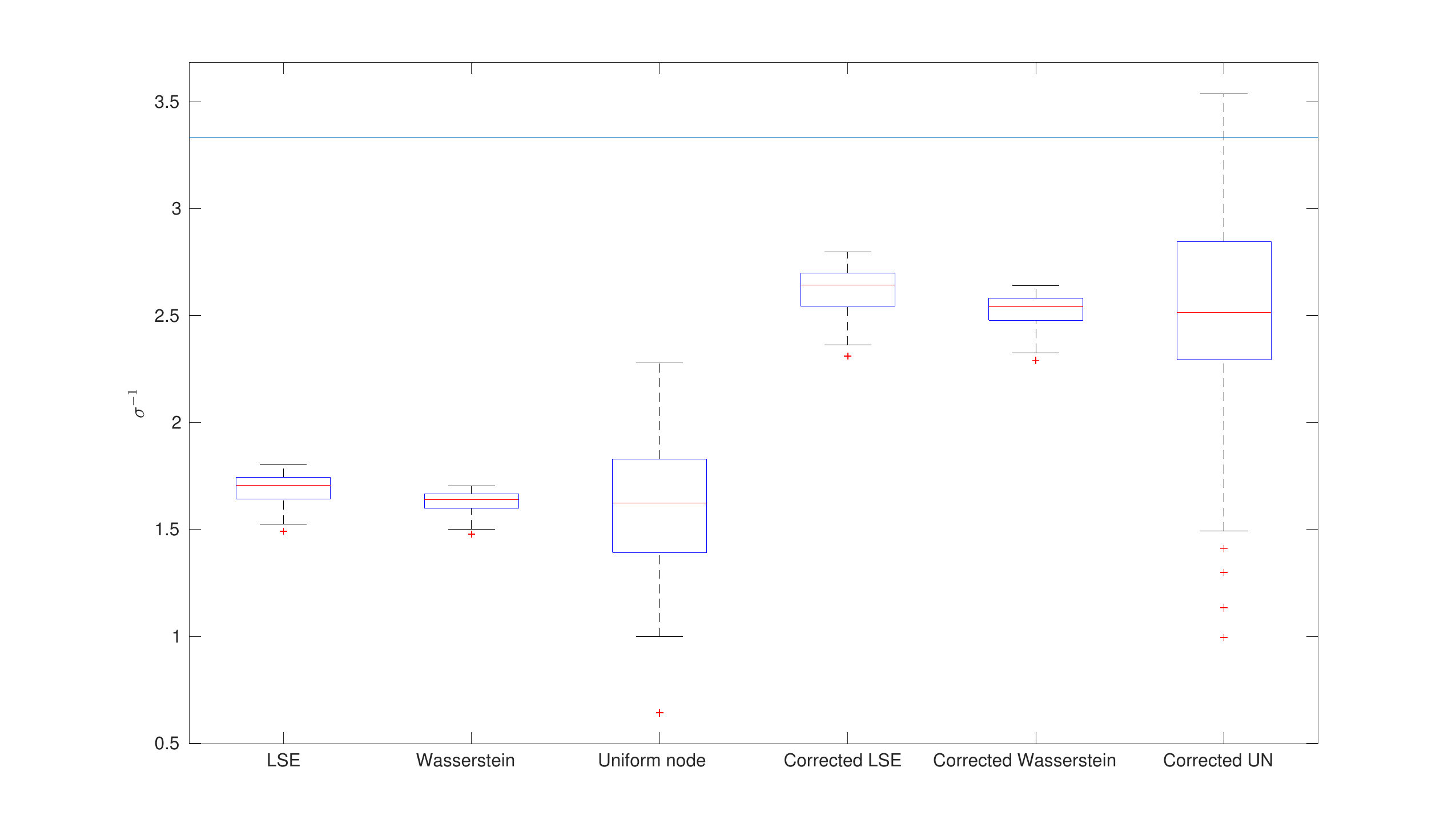}\includegraphics[width=5.5cm,height=3cm]{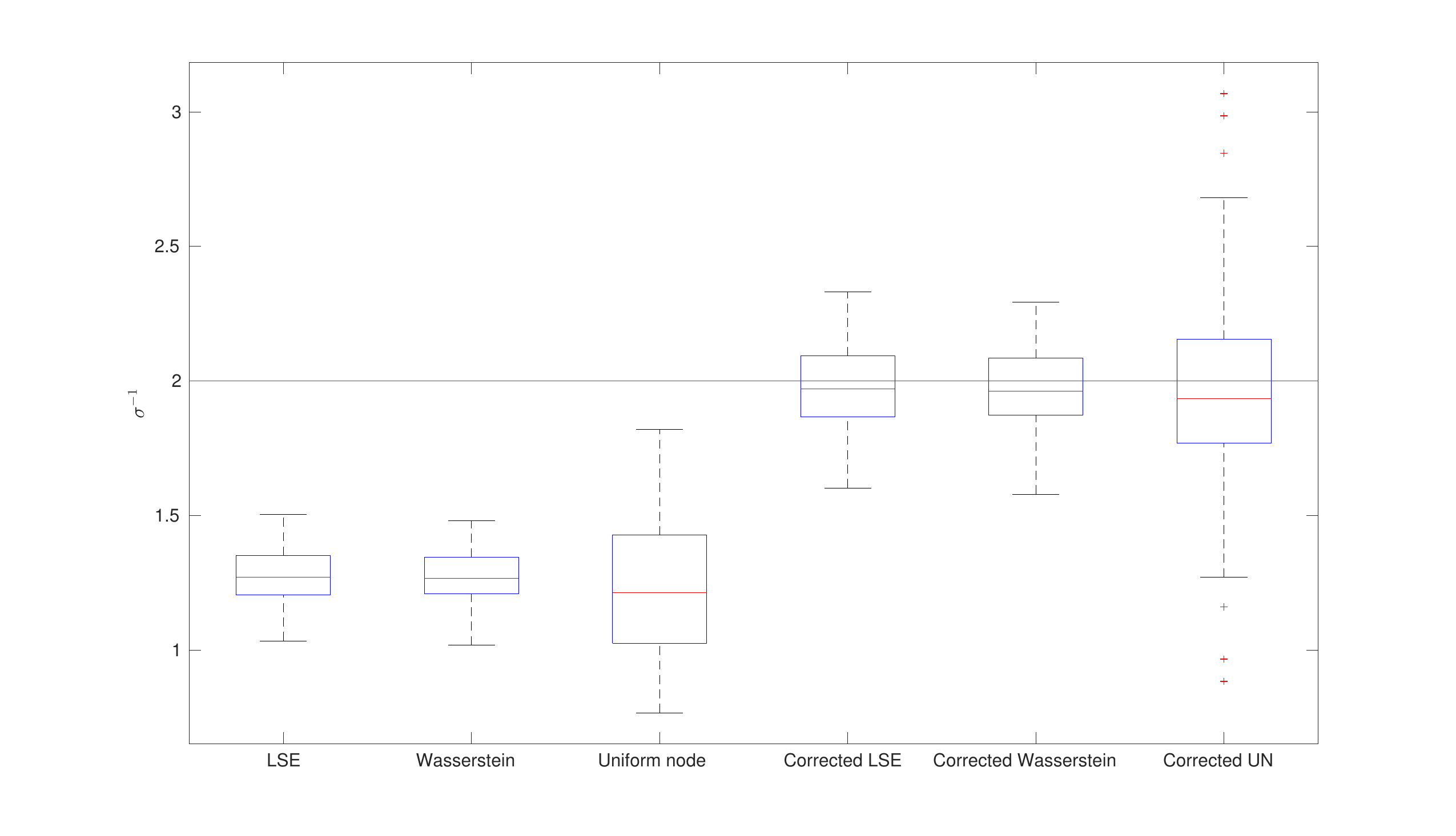}\\
	\includegraphics[width=5.5cm,height=3cm]{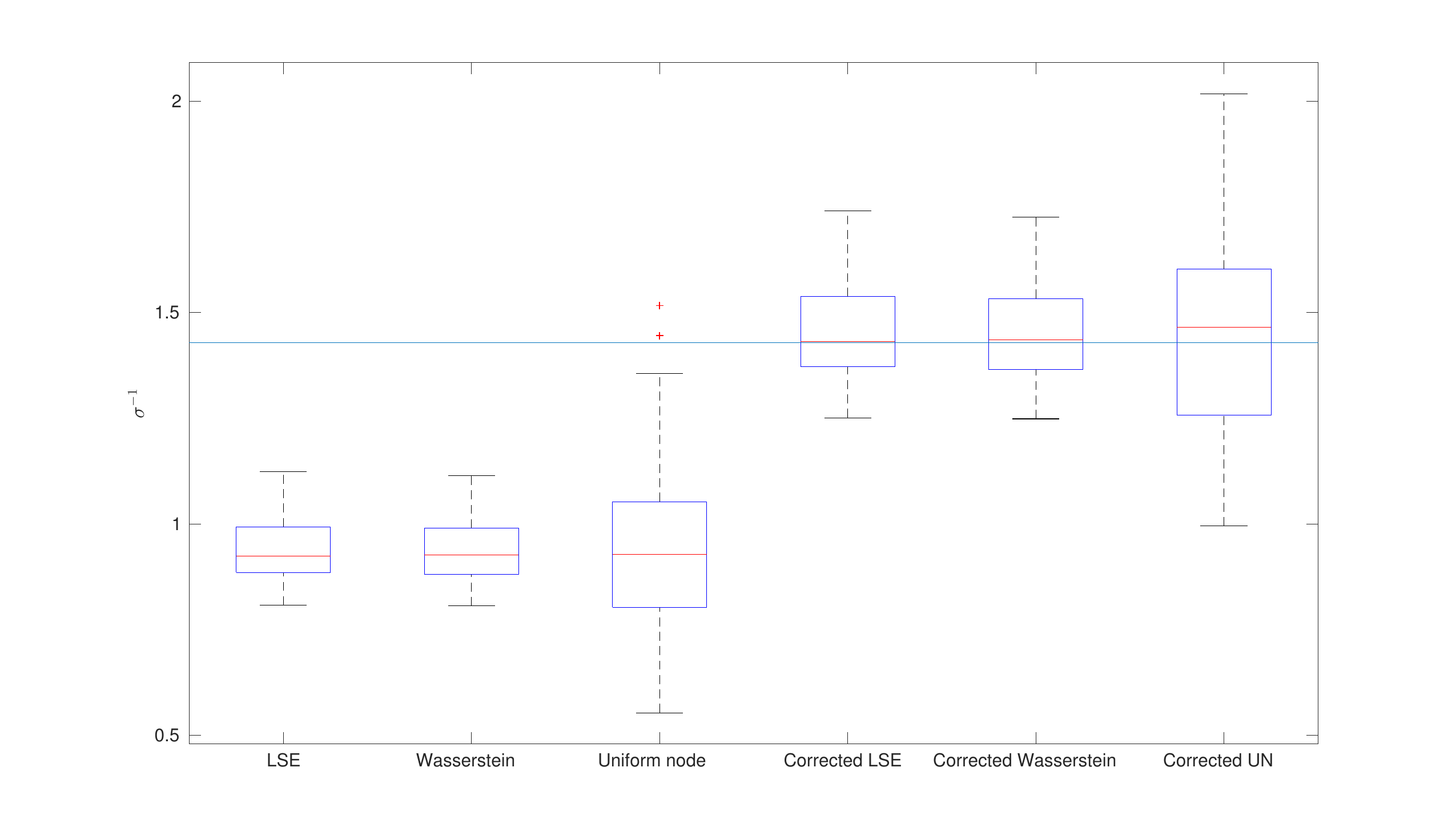}\includegraphics[width=5.5cm,height=3cm]{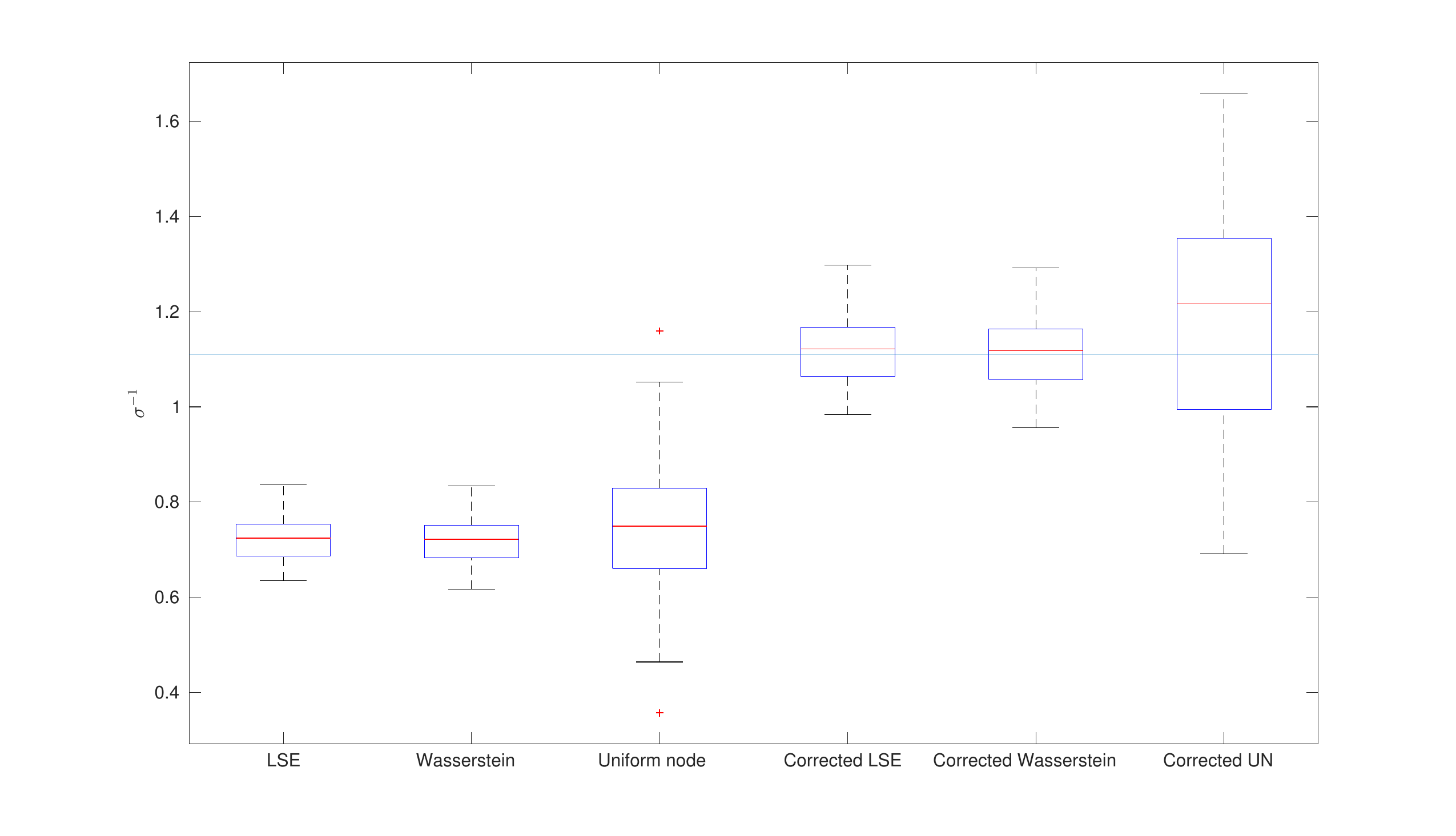}
	\caption{Estimation and bias correction for forests of $10$ trees with $20$ nodes for $\sigma$ equals to $0.3$ (top, left) $0.5$ (top, right), $0.7$ (bottom, left) and $0.9$ (bottom right). Boxplots have been drawn from $100$ replicates each.}
	\label{fig:biascorrection}
\end{figure}

\begin{figure}[!h]
	\centering
	\includegraphics[width=3.5cm,height=3cm]{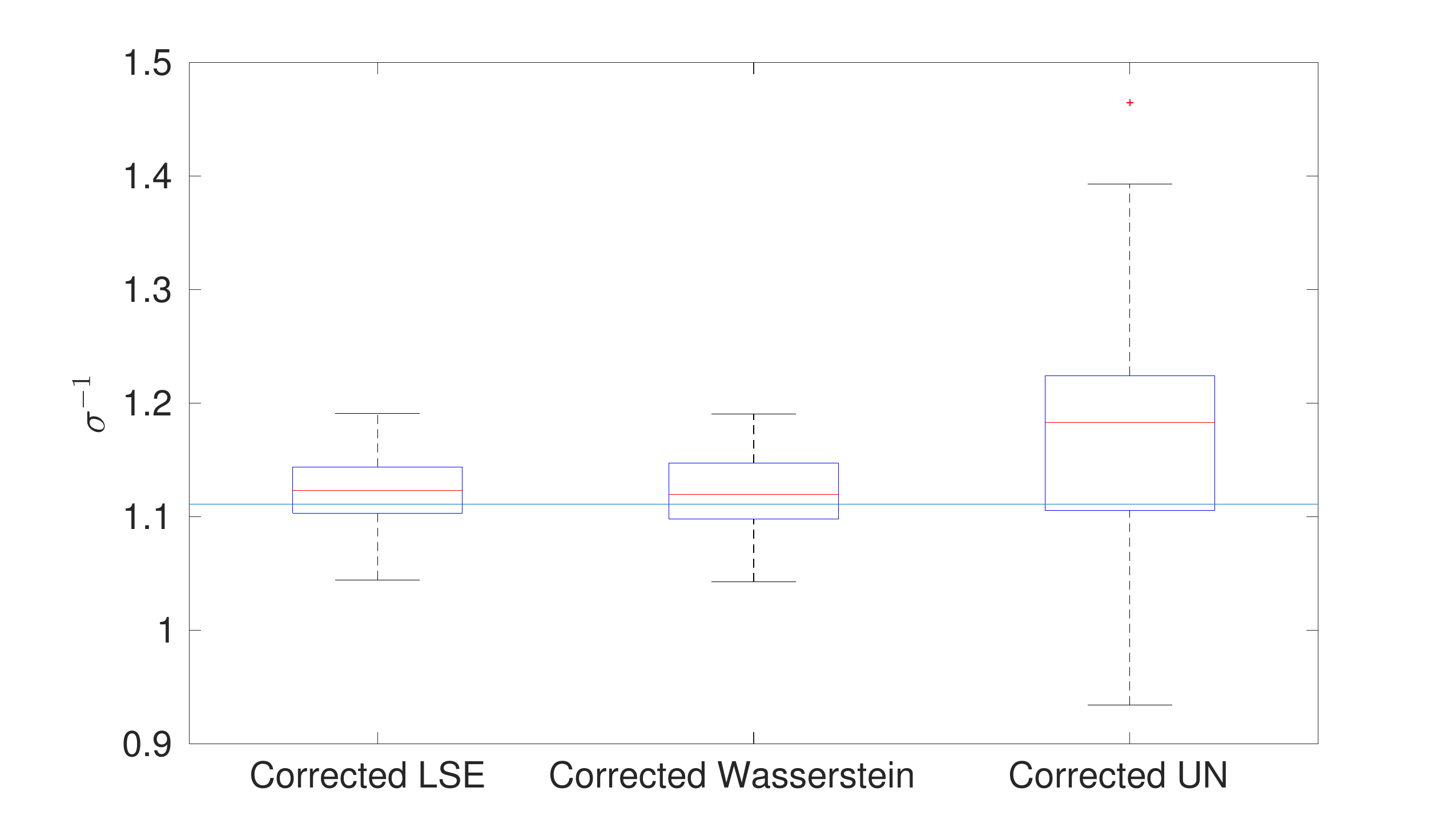}\includegraphics[width=3.5cm,height=3cm]{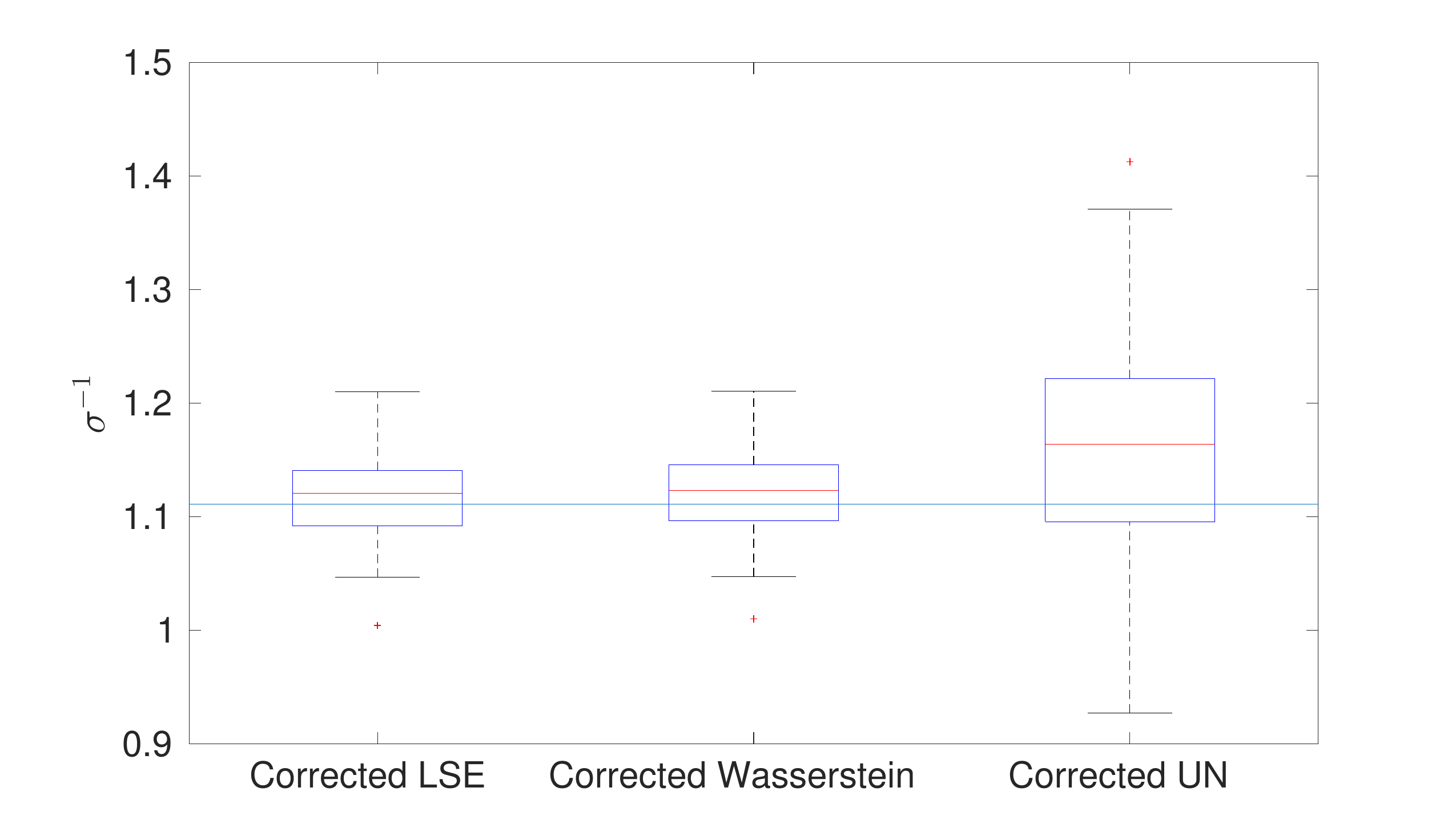}\includegraphics[width=3.5cm,height=3cm]{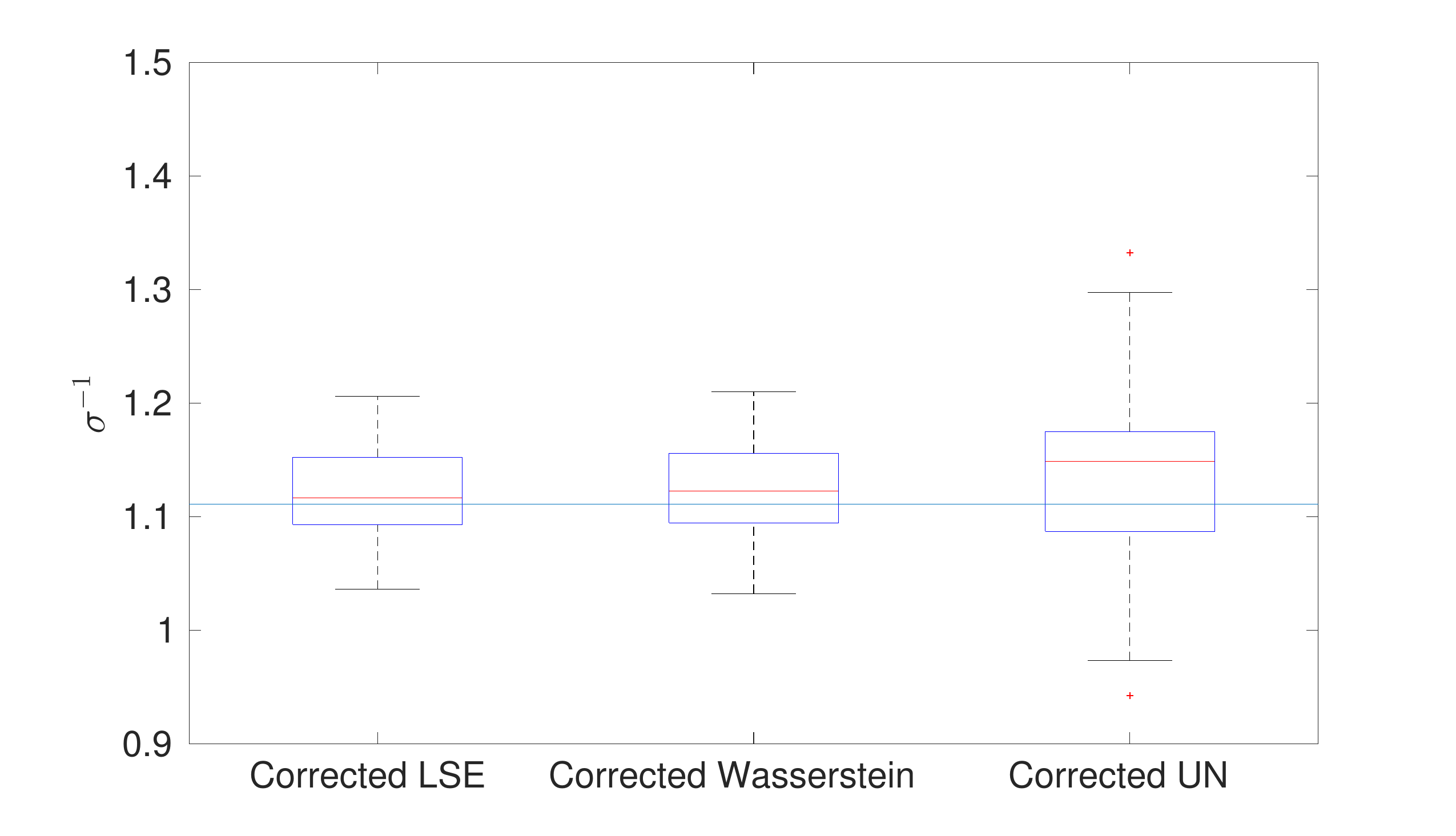}
	\caption{Influence of the size of the trees for $\sigma$ equals to $0.9$: tree sizes varying from $20$ nodes (left), $50$ nodes (center), to $100$ nodes (right). Forests of $50$ trees. Boxplots have been drawn from $100$ replicates each.}
	\label{fig:sizecomp}
\end{figure}

\begin{figure}[!h]
	\centerfloat
	\includegraphics[width=3.5cm,height=3cm]{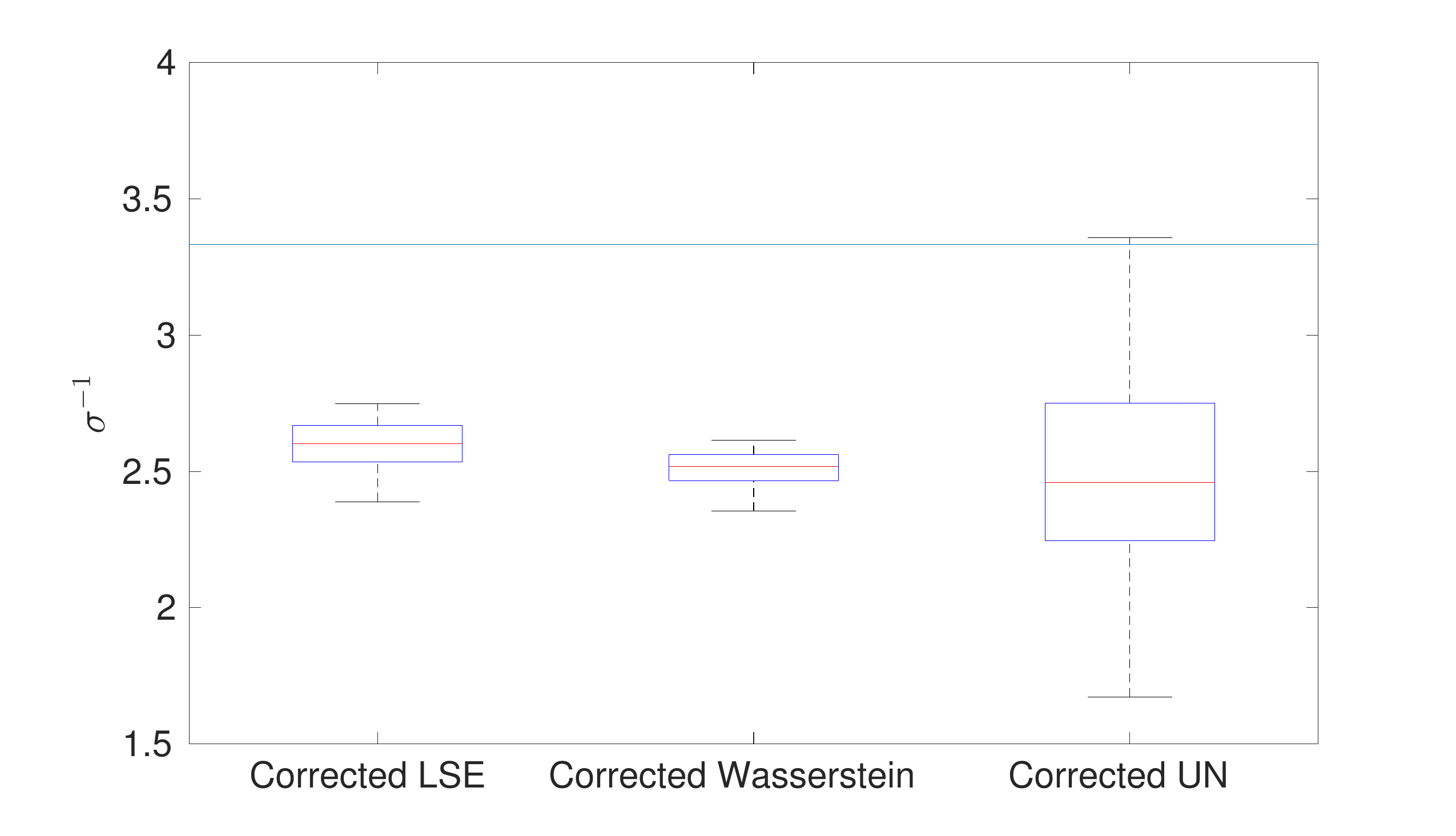}\includegraphics[width=3.5cm,height=3cm]{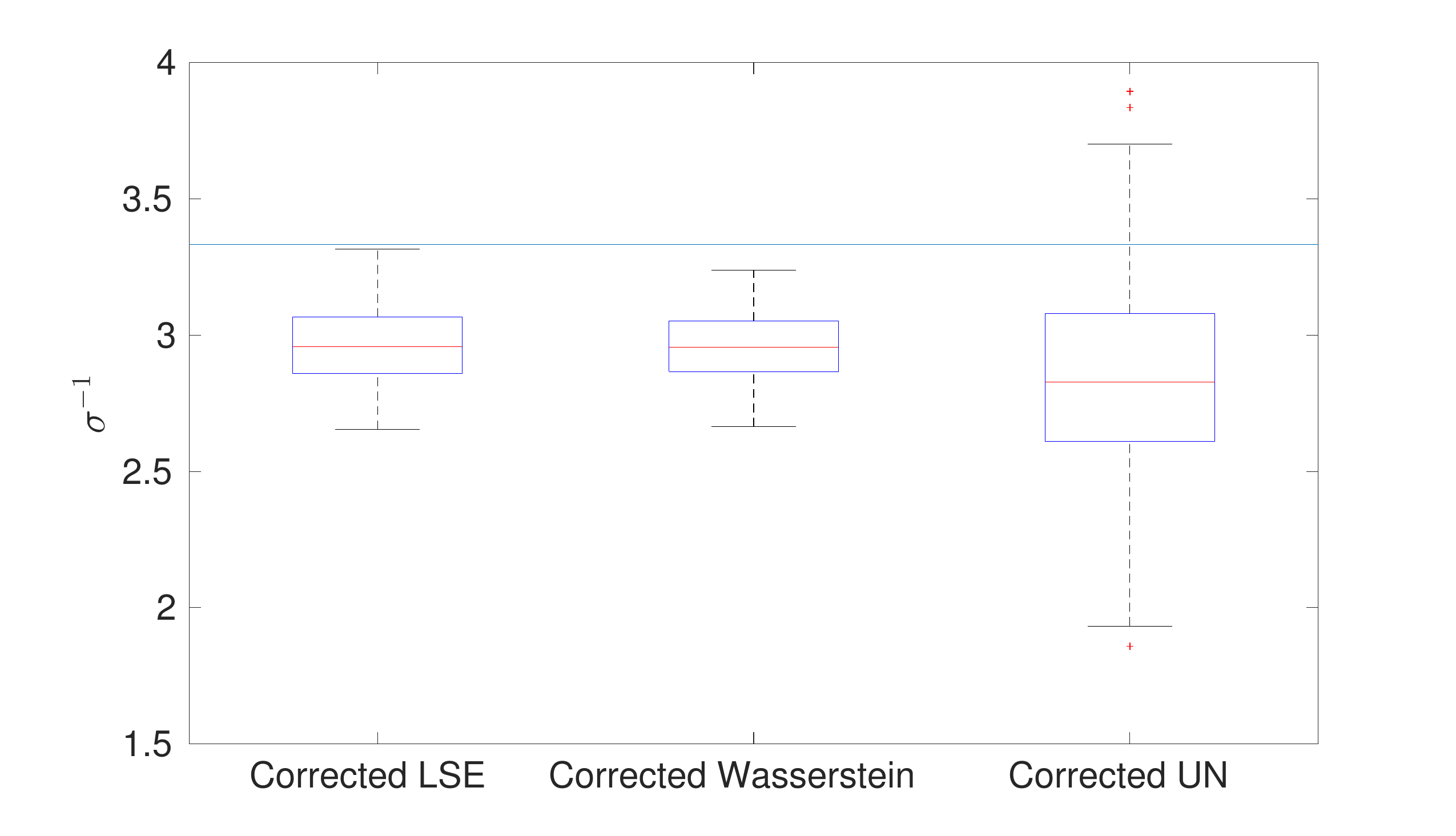}\includegraphics[width=3.5cm,height=3cm]{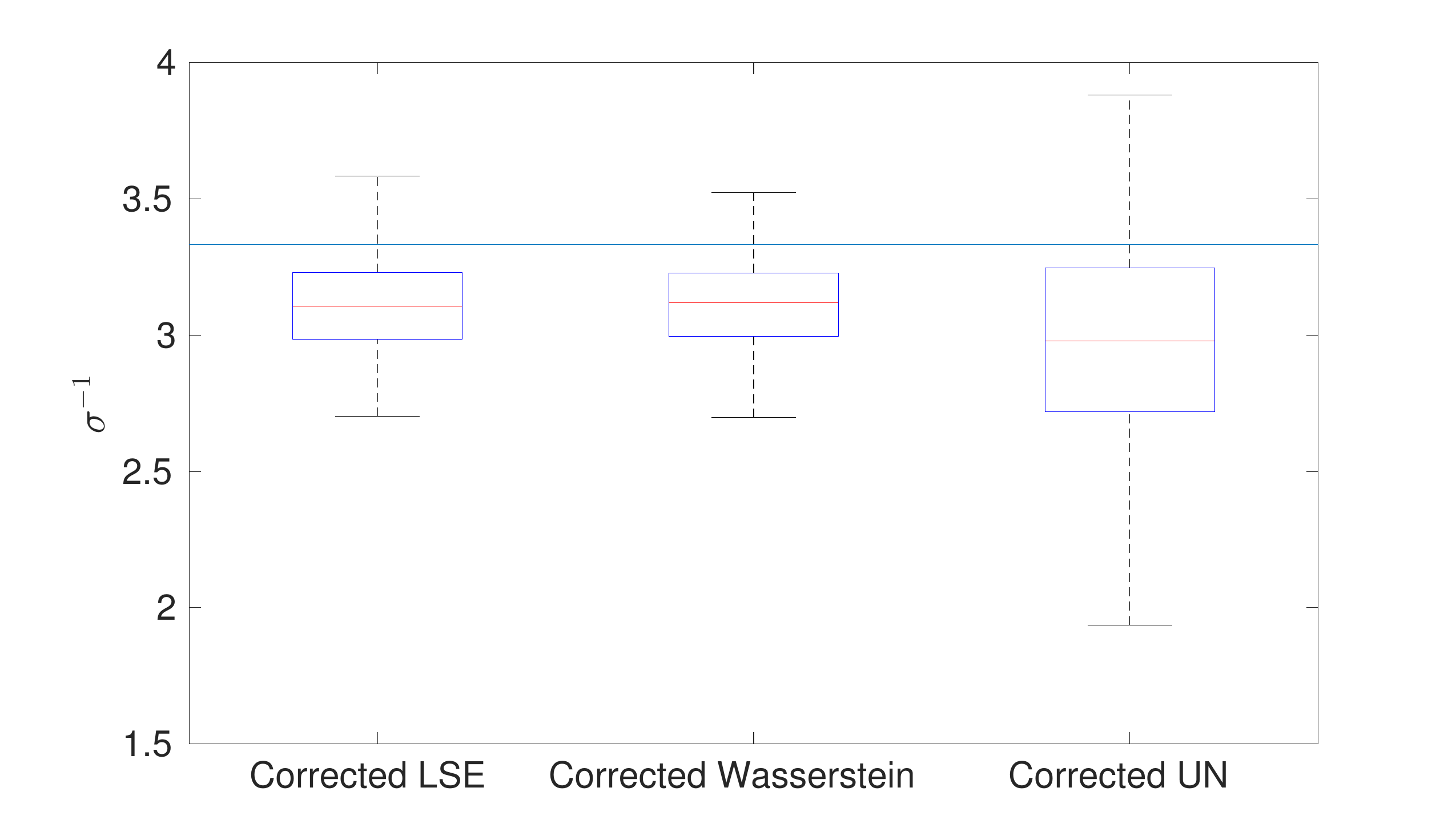}
	\caption{Influence of the size of the trees for $\sigma$ equals to $0.3$: tree sizes varying from $20$ nodes (left), $50$ nodes (center), to $100$ nodes (right). Forests of $50$ trees. Boxplots have been drawn from $100$ replicates each.}
	\label{fig:sizecompsmallsigma}
\end{figure}

For improved comparison, we also compute the estimator $\widehat{\lambda}_{un}[\mathcal{F}]$ of $\sigma^{-1}$ based on the work \cite{KPDRV14} (see Subsection \ref{ssec:adequacy}) given by
$$\widehat{\lambda}_{un}[\mathcal{F}] = \frac{1}{N}\sum_{i=1}^N \widehat{\delta}[\tau^i],$$
where $\widehat{\delta}[\tau^i]$ is defined (see equations \eqref{eq:delta1} and \eqref{eq:delta2}) from a node $v$ randomly chosen in $\tau^i$ by
$$\widehat{\delta}[\tau^i] = \frac{\sqrt{2}\,h(v)}{\sqrt{\pi \#\tau^i}}.$$
The estimator $\widehat{\lambda}_{un}[\mathcal{F}]$ is expected to present the bias due to the approximation of Harris paths by their expected limit. We correct it by the aforementioned method,
$$\widehat{\lambda}_{un}^c[\mathcal{F}] = \frac{1}{N}\sum_{i=1}^N \widehat{\eta}(\#\tau^i)\widehat{\delta}[\tau^i].$$
In Figures \ref{fig:biascorrection}, \ref{fig:sizecomp} and \ref{fig:sizecompsmallsigma}, estimators $\widehat{\lambda}_{ls}[\mathcal{F}]$, $\widehat{\lambda}_{W}[\mathcal{F}]$ and $\widehat{\lambda}_{un}[\mathcal{F}]$ are denoted by ``LSE'', ``Wasserstein'and ``Uniform node'(or ``UN'in short), respectively.

The study of Figure \ref{fig:biascorrection} shows that for values of $\sigma$ greater than $0.5$, the bias correction works properly. Moreover, it also shows that the approach developed in \cite{KPDRV14} presents the same kind of bias as ours. In the case of small parameter $\sigma$, the bias correction is not as accurate. This was expected because the bias corrector does not fit as well to the bias curve  for small values of sigma as it does for greater values of $\sigma$.

Since we have an estimation procedure which seems to work, the natural further study is to see how the quality of our estimators varies as the characteristics of the forest change. We begin by looking at the variations when the sizes of the trees increase. 
A priori, the sizes of the trees in the considered forest should not have influence on the variability of the estimators. Indeed, our estimation strategy is based on the approximation of the Harris path of a finite tree by its limit. As a consequence, the size parameter only governs the quality of this approximation. Whatever the sizes of the trees, the variability will be given by the variance of the limit distribution $\Lambda_{\infty}$. As expected Figure \ref{fig:sizecomp} shows that the variability of the estimators does not change as the sizes of the trees change when $\sigma$ takes great values. Similarly, as shown in Figure \ref{fig:sizecompsmallsigma}, for small values of $\sigma$, the sizes of the trees do not influence the dispersion of the estimator. However, Figure \ref{fig:sizecompsmallsigma} also shows that the sizes of the trees have a positive influence on the bias of the estimators.

\begin{figure}[!h]
	\centerfloat
	\includegraphics[scale=0.3]{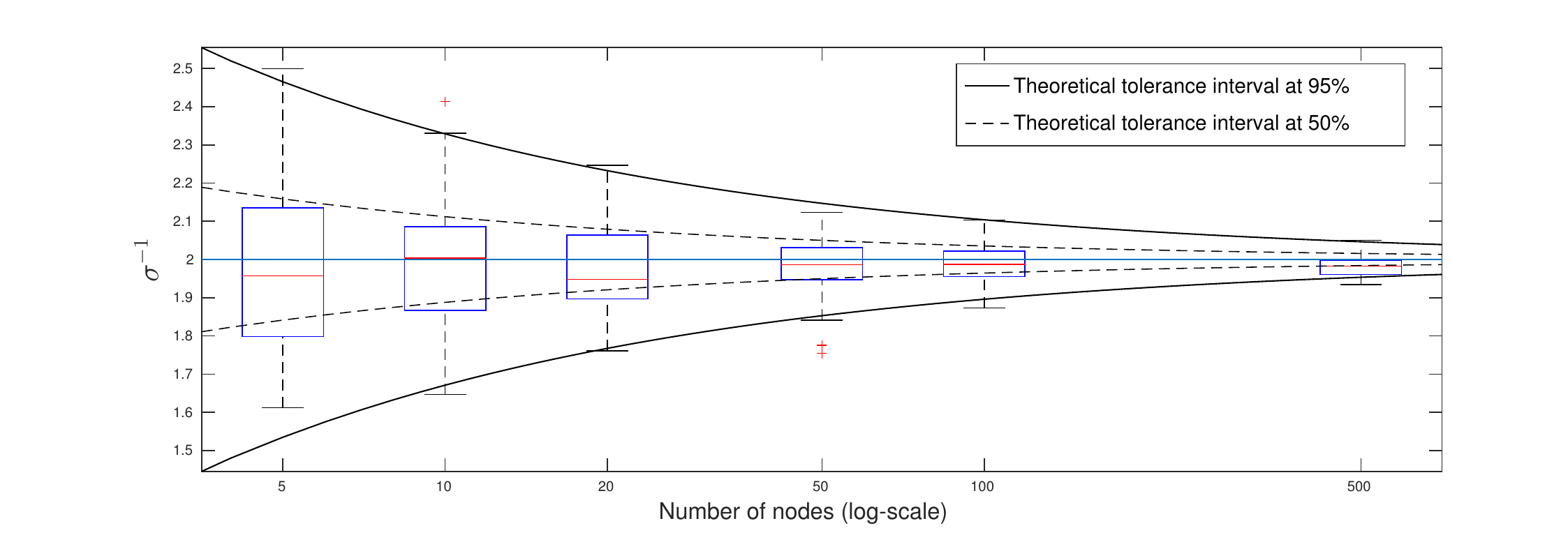}
	\caption{Least square estimation of $\sigma^{-1}$ for different sizes of forests ($\sigma=0.5$ with trees of size $20$). Boxplots have been drawn from $100$ replicates each.}
	\label{fig:Fsize}
\end{figure}

\begin{figure}[!h]
	\centerfloat
	\includegraphics[scale=0.3]{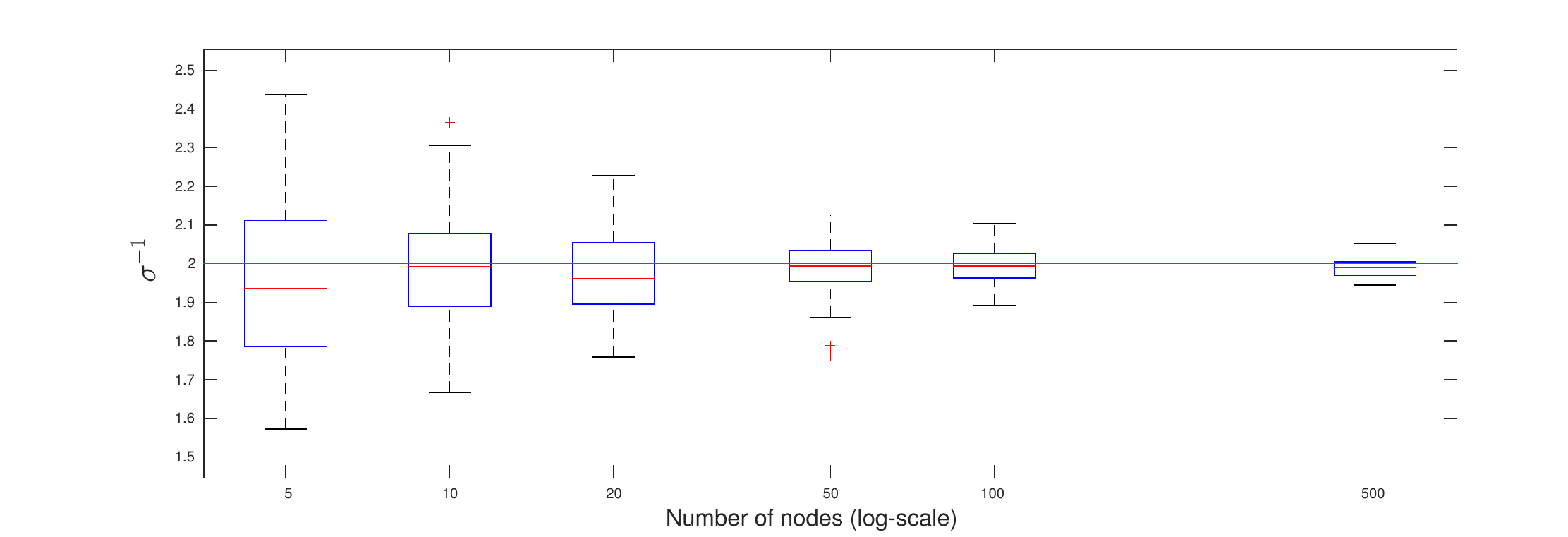}
	\caption{Wasserstein estimation of $\sigma^{-1}$ for different sizes of forests ($\sigma=0.5$ with trees of size $20$). Boxplots have been drawn from $100$ replicates each.}
	\label{fig:FsizeW}
\end{figure}

Finally, Figure \ref{fig:Fsize} shows the variation of the dispersion of the least square estimator as the size of the forest changes. It appears to be consistent with the theoretical tolerance intervals given by the central limit theorem. Similar results have been obtained from the Wasserstein method (see Figure \ref{fig:FsizeW}).

\subsection{Missing or noisy data}
\label{ss:missingnoisydata}

We focus here on the application of our statistical methods in the framework of missing or noisy data. This section is key for this paper because we exhibit difficult contexts in which Harris paths-based estimators perform well while the empirical variance is biased or can not be computed.

\subsubsection{Estimation with outliers}

We assume that the forest $\mathcal{F}=(\tau^i)_{1\leq i\leq N}$ is mainly composed of binary conditioned Galton-Watson trees with the same variance $\sigma^2$. However the forest also contains trees $\tau^i$ that have not been generated from the model and such that $\widehat{\lambda}[\tau^i]$ is significantly lesser or greater than the true parameter $\sigma^{-1}$. As a consequence, the estimators $\widehat{\lambda}_{ls}[\mathcal{F}]$ and $\widehat{\lambda}_{W}[\mathcal{F}]$ should be disturbed by these outliers. In the following simulation experiments, the forest contains $500$ conditioned Galton-Watson trees and $50$ outliers (addition of $10\%$ outliers).

\begin{figure}[!h]
	\centerfloat
	\includegraphics[scale=0.4,trim= 50mm 85mm 50mm 85mm]{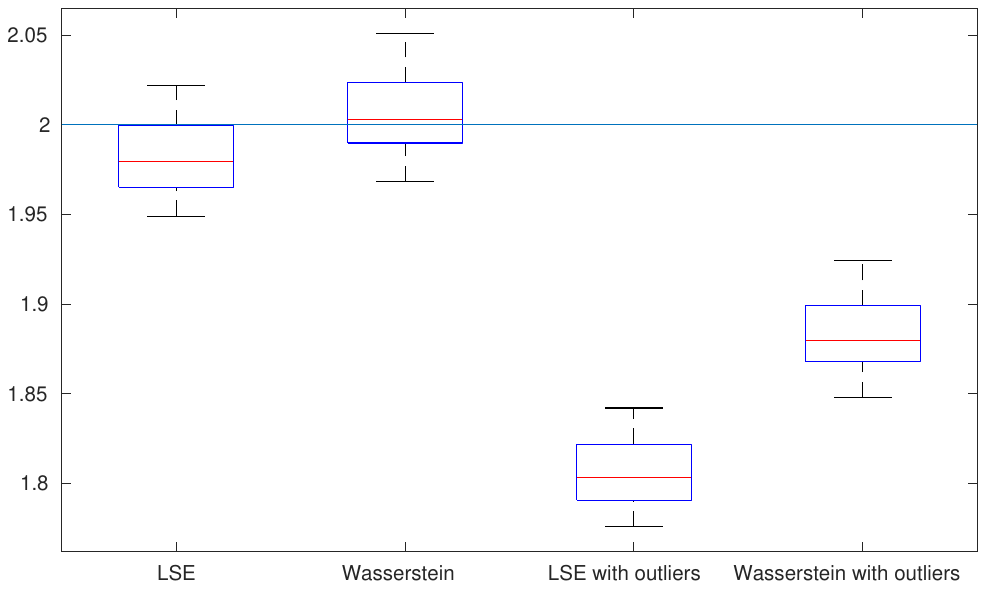}
	\caption{Boxplots of $\widehat{\lambda}_{ls}[\mathcal{F}]$ and $\widehat{\lambda}_{W}[\mathcal{F}]$ from forests containing (right) or not (left) small outliers impacting the quality of the estimation.}
	\label{fig:outliers:fav:wass}
\end{figure}

\begin{figure}[!h]
	\centerfloat
	\includegraphics[scale=0.4,trim= 50mm 85mm 50mm 85mm]{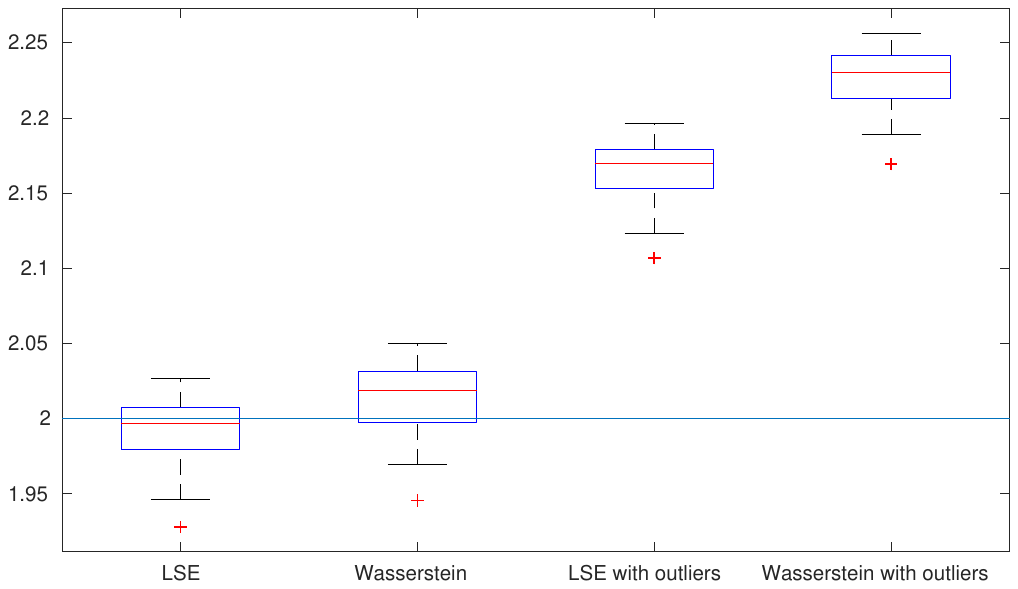}
	\caption{Boxplots of $\widehat{\lambda}_{ls}[\mathcal{F}]$ and $\widehat{\lambda}_{W}[\mathcal{F}]$ from forests containing (right) or not (left) large outliers impacting the quality of the estimation.}
	\label{fig:outliers:fav:lse}
\end{figure}

First, we consider the presence of outliers $\widehat{\lambda}[\tau^i]\simeq0.03$ smaller than the expected parameter $\sigma^{-1}=2$ and we compare our estimation strategies (see Figure \ref{fig:outliers:fav:wass}). One may observe that the Wasserstein method is less sensitive to these outliers than the least square estimator. This may be explained by the small weights $\int_{\frac{i-1}{N}}^{\frac{i}{N}}F_{\Lambda_\infty}^{-1}(s)\dd s$ given to the smallest values of $\widehat{\lambda}[\tau^{(i)}]$, that is to say, to the outliers, in formula \eqref{eq:wasserstein:formula} of $\widehat{\lambda}_{W}[\mathcal{F}]$. For the same reason, the Wasserstein estimator is more sensitive to large outliers $\widehat{\lambda}[\tau^i]\simeq2.9$ (see Figure \ref{fig:outliers:fav:lse}) than the least square strategy.

The numerical results of Subsection \ref{ss:simulateddata} show that the two strategies developed in this paper perform in a similar way on a dataset without outliers, which is clearly a benefit for the least square estimator $\widehat{\lambda}_{ls}[\mathcal{F}]$ easier to compute than $\widehat{\lambda}_{ls}[\mathcal{F}]$ (see Remark \ref{rem:Flambdamoins1}). Nevertheless, if one suspects the presence of small (large, respectively) outliers, these numerical experiments yield that the Wasserstein (least square, respectively) estimator should be privileged. In addition, it seems that the two estimators behave differently only under the presence of outliers. This observation could be used to detect suspicious data.

\subsubsection{Missing leaves}

\begin{figure}[!h]
	\centering
	\includegraphics[scale=0.25]{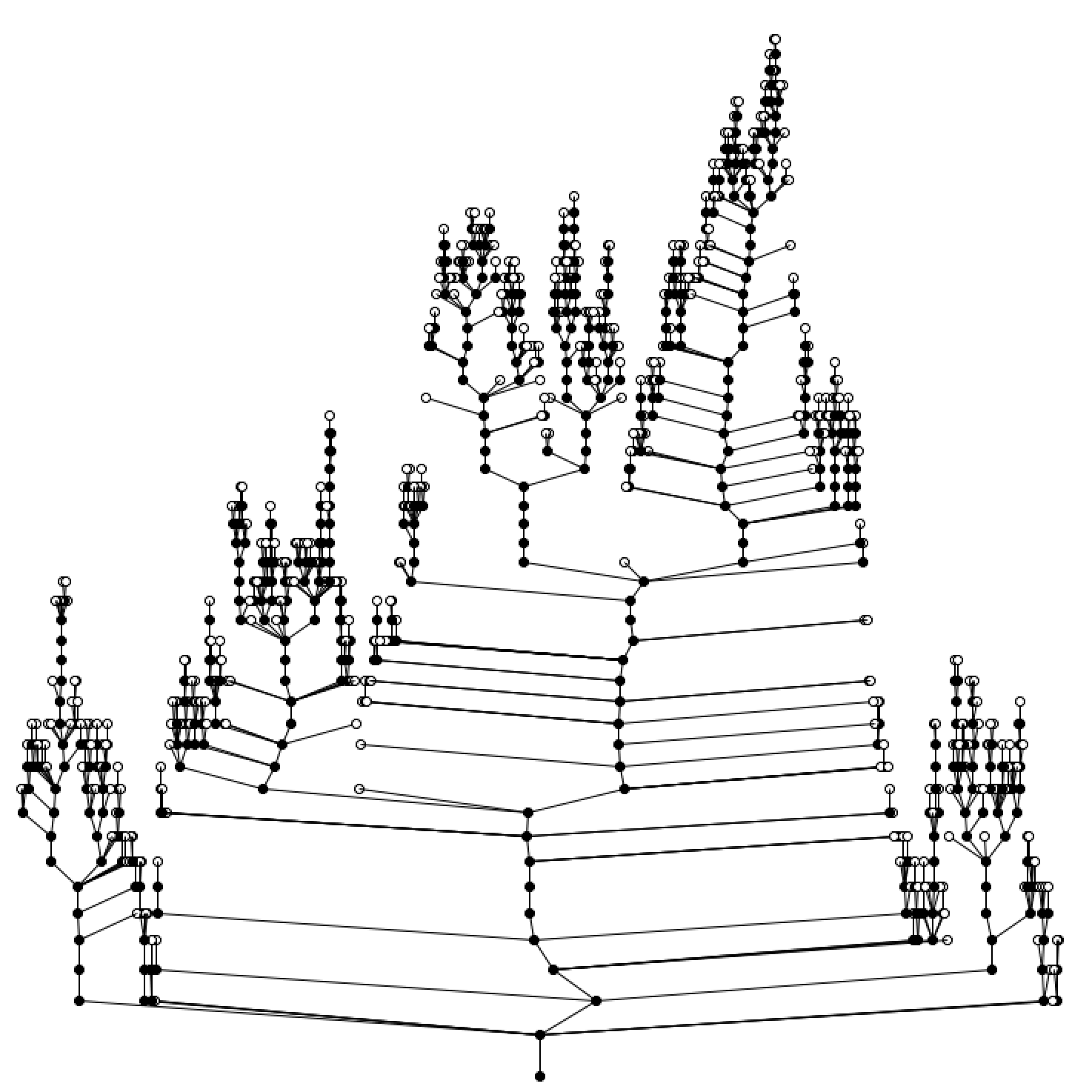}\\
		\centerfloat
	\includegraphics[scale=0.25,trim= 20mm 0mm 0mm 0mm]{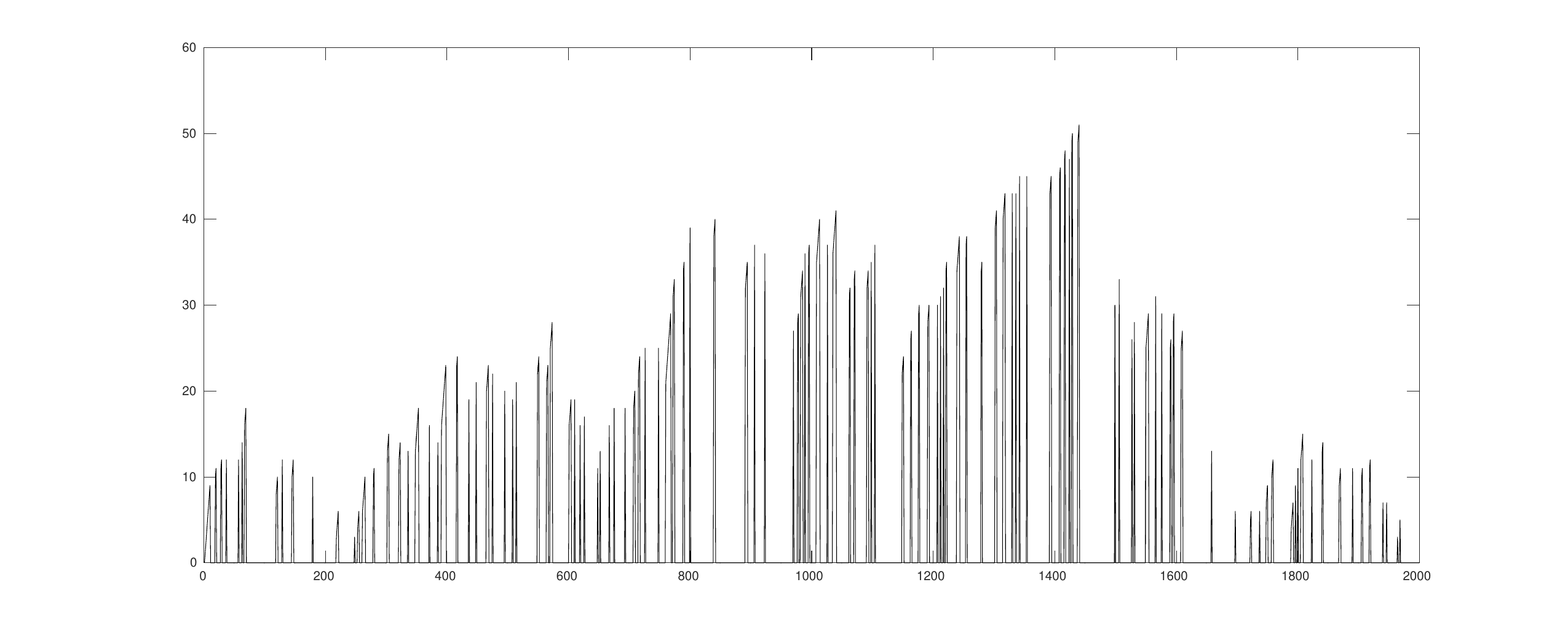}
	\caption{The conditioned Galton-Watson tree of Figure \ref{fig:exempleGW} (top) and its Harris path (bottom) are partially observed: leaves (in white) are missing.}
	\label{fig:feuillesblanches}
\end{figure}

\begin{figure}[!h]
	\centering
	\includegraphics[scale=0.4,trim=50mm 80mm 50mm 100mm]{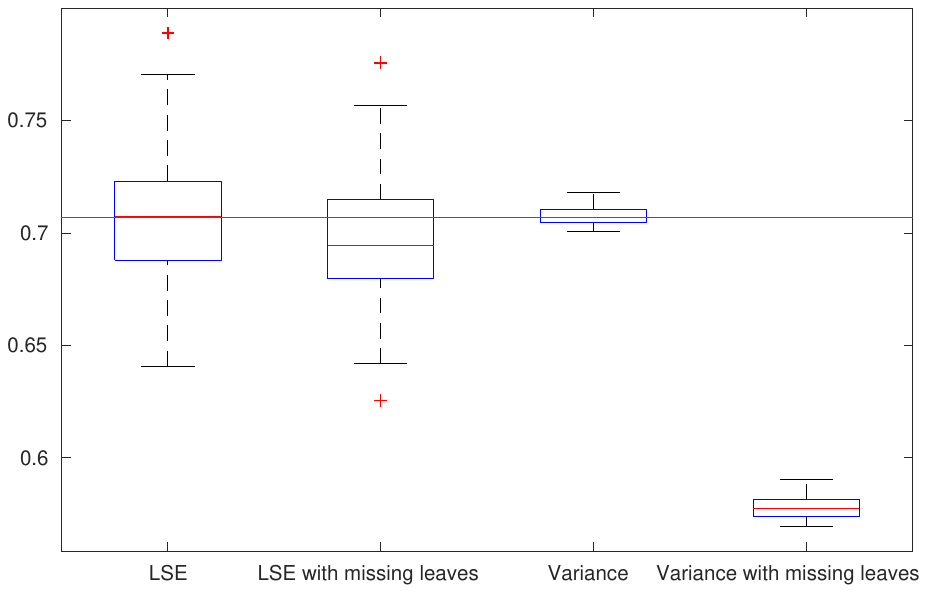}
	\caption{Least square and empirical estimators of $\sigma^{-1}$ from complete trees and trees with missing leaves.}
	\label{fig:missingleaves}
\end{figure}

We assume that we observe a conditioned critical Galton-Watson tree $\tau$ through a noise hiding its leaves. It means that, in the Harris path $\mathcal{H}[\tau]$, all the sections $[i-1,i+1]$ with
$$\mathcal{H}[\tau](i-1)=\mathcal{H}[\tau](i+1) = \mathcal{H}[\tau](i)-1,$$
i.e., sections corresponding to leaves, are also hidden. We refer the reader to Figure \ref{fig:feuillesblanches} for an example of conditioned Galton-Watson tree observed through this deleting noise. From a tree $\tau$ with $n$ nodes, we propose to estimate $\widehat{\lambda}[\tau]$ from the partially observed Harris path as follows.
Let us denote $\zeta\subset[0,2n]$ the union of the unobserved sections of the Harris path and $\eta=\frac{\zeta}{2n}\subset[0,1]$. Even if $\widehat{\lambda}[\tau]$ is uncomputable, we can approximate it by the solution $\widetilde{\lambda}[\tau]$ of the least square problem on $\eta^c$,
$$\lambda\mapsto\left\|\left(\frac{\mathcal{H}[\tau](2n\cdot)}{\sqrt{n}}- 2\lambda E\right)\mathbf{1}_{\eta^c}\right\|_2^2.$$
Mimicking the expression \eqref{eq:def:lambdahat} of $\widehat{\lambda}[\tau]$, $\widetilde{\lambda}[\tau]$ is given by
$$\widetilde{\lambda}[\tau] = \frac{\int_{[0,1]\cap\,\eta^c} \mathcal{H}[\tau](2nt) E_t\,\dd t}{2\sqrt{n}~\int_{[0,1]\cap\,\eta^c} E_t^2\,\dd t} .$$
From a forest of conditioned Galton-Watson trees $\tau$, we compute the values of $\widehat{\lambda}[\tau]$ and $\widetilde{\lambda}[\tau]$ in order to evaluate the influence of the noise on the quality of the estimation. We also compute the empirical variance $V[\tau]$ of the numbers of children appearing in $\tau$ from the complete tree and from its noisy version, $\sigma^{-1}$ being estimated by $V[\tau]^{-1/2}$. Some numerical results are presented in Figure \ref{fig:missingleaves}.


First, we remark that the distributions of $\widehat{\lambda}[\tau]$ and $\widetilde{\lambda}[\tau]$ are quite close, whereas the behavior of the empirical variance is highly disturbed by the absence of zeros in the set of numbers of children. As a consequence, statistical estimators computed from the Harris path seem to be more robust than empirical estimators of the birth distribution, even when the Harris path is largely hidden (see the example of Figure \ref{fig:feuillesblanches}).

\subsubsection{Partial observation of the Harris path}

\begin{figure}[!h]
	\centering
	\includegraphics[scale=0.18,trim= 20mm 0mm 0mm 0mm]{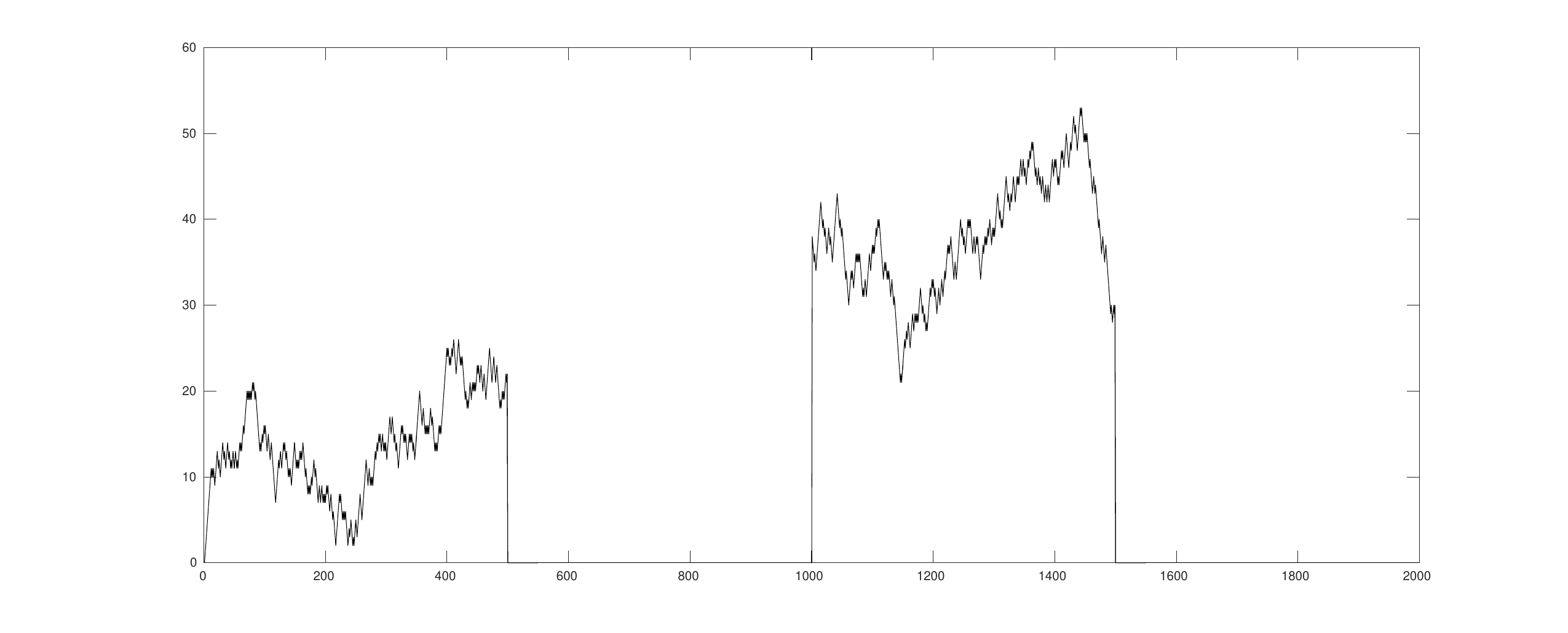}
	\includegraphics[scale=0.25,trim= 30mm 85mm 30mm 100mm]{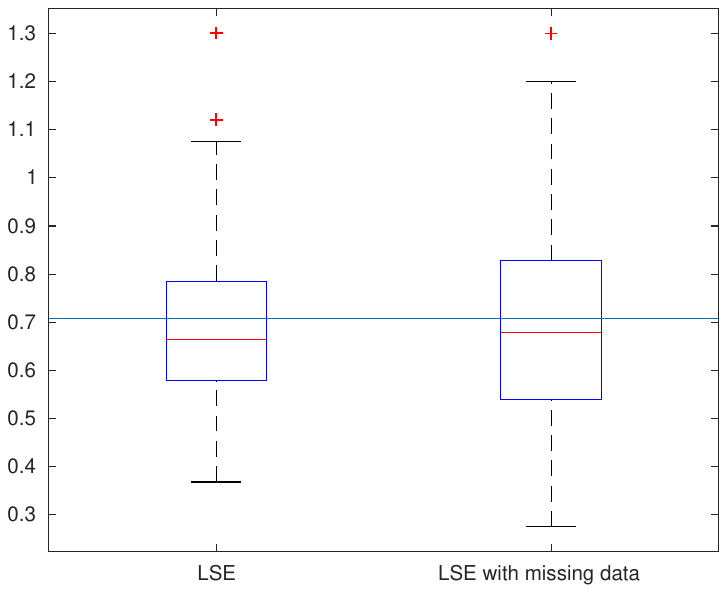}
	\caption{The Harris path of the conditioned Galton-Watson tree of Figure \ref{fig:exempleGW} observed only on the intervals $[0,500]$ and $[1000,1500]$ (left) and boxplots of $\widehat{\lambda}[\tau]$ from complete and partially observed Harris paths (right).}
	\label{fig:harrispath:trous}
\end{figure}

\begin{figure}[!h]
	\centering
	\includegraphics[scale=0.18,trim= 20mm 0mm 0mm 0mm]{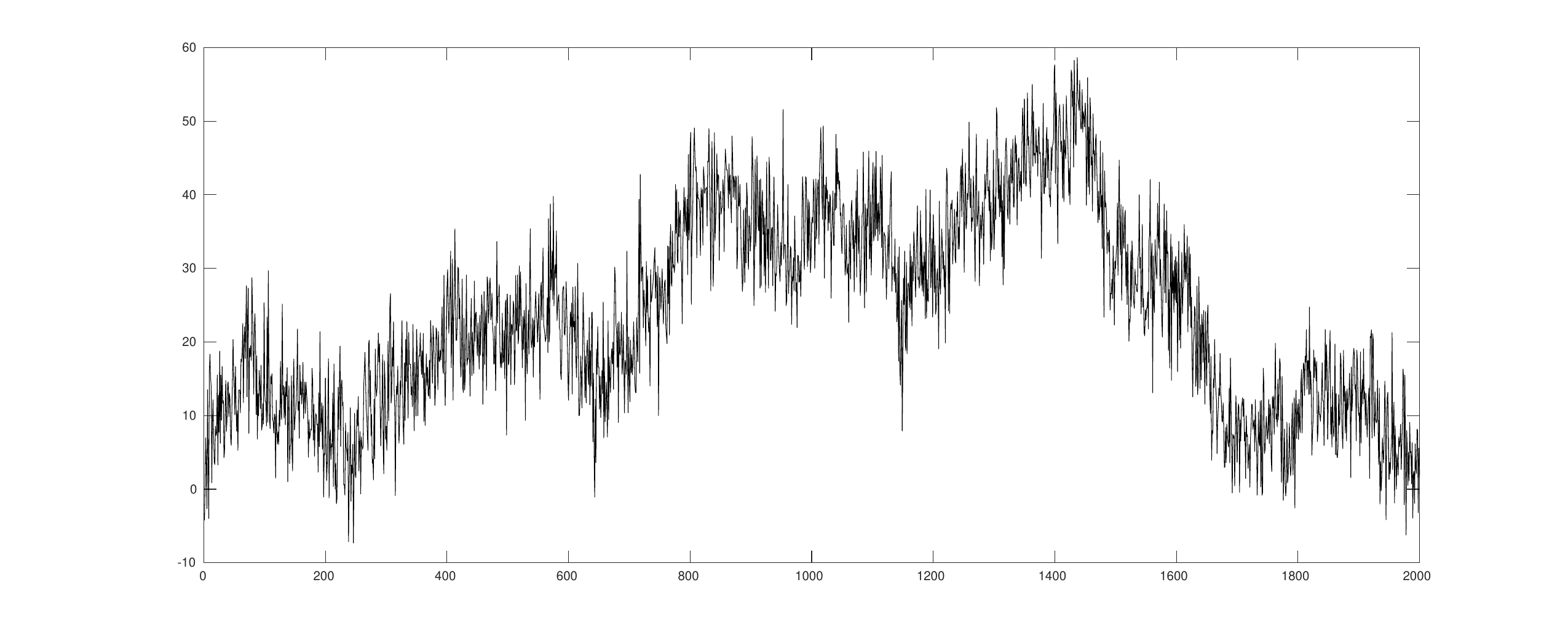}
	\includegraphics[scale=0.25,trim= 30mm 85mm 30mm 100mm]{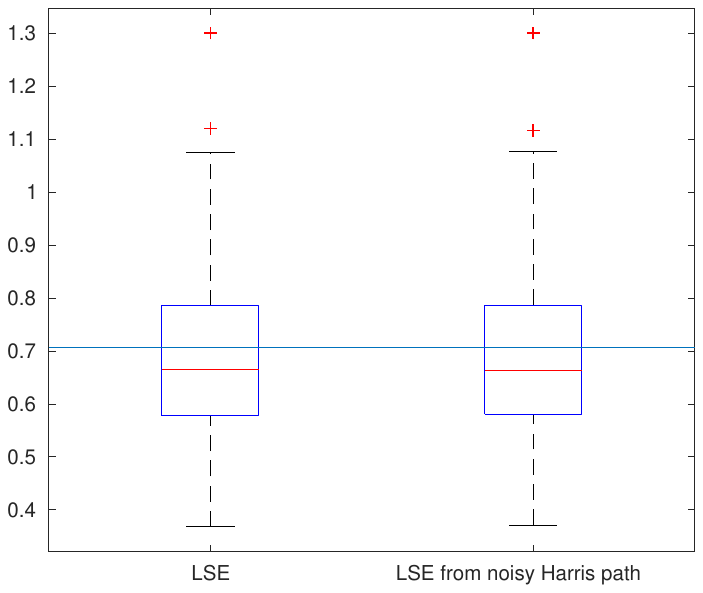}
	\caption{The Harris path of the conditioned Galton-Watson tree of Figure \ref{fig:exempleGW} observed through a Gaussian noise with standard deviation $5$ (left) and boxplots of $\widehat{\lambda}[\tau]$ from complete and partially observed Harris paths (right).}
	\label{fig:harrispath:gaussiannoise}
\end{figure}

Here we assume that the Harris path is partially observed in such a way that the underlying tree can not be reconstructed. This kind of disturbance may appear in data transmission where unwanted electromagnetic energy can degrade the quality of the signal. We consider two types of partial observation: (i) large sections of the Harris path are hidden (see Figure \ref{fig:harrispath:trous}) and (ii) the Harris path is observed through an additive Gaussian noise (see Figure \ref{fig:harrispath:gaussiannoise}).

\noindent
Since the tree can not be deduced from these noisy observations, empirical estimators of the birth distribution can not be computed, while statistical methods based on the Harris path are still feasible. The numerical results of Figures \ref{fig:harrispath:trous} and \ref{fig:harrispath:gaussiannoise} show that the distribution of $\widehat{\lambda}[\tau]$ is only slightly disturbed by the noise proving again the robustness of statistical estimators computed on Harris paths.


\section{Real data analysis: history of Wikipedia webpages}
\label{ss:realdata}

\begin{figure}[!h]
	\centering
	\subfloat{\includegraphics[scale=0.4, trim= 10mm 75mm 30mm 10mm]{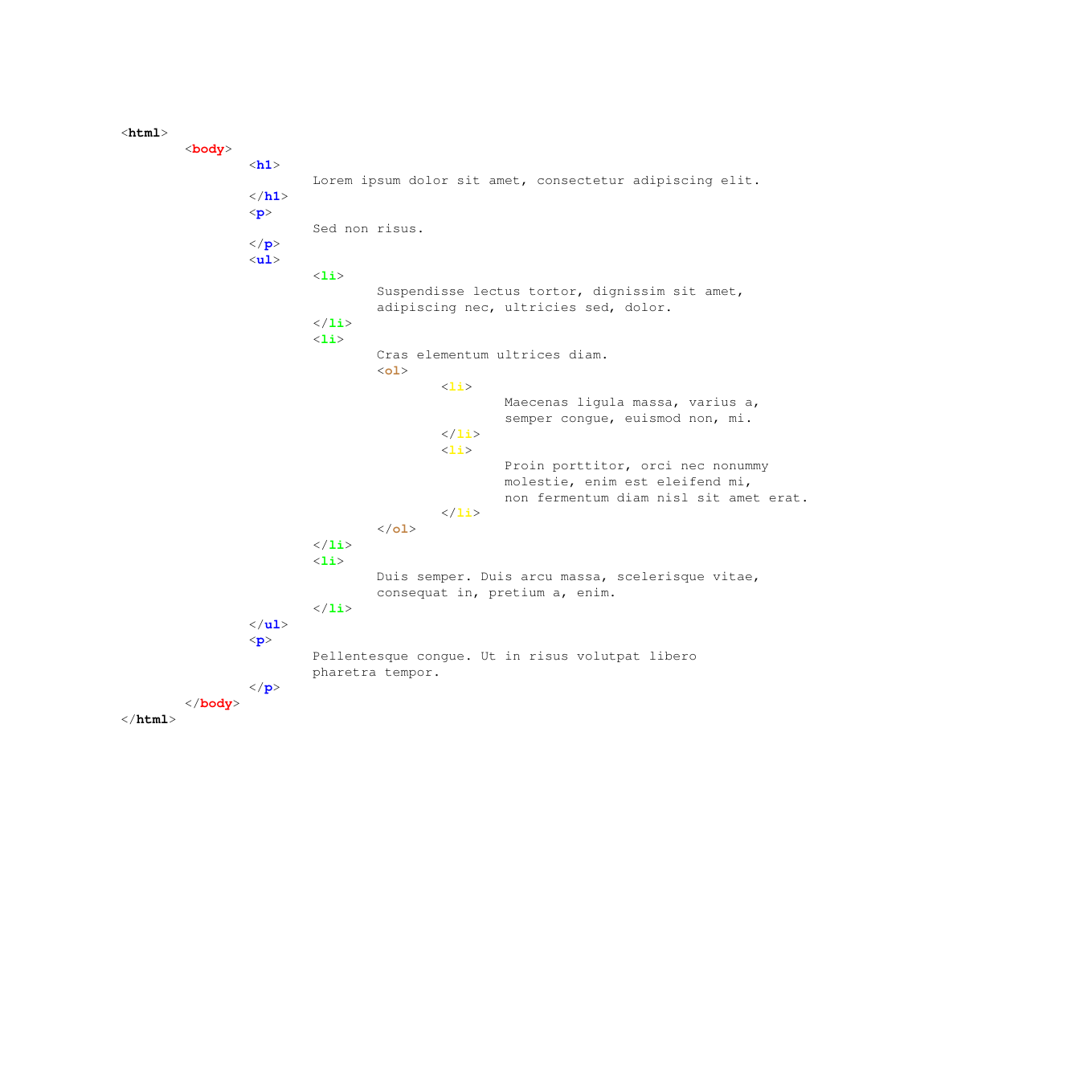}}
	\quad
	\subfloat{\begin{tikzpicture}[xscale=0.06,yscale=0.19]
		\tikzstyle{fleche}=[-,>=latex,thick]
		\tikzstyle{noeud}=[fill=black,circle,draw,scale=0.8]
		\tikzstyle{blanc}=[fill=white,circle,draw=none,scale=0.8]
		\tikzstyle{etiquette}=[midway]
		
		\def\DistanceInterNiveaux{5}
		\def\DistanceInterFeuilles{5}
		
		\def\NiveauA{(-5)*\DistanceInterNiveaux}
		\def\NiveauB{(-4)*\DistanceInterNiveaux}
		\def\NiveauC{(-3)*\DistanceInterNiveaux}
		\def\NiveauD{(-2)*\DistanceInterNiveaux}
		\def\NiveauE{(-1)*\DistanceInterNiveaux}
		\def\NiveauF{(0)*\DistanceInterNiveaux}
		\def\InterFeuilles{(1)*\DistanceInterFeuilles}
		
		\node[noeud] (A) at ({(0)*\InterFeuilles},{\NiveauA}) {};
		
		\node[noeud,color=red] (B) at ({(0)*\InterFeuilles},{\NiveauB}) {};
		\node[noeud,color=blue] (C1) at ({(-3)*\InterFeuilles},{\NiveauC}) {};
		\node[noeud,color=blue] (C2) at ({(-1)*\InterFeuilles},{\NiveauC}) {};
		\node[noeud,color=blue] (C3) at ({(1)*\InterFeuilles},{\NiveauC}) {};
		\node[noeud,color=blue] (C4) at ({(3)*\InterFeuilles},{\NiveauC}) {};
		
		\node[noeud,color=black!30!green] (D1) at ({(0)*\InterFeuilles},{\NiveauD}) {};
		\node[noeud,color=black!30!green] (D2) at ({(1)*\InterFeuilles},{\NiveauD}) {};
		\node[noeud,color=black!30!green] (D3) at ({(2)*\InterFeuilles},{\NiveauD}) {};
		
		\node[noeud,color=brown] (E) at ({(1)*\InterFeuilles},{\NiveauE}) {};
		
		\node[noeud,color=yellow] (F1) at ({(0.5)*\InterFeuilles},{\NiveauF}) {};
		\node[noeud,color=yellow] (F2) at ({(1.5)*\InterFeuilles},{\NiveauF}) {};
		
		\draw[fleche] (A)--(B) node[etiquette] {};
		\draw[fleche] (B)--(C1) node[etiquette] {};
		\draw[fleche] (B)--(C2) node[etiquette] {};
		\draw[fleche] (B)--(C3) node[etiquette] {};
		\draw[fleche] (B)--(C4) node[etiquette] {};
		\draw[fleche] (C3)--(D1) node[etiquette] {};
		\draw[fleche] (C3)--(D2) node[etiquette] {};
		\draw[fleche] (C3)--(D3) node[etiquette] {};
		\draw[fleche] (D2)--(E) node[etiquette] {};
		\draw[fleche] (E)--(F1) node[etiquette] {};
		\draw[fleche] (E)--(F2) node[etiquette] {};
		\end{tikzpicture}}
	\cprotect\caption{Underlying ordered tree structure (right) present in an \verb+HTML+ document (left). Each level in the tree is colored in the same way as the corresponding tags in the document. Natural order from top to bottom in the \verb|HTML| document corresponds to left-to-right order in the tree.}
	\label{fig:htmltree}
\end{figure}
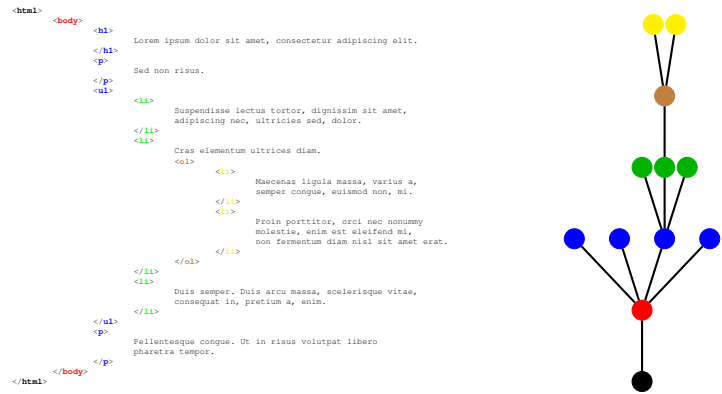

The aim of this section is to show that the methodology developed in this paper can be used to analyze the history of some real hierarchical data.
More precisely, we focus on the evolution over time of a given webpage on the World Wide Web.
\verb+HTML+ is the standard markup language for creating webpages.
Documents encoded in a markup language naturally presents a tree structure:
the area delimited by opening and closing tags represents a node of the tree;
the children of this node are given by the tags directly found in this area in the order they appear
(see Figure \ref{fig:htmltree} for an example of \verb+HTML+ document and the corresponding ordered tree structure).
It should be noted that the ordered tree representing an \verb+HTML+ document does not take into account the text between tags but only the hierarchical structure.

\lstset{
	keywords=[2]{body},
	keywords=[3]{html},
	keywordstyle=[2]\color{red}\bfseries,
	keywordstyle=[3]\color{black}\bfseries,
	keywords=[4]{h1,p,ul},
	keywords=[5]{li},
	keywordstyle=[4]\color{blue}\bfseries,
	keywordstyle=[5]\color{black!30!green}\bfseries,
	keywords=[6]{ol},
	keywordstyle=[6]\color{brown}\bfseries,
}
\lstset{aboveskip=0pt,belowskip=0pt}

\begin{figure}[!h]
	\centering
	\includegraphics[scale=0.25]{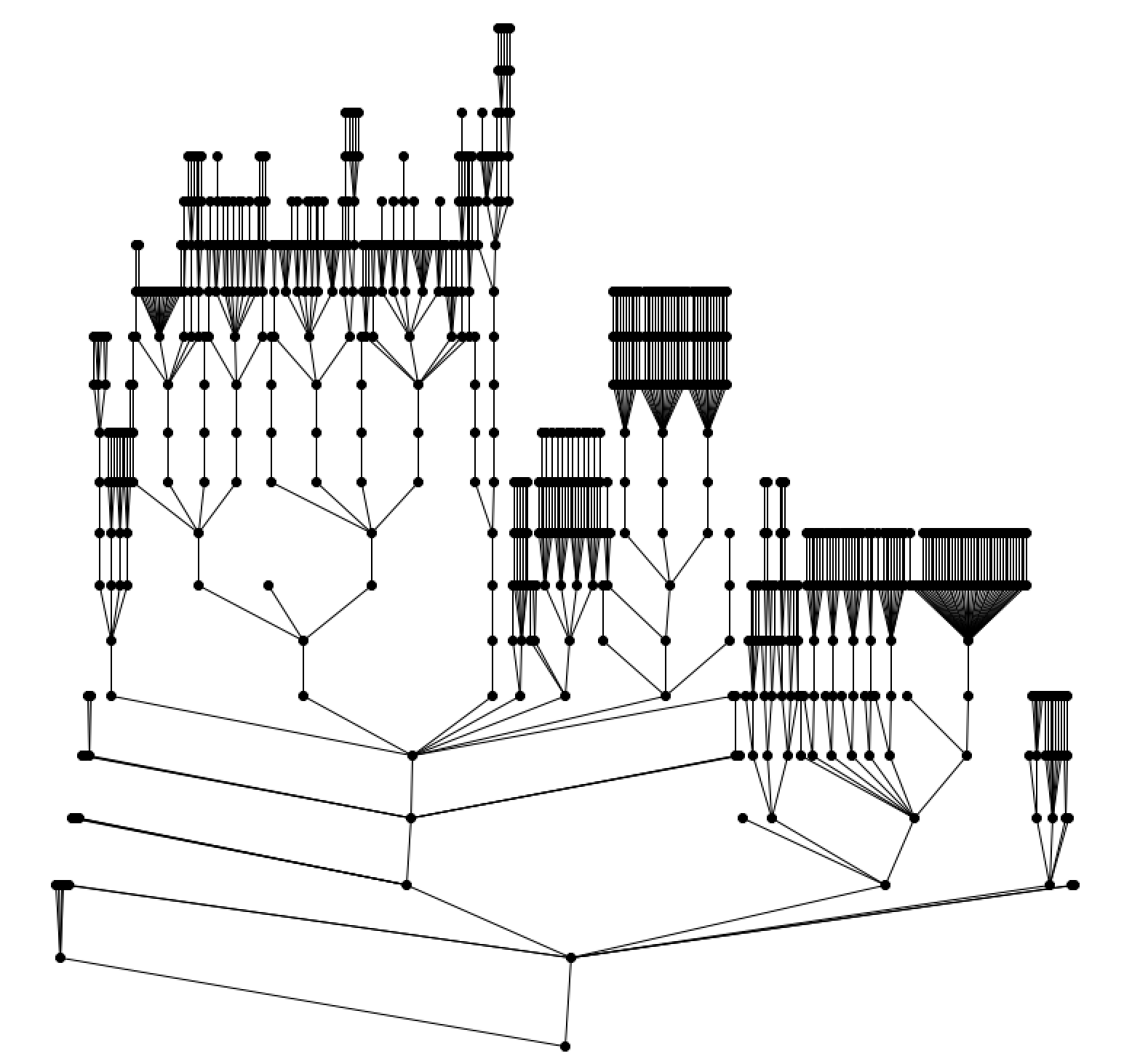}\\
		\centerfloat
	\includegraphics[scale=0.25,trim= 20mm 0mm 0mm 0mm]{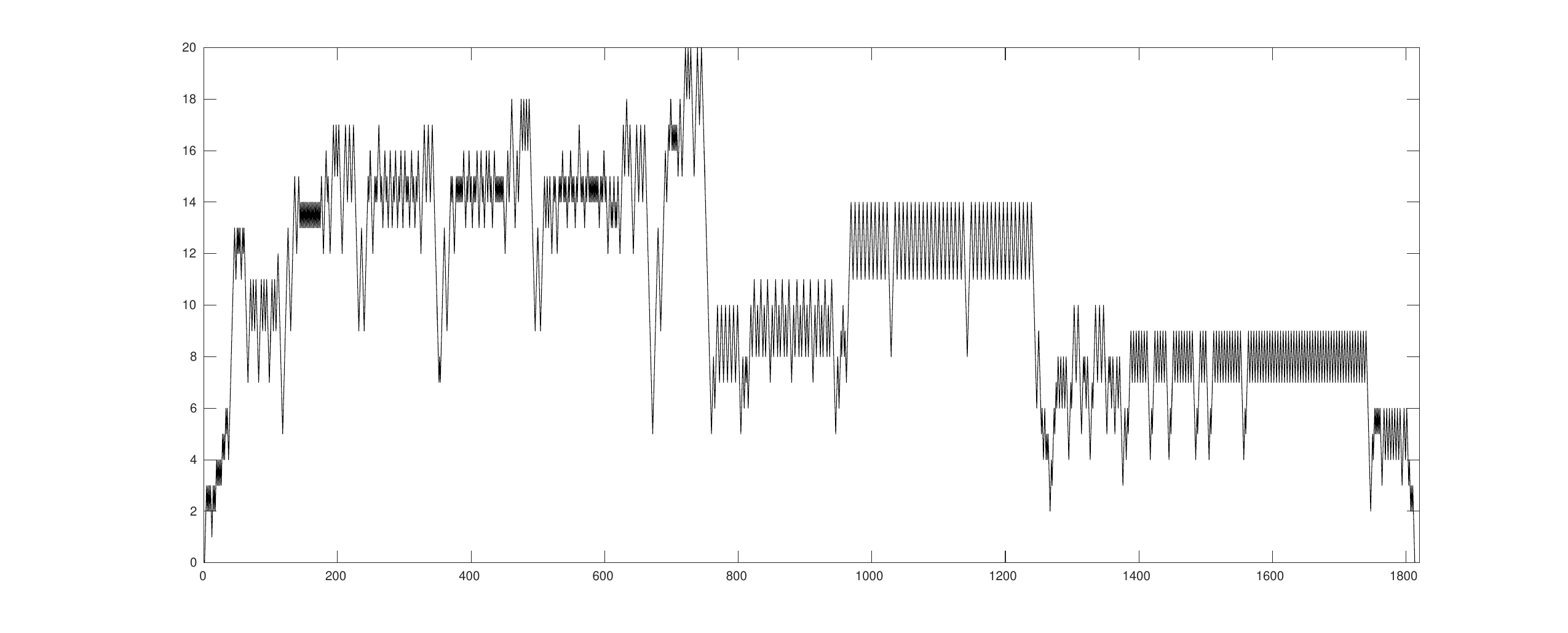}
	\caption{Underlying tree of the main page of Wikipedia accessed on April 12 2017 with 906 nodes (top) and its Harris path (bottom).}
	\label{ex:tree:wiki}
\end{figure}

Nowadays, Wikipedia is probably the most famous free Internet encyclopedia.
It allows its users to create and edit almost any article.
All past changes are listed in reverse-chronological order and are accessible from the current version of the Wikipedia webpage.
Consequently, each article forms a time series composed of hundreds of revisions.
The analysis of this chronological dataset is difficult because of the complex structure of the data which has no representation in an Euclidean state space.
We propose to apply the strategy presented in this paper to investigate this question and obtain informations on the history of articles.
Wikipedia webpages do not look like conditioned Galton-Watson trees (see Figure \ref{ex:tree:wiki} for a typical Wikipedia webpage and its Harris path which should be compared with the conditioned Galton-Watson tree with a comparable number of nodes of Figure \ref{fig:exempleGW}) but they share the same structure with a typical layout that consists in standardized \verb+HTML/CSS+ files on which articles are based.
Thus webpages at hand might not be differentiated by considering their shape but some scale parameter, as it is the case for conditioned Galton-Watson trees.
We claim that the quantity $\widehat{\lambda}[\tau]$, where $\tau$ is the underlying tree of a given webpage, is a good estimate of its relative scale and may be used to represent the revision history.

\begin{figure}[!h]
		\centerfloat
	\includegraphics[scale=0.25]{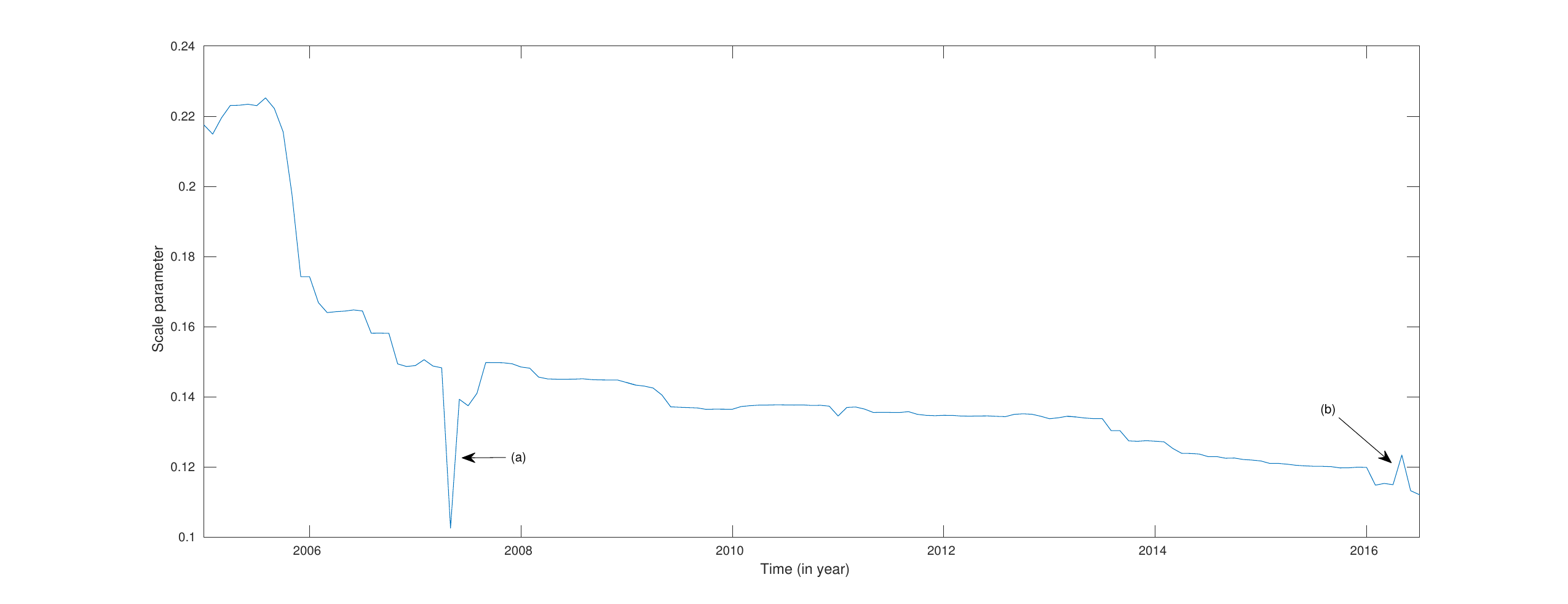}
	\caption{History of $\widehat{\lambda}_{ls}$ between January 2005 and June 2016 for the Wikipedia article \textit{Gravitational wave}. Events $(a)$ and $(b)$ are related to vandalism.}
	\label{fig:wiki:1}
\end{figure}

\begin{figure}[!h]
		\centerfloat
	\includegraphics[scale=0.2, trim= 20mm 0mm 160mm 0mm]{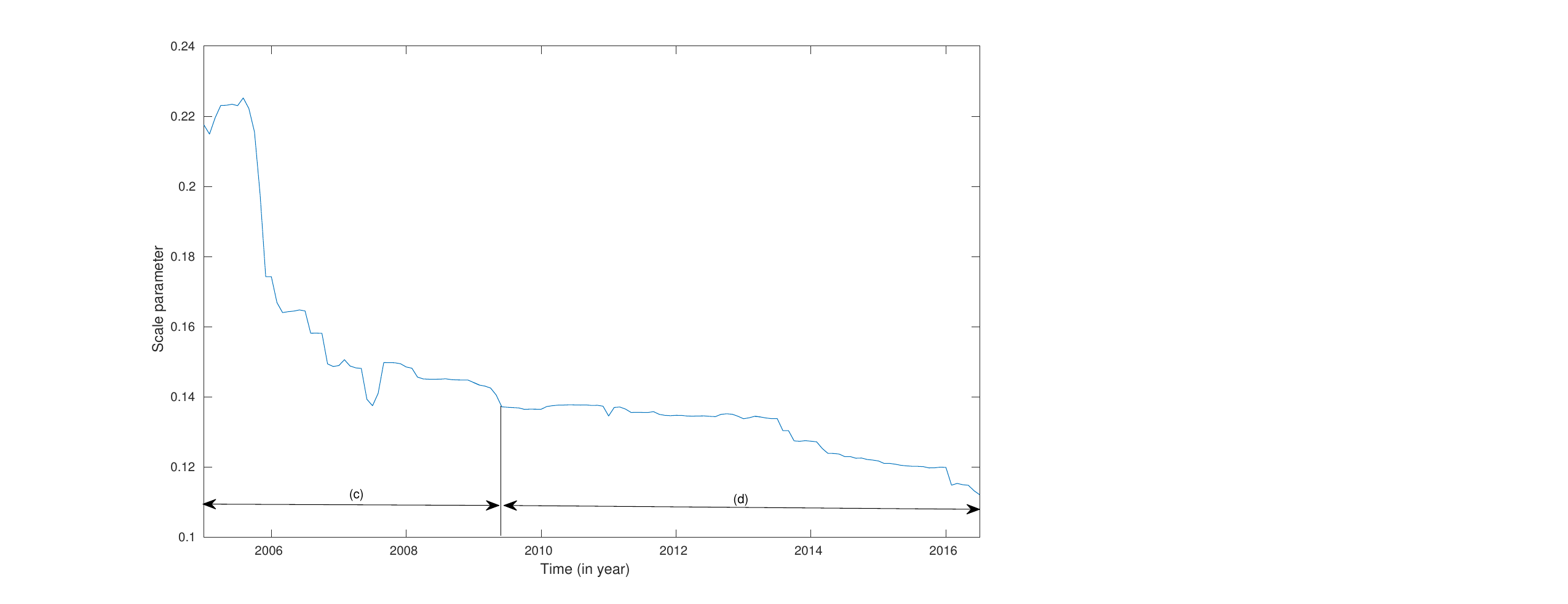}\includegraphics[scale=0.2, trim= 0mm 0mm 190mm 0mm]{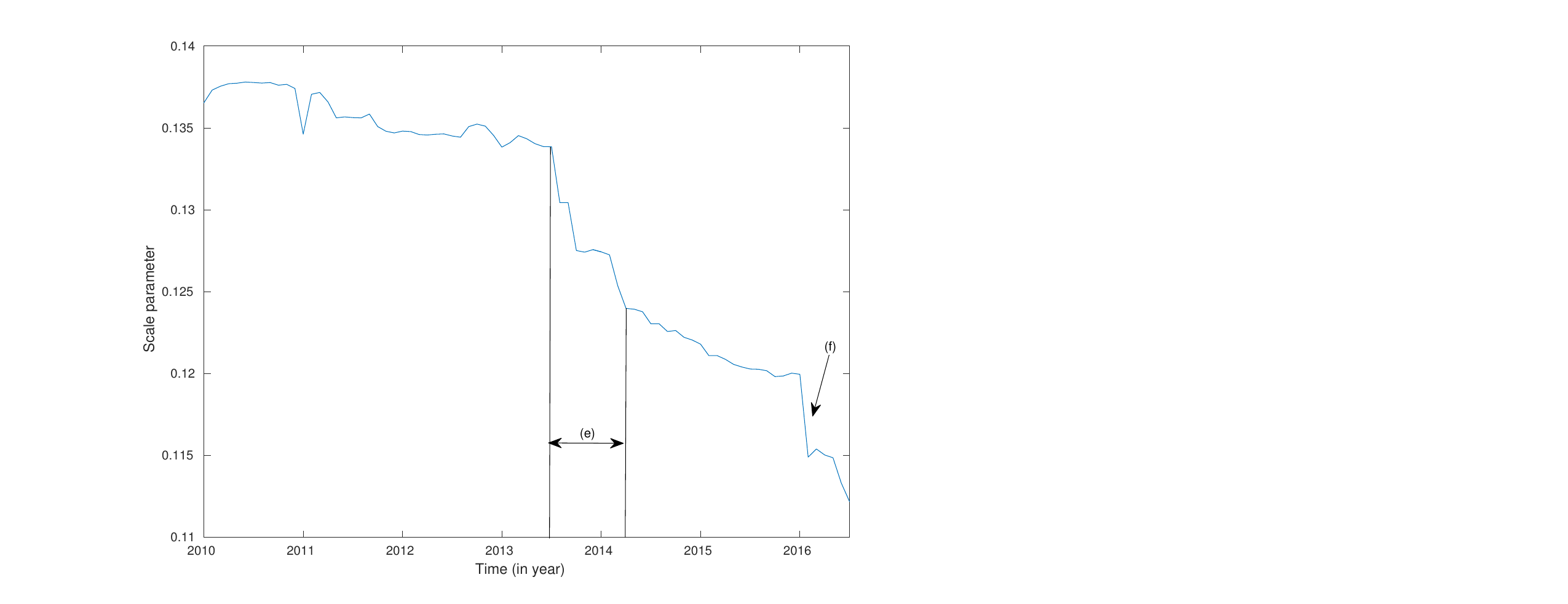}
	\caption{History of $\widehat{\lambda}_{ls}$ for the Wikipedia article \textit{Gravitational wave} without taking into account the two vandalism pages related to $(a)$ and $(b)$ between 2005 and 2016 (left) and 2010 and 2016 (right).}
	\label{fig:wiki:2}
\end{figure}

We begin with the English version of the Wikipedia article \textit{Gravitational wave}\footnote{Wikipedia article \textit{Gravitational wave}: \url{https://en.wikipedia.org/wiki/Gravitational_wave}}. This article has been edited 2001 times by 810 Wikipedians since its creation on September 3rd 2001 (information acquired on August 11 2016).
For each month since January 2005, we compute $\widehat{\lambda}_{ls}$ from the forest of the versions revised during this month.
If no revision has been found during this period, $\widehat{\lambda}_{ls}$ is equal to the estimate of the previous month, and recursively.
Figure \ref{fig:wiki:1} displays the evolution of $\widehat{\lambda}_{ls}$ over time.
First, we remark two spikes $(a)$, negative in May 2007, and $(b)$, positive in May 2016.
Both these spikes correspond to massive vandalism of the article on May 7 2007 (addition of 720 pointless sections with random text) and May 23 2016 (complete deletion of the article) by malicious people.
Indeed, if we do not consider these two vandalized webpages in our estimation, we obtain the graph of Figure \ref{fig:wiki:2} (left) that has no spikes.
In addition, we observe in Figure \ref{fig:wiki:2} (left) that the time series of $\widehat{\lambda}_{ls}$ has roughly two regimes $(c)$ and $(d)$.
The first period $(c)$ corresponds to the ``running in'required to find the adequate structure of the article.
In this period, the webpage is subject to major changes that are most often additions of new sections or paragraphs but may be deletions of inappropriate content.
When a good structure arises, the webpage is then slowly broadened during the second regime $(d)$.
It should be remarked in Figure \ref{fig:wiki:2} (right) that two important modifications occur in the period $(d)$: $(e)$ between July 2013 and April 2014 and $(f)$ in February 2016.
The $(e)$ period is related to major changes in the webpage (mainly addition of references and reorganization of some sections) especially following advances in this field.
The second event $(f)$ corresponds to extensive adding following the announce of the first observation of gravitational waves using the Advance LIGO detectors.

\begin{figure}[!h]
	\centerfloat
	\includegraphics[scale=0.3]{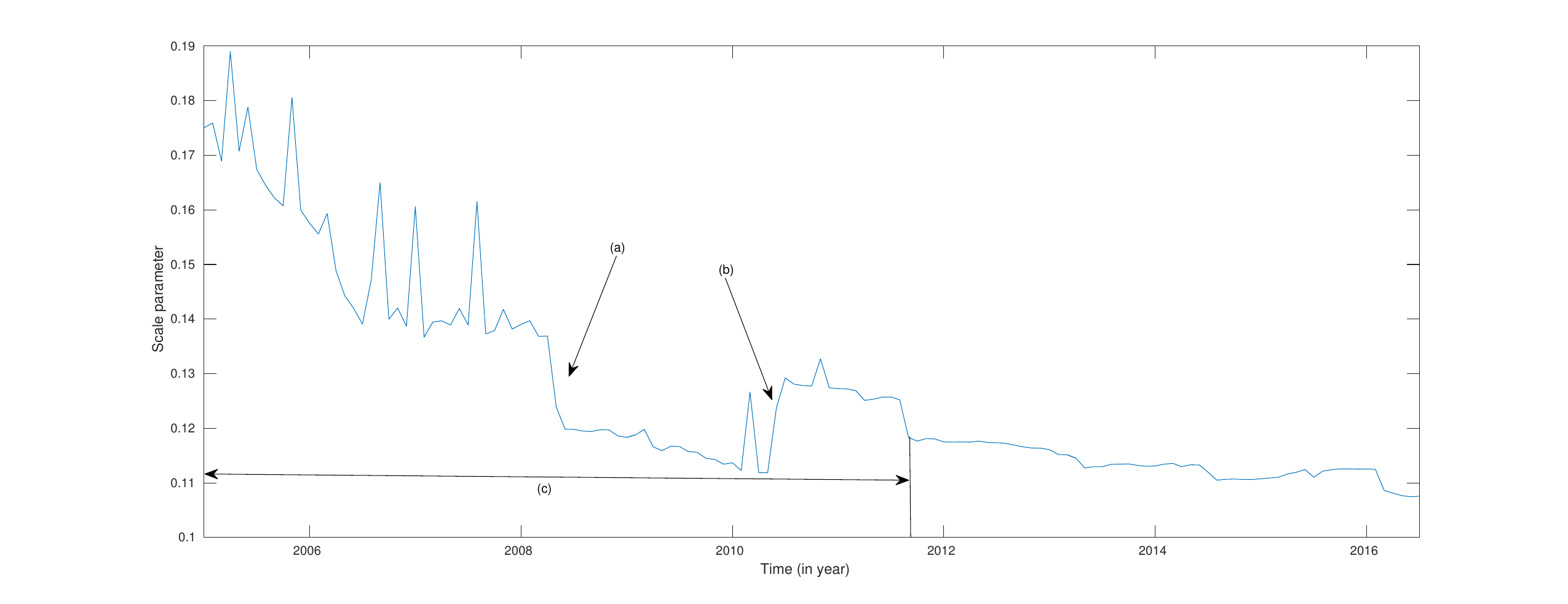}
	\caption{History of $\widehat{\lambda}_{ls}$ between January 2005 and June 2016 for the Wikipedia article \textit{Chocolate}. All spikes are related to vandalism.}
	\label{fig:wiki:3}
\end{figure}

We perform the same methodology on the history of the Wikipedia article \textit{Chocolate}\footnote{Wikipedia article \textit{Chocolate}: \url{https://en.wikipedia.org/wiki/Chocolate}} (see Figure \ref{fig:wiki:3}).
This article has been edited 6332 times by 3105 Wikipedians since its creation on
November 13 2001 (information acquired on August 11 2016).
All the spikes observed on the graph of Figure \ref{fig:wiki:3} correspond to acts of vandalism (deletion of substantial content). For the sake of example, we highlight two
major events $(a)$ (in May 2008) and $(b)$ (in June 2010) occuring during the ``running in'period $(c)$: $(a)$ corresponds to important additions in the article (sections \textit{Etymology}, \textit{Holydays} and \textit{Manufacturers} have been added), while $(b)$ is related to the creation of the parallel article {\it Health effects of chocolate} leading to deletion of the corresponding sections in the main article.

For both examples, we empirically observe that, when $\widehat{\lambda}_{ls}$ decreases, some content has been added to the webpage, and conversely, when $\widehat{\lambda}_{ls}$ increases, some parts of the article have been removed.
Our analysis shows that, starting from their creation, these Wikipedia articles are broadened over time after a long ``running in'period used to unconsciously find the
adequate structure. One may also detect vandalism on Wikipedia articles by identifying spikes a posteriori. 
Vandalism is usually removed by dedicated individuals who patrol Wikipedia webpages, but this is an onerous task with a rate of $10$ edits per second\footnote{Wikipedia statistics (last consulted on August 11 2016): \url{https://en.wikipedia.org/wiki/Wikipedia:Statistics}} and around $7\%$ of edits have been estimated to be vandalism \cite{Potthast:2010:CWV:1835449.1835617}. Vandalism detection is often based on a combination of various indicating features \cite{Adler2011,Mola2010}. Our algorithm might be used as a new feature for identifying acts of vandalism on the structure of the article.

\section*{Acknowledgment}
The authors would like to show their gratitude to an anonymous reviewer who provided many relevant comments on the manuscript.

\bibliographystyle{plain}
\bibliography{main}

\end{document}